
\ifx\shlhetal\undefinedcontrolsequence\let\shlhetal\relax\fi
\def\fmtname{AmS-TeX}

\def\fmtversion{2.1}
\catcode`\@=11
\ifx\amstexloaded@\relax\catcode`\@=\active
  \endinput\else\let\amstexloaded@\relax\fi
\newlinechar=`\^^J
\def\W@{\immediate\write\sixt@@n}
\def\CR@{\W@{^^J\fmtname - Version \fmtversion^^J%
COPYRIGHT 1985, 1990, 1991 - AMERICAN MATHEMATICAL SOCIETY^^J%
Use of this macro package is not restricted provided^^J%
each use is acknowledged upon publication.^^J}}
\CR@ \everyjob{\CR@}
\message{Loading definitions for}
\message{misc utility macros,}
\toksdef\toks@@=2
\long\def\rightappend@#1\to#2{\toks@{\\{#1}}\toks@@
 =\expandafter{#2}\xdef#2{\the\toks@@\the\toks@}\toks@{}\toks@@{}}
\def\alloclist@{}
\newif\ifalloc@
\def\showallocations{{\def\\{\immediate\write\m@ne}\alloclist@}\alloc@true}
\def\alloc@#1#2#3#4#5{\global\advance\count1#1by\@ne
 \ch@ck#1#4#2\allocationnumber=\count1#1
 \global#3#5=\allocationnumber
 \edef\next@{\string#5=\string#2\the\allocationnumber}%
 \expandafter\rightappend@\next@\to\alloclist@}
\newcount\count@@
\newcount\count@@@
\def\FN@{\futurelet\next}
\def\DN@{\def\next@}
\def\DNii@{\def\nextii@}
\def\RIfM@{\relax\ifmmode}
\def\RIfMIfI@{\relax\ifmmode\ifinner}
\def\setboxz@h{\setbox\z@\hbox}
\def\wdz@{\wd\z@}
\def\boxz@{\box\z@}
\def\setbox@ne{\setbox\@ne}
\def\wd@ne{\wd\@ne}
\def\iterate{\body\expandafter\iterate\else\fi}
\def\err@#1{\errmessage{AmS-TeX error: #1}}
\newhelp\defaulthelp@{Sorry, I already gave what help I could...^^J
Maybe you should try asking a human?^^J
An error might have occurred before I noticed any problems.^^J
``If all else fails, read the instructions.''}
\def\Err@{\errhelp\defaulthelp@\err@}
\def\eat@#1{}
\def\in@#1#2{\def\in@@##1#1##2##3\in@@{\ifx\in@##2\in@false\else\in@true\fi}%
 \in@@#2#1\in@\in@@}
\newif\ifin@
\def\space@.{\futurelet\space@\relax}
\space@. %
\newhelp\athelp@
{Only certain combinations beginning with @ make sense to me.^^J
Perhaps you wanted \string\@\space for a printed @?^^J
I've ignored the character or group after @.}
{\catcode`\~=\active 
 \lccode`\~=`\@ \lowercase{\gdef~{\FN@\at@}}}
\def\at@{\let\next@\at@@
 \ifcat\noexpand\next a\else\ifcat\noexpand\next0\else
 \ifcat\noexpand\next\relax\else
   \let\next\at@@@\fi\fi\fi
 \next@}
\def\at@@#1{\expandafter
 \ifx\csname\space @\string#1\endcsname\relax
  \expandafter\at@@@ \else
  \csname\space @\string#1\expandafter\endcsname\fi}
\def\at@@@#1{\errhelp\athelp@ \err@{\Invalid@@ @}}
\def\atdef@#1{\expandafter\def\csname\space @\string#1\endcsname}
\newhelp\defahelp@{If you typed \string\define\space cs instead of
\string\define\string\cs\space^^J
I've substituted an inaccessible control sequence so that your^^J
definition will be completed without mixing me up too badly.^^J
If you typed \string\define{\string\cs} the inaccessible control sequence^^J
was defined to be \string\cs, and the rest of your^^J
definition appears as input.}
\newhelp\defbhelp@{I've ignored your definition, because it might^^J
conflict with other uses that are important to me.}
\def\define{\FN@\define@}
\def\define@{\ifcat\noexpand\next\relax
 \expandafter\define@@\else\errhelp\defahelp@                               
 \err@{\string\define\space must be followed by a control
 sequence}\expandafter\def\expandafter\nextii@\fi}                          
\def\undefined@@@@@@@@@@{}
\def\preloaded@@@@@@@@@@{}
\def\next@@@@@@@@@@{}
\def\define@@#1{\ifx#1\relax\errhelp\defbhelp@                              
 \err@{\string#1\space is already defined}\DN@{\DNii@}\else
 \expandafter\ifx\csname\expandafter\eat@\string                            
 #1@@@@@@@@@@\endcsname\undefined@@@@@@@@@@\errhelp\defbhelp@
 \err@{\string#1\space can't be defined}\DN@{\DNii@}\else
 \expandafter\ifx\csname\expandafter\eat@\string#1\endcsname\relax          
 \global\let#1\undefined\DN@{\def#1}\else\errhelp\defbhelp@
 \err@{\string#1\space is already defined}\DN@{\DNii@}\fi
 \fi\fi\next@}

\def\predefine#1#2{\let#1#2}
\def\undefine#1{\let#1\undefined}
\message{page layout,}
\newdimen\captionwidth@
\captionwidth@\hsize
\advance\captionwidth@-1.5in
\def\pagewidth#1{\hsize#1\relax
 \captionwidth@\hsize\advance\captionwidth@-1.5in}
\def\pageheight#1{\vsize#1\relax}
\def\hcorrection#1{\advance\hoffset#1\relax}
\def\vcorrection#1{\advance\voffset#1\relax}
\message{accents/punctuation,}

\let\graveaccent\`
\let\acuteaccent\'
\let\tildeaccent\~
\let\hataccent\^
\let\underscore\_
\let\B\=
\let\D\.
\let\ic@\/
\def\/{\unskip\ic@}
\def\textfonti{\the\textfont\@ne}
\def\t#1#2{{\edef\next@{\the\font}\textfonti\accent"7F \next@#1#2}}
\def~{\unskip\nobreak\ \ignorespaces}
\def\.{.\spacefactor\@m}
\atdef@;{\leavevmode\null;}
\atdef@:{\leavevmode\null:}
\atdef@?{\leavevmode\null?}
\edef\@{\string @}
\def\relaxnext@{\let\next\relax}
\atdef@-{\relaxnext@\leavevmode
 \DN@{\ifx\next-\DN@-{\FN@\nextii@}\else
  \DN@{\leavevmode\hbox{-}}\fi\next@}%
 \DNii@{\ifx\next-\DN@-{\leavevmode\hbox{---}}\else
  \DN@{\leavevmode\hbox{--}}\fi\next@}%
 \FN@\next@}
\def\srdr@{\kern.16667em}
\def\drsr@{\kern.02778em}
\def\sldl@{\drsr@}
\def\dlsl@{\srdr@}
\atdef@"{\unskip\relaxnext@
 \DN@{\ifx\next\space@\DN@. {\FN@\nextii@}\else
  \DN@.{\FN@\nextii@}\fi\next@.}%
 \DNii@{\ifx\next`\DN@`{\FN@\nextiii@}\else
  \ifx\next\lq\DN@\lq{\FN@\nextiii@}\else
  \DN@####1{\FN@\nextiv@}\fi\fi\next@}%
 \def\nextiii@{\ifx\next`\DN@`{\sldl@``}\else\ifx\next\lq
  \DN@\lq{\sldl@``}\else\DN@{\dlsl@`}\fi\fi\next@}%
 \def\nextiv@{\ifx\next'\DN@'{\srdr@''}\else
  \ifx\next\rq\DN@\rq{\srdr@''}\else\DN@{\drsr@'}\fi\fi\next@}%
 \FN@\next@}

\def\textfontii{\the\textfont\tw@}
\def\lbrace@{\delimiter"4266308 }
\def\rbrace@{\delimiter"5267309 }
\def\{{\RIfM@\lbrace@\else{\textfontii f}\spacefactor\@m\fi}
\def\}{\RIfM@\rbrace@\else
 \let\@sf\empty\ifhmode\edef\@sf{\spacefactor\the\spacefactor}\fi
 {\textfontii g}\@sf\relax\fi}
\let\lbrace\{
\let\rbrace\}
\def\AmSTeX{{\textfontii A\kern-.1667em%
  \lower.5ex\hbox{M}\kern-.125emS}-\TeX}
\message{line and page breaks,}
\def\vmodeerr@#1{\Err@{\string#1\space not allowed between paragraphs}}
\def\mathmodeerr@#1{\Err@{\string#1\space not allowed in math mode}}
\def\linebreak{\RIfM@\mathmodeerr@\linebreak\else
 \ifhmode\unskip\unkern\break\else\vmodeerr@\linebreak\fi\fi}

\newskip\saveskip@
\def\allowlinebreak{\RIfM@\mathmodeerr@\allowlinebreak\else
 \ifhmode\saveskip@\lastskip\unskip
 \allowbreak\ifdim\saveskip@>\z@\hskip\saveskip@\fi
 \else\vmodeerr@\allowlinebreak\fi\fi}
\def\nolinebreak{\RIfM@\mathmodeerr@\nolinebreak\else
 \ifhmode\saveskip@\lastskip\unskip
 \nobreak\ifdim\saveskip@>\z@\hskip\saveskip@\fi
 \else\vmodeerr@\nolinebreak\fi\fi}
\def\newline{\relaxnext@
 \DN@{\RIfM@\expandafter\mathmodeerr@\expandafter\newline\else
  \ifhmode\ifx\next\par\else
  \expandafter\unskip\expandafter\null\expandafter\hfill\expandafter\break\fi
  \else
  \expandafter\vmodeerr@\expandafter\newline\fi\fi}%
 \FN@\next@}
\def\dmatherr@#1{\Err@{\string#1\space not allowed in display math mode}}
\def\nondmatherr@#1{\Err@{\string#1\space not allowed in non-display math
 mode}}
\def\onlydmatherr@#1{\Err@{\string#1\space allowed only in display math mode}}
\def\nonmatherr@#1{\Err@{\string#1\space allowed only in math mode}}
\def\mathbreak{\RIfMIfI@\break\else
 \dmatherr@\mathbreak\fi\else\nonmatherr@\mathbreak\fi}
\def\nomathbreak{\RIfMIfI@\nobreak\else
 \dmatherr@\nomathbreak\fi\else\nonmatherr@\nomathbreak\fi}
\def\allowmathbreak{\RIfMIfI@\allowbreak\else
 \dmatherr@\allowmathbreak\fi\else\nonmatherr@\allowmathbreak\fi}
\def\pagebreak{\RIfM@
 \ifinner\nondmatherr@\pagebreak\else\postdisplaypenalty-\@M\fi
 \else\ifvmode\removelastskip\break\else\vadjust{\break}\fi\fi}
\def\nopagebreak{\RIfM@
 \ifinner\nondmatherr@\nopagebreak\else\postdisplaypenalty\@M\fi
 \else\ifvmode\nobreak\else\vadjust{\nobreak}\fi\fi}
\def\nonvmodeerr@#1{\Err@{\string#1\space not allowed within a paragraph
 or in math}}
\def\vnonvmode@#1#2{\relaxnext@\DNii@{\ifx\next\par\DN@{#1}\else
 \DN@{#2}\fi\next@}%
 \ifvmode\DN@{#1}\else
 \DN@{\FN@\nextii@}\fi\next@}
\def\newpage{\vnonvmode@{\vfill\break}{\nonvmodeerr@\newpage}}
\def\smallpagebreak{\vnonvmode@\smallbreak{\nonvmodeerr@\smallpagebreak}}
\def\medpagebreak{\vnonvmode@\medbreak{\nonvmodeerr@\medpagebreak}}
\def\bigpagebreak{\vnonvmode@\bigbreak{\nonvmodeerr@\bigpagebreak}}
\def\NoBlackBoxes{\global\overfullrule\z@}
\def\BlackBoxes{\global\overfullrule5\p@}
\def\Invalid@#1{\def#1{\Err@{\Invalid@@\string#1}}}
\def\Invalid@@{Invalid use of }
\message{figures,}
\Invalid@\caption
\Invalid@\captionwidth
\newdimen\smallcaptionwidth@
\def\topspace{\mid@false\ins@}
\def\midspace{\mid@true\ins@}
\newif\ifmid@
\def\captionfont@{}
\def\ins@#1{\relaxnext@\allowbreak
 \smallcaptionwidth@\captionwidth@\gdef\thespace@{#1}%
 \DN@{\ifx\next\space@\DN@. {\FN@\nextii@}\else
  \DN@.{\FN@\nextii@}\fi\next@.}%
 \DNii@{\ifx\next\caption\DN@\caption{\FN@\nextiii@}%
  \else\let\next@\nextiv@\fi\next@}%
 \def\nextiv@{\vnonvmode@
  {\ifmid@\expandafter\midinsert\else\expandafter\topinsert\fi
   \vbox to\thespace@{}\endinsert}
  {\ifmid@\nonvmodeerr@\midspace\else\nonvmodeerr@\topspace\fi}}%
 \def\nextiii@{\ifx\next\captionwidth\expandafter\nextv@
  \else\expandafter\nextvi@\fi}%
 \def\nextv@\captionwidth##1##2{\smallcaptionwidth@##1\relax\nextvi@{##2}}%
 \def\nextvi@##1{\def\thecaption@{\captionfont@##1}%
  \DN@{\ifx\next\space@\DN@. {\FN@\nextvii@}\else
   \DN@.{\FN@\nextvii@}\fi\next@.}%
  \FN@\next@}%
 \def\nextvii@{\vnonvmode@
  {\ifmid@\expandafter\midinsert\else
  \expandafter\topinsert\fi\vbox to\thespace@{}\nobreak\smallskip
  \setboxz@h{\noindent\ignorespaces\thecaption@\unskip}%
  \ifdim\wdz@>\smallcaptionwidth@\centerline{\vbox{\hsize\smallcaptionwidth@
   \noindent\ignorespaces\thecaption@\unskip}}%
  \else\centerline{\boxz@}\fi\endinsert}
  {\ifmid@\nonvmodeerr@\midspace
  \else\nonvmodeerr@\topspace\fi}}%
 \FN@\next@}
\message{comments,}
\def\newcodes@{\catcode`\\12\catcode`\{12\catcode`\}12\catcode`\#12%
 \catcode`\%12\relax}
\def\oldcodes@{\catcode`\\0\catcode`\{1\catcode`\}2\catcode`\#6%
 \catcode`\%14\relax}
\def\comment{\newcodes@\endlinechar=10 \comment@}
{\lccode`\0=`\\
\lowercase{\gdef\comment@#1^^J{\comment@@#10endcomment\comment@@@}%
\gdef\comment@@#10endcomment{\FN@\comment@@@}%
\gdef\comment@@@#1\comment@@@{\ifx\next\comment@@@\let\next\comment@
 \else\def\next{\oldcodes@\endlinechar=`\^^M\relax}%
 \fi\next}}}
\def\pr@m@s{\ifx'\next\DN@##1{\prim@s}\else\let\next@\egroup\fi\next@}
\def\prime{{\null\prime@\null}}
\mathchardef\prime@="0230
\let\dsize\displaystyle

\let\ssize\scriptstyle

\message{math spacing,}
\def\,{\RIfM@\mskip\thinmuskip\relax\else\kern.16667em\fi}
\def\!{\RIfM@\mskip-\thinmuskip\relax\else\kern-.16667em\fi}
\let\thinspace\,
\let\negthinspace\!
\def\medspace{\RIfM@\mskip\medmuskip\relax\else\kern.222222em\fi}
\def\negmedspace{\RIfM@\mskip-\medmuskip\relax\else\kern-.222222em\fi}
\def\thickspace{\RIfM@\mskip\thickmuskip\relax\else\kern.27777em\fi}
\let\;\thickspace
\def\negthickspace{\RIfM@\mskip-\thickmuskip\relax\else
 \kern-.27777em\fi}
\atdef@,{\RIfM@\mskip.1\thinmuskip\else\leavevmode\null,\fi}
\atdef@!{\RIfM@\mskip-.1\thinmuskip\else\leavevmode\null!\fi}
\atdef@.{\RIfM@&&\else\leavevmode.\spacefactor3000 \fi}
\def\and{\DOTSB\;\mathchar"3026 \;}

\message{fractions,}
\def\frac#1#2{{#1\over#2}}

\newdimen\ex@
\ex@.2326ex
\Invalid@\thickness
\def\thickfrac{\relaxnext@
 \DN@{\ifx\next\thickness\let\next@\nextii@\else
 \DN@{\nextii@\thickness1}\fi\next@}%
 \DNii@\thickness##1##2##3{{##2\above##1\ex@##3}}%
 \FN@\next@}

\def\thickfracwithdelims#1#2{\relaxnext@\def\ldelim@{#1}\def\rdelim@{#2}%
 \DN@{\ifx\next\thickness\let\next@\nextii@\else
 \DN@{\nextii@\thickness1}\fi\next@}%
 \DNii@\thickness##1##2##3{{##2\abovewithdelims
 \ldelim@\rdelim@##1\ex@##3}}%
 \FN@\next@}
\def\binom#1#2{{#1\choose#2}}

\def\:{\nobreak\hskip.1111em\mathpunct{}\nonscript\mkern-\thinmuskip{:}\hskip
 .3333emplus.0555em\relax}
\def\snug{\unskip\kern-\mathsurround}
\message{smash commands,}
\def\topsmash{\top@true\bot@false\smash@}
\def\botsmash{\top@false\bot@true\smash@}
\newif\iftop@
\newif\ifbot@
\def\smash{\top@true\bot@true\smash@}
\def\smash@{\RIfM@\expandafter\mathpalette\expandafter\mathsm@sh\else
 \expandafter\makesm@sh\fi}
\def\finsm@sh{\iftop@\ht\z@\z@\fi\ifbot@\dp\z@\z@\fi\leavevmode\boxz@}
\message{large operator symbols,}
\def\LimitsOnSums{\global\let\slimits@\displaylimits}
\def\NoLimitsOnSums{\global\let\slimits@\nolimits}
\LimitsOnSums
\mathchardef\coprod@="1360       \def\coprod{\DOTSB\coprod@\slimits@}
\mathchardef\bigvee@="1357       \def\bigvee{\DOTSB\bigvee@\slimits@}
\mathchardef\bigwedge@="1356     \def\bigwedge{\DOTSB\bigwedge@\slimits@}
\mathchardef\biguplus@="1355     \def\biguplus{\DOTSB\biguplus@\slimits@}
\mathchardef\bigcap@="1354       \def\bigcap{\DOTSB\bigcap@\slimits@}
\mathchardef\bigcup@="1353       \def\bigcup{\DOTSB\bigcup@\slimits@}
\mathchardef\prod@="1351         \def\prod{\DOTSB\prod@\slimits@}
\mathchardef\sum@="1350          \def\sum{\DOTSB\sum@\slimits@}
\mathchardef\bigotimes@="134E    \def\bigotimes{\DOTSB\bigotimes@\slimits@}
\mathchardef\bigoplus@="134C     \def\bigoplus{\DOTSB\bigoplus@\slimits@}
\mathchardef\bigodot@="134A      \def\bigodot{\DOTSB\bigodot@\slimits@}
\mathchardef\bigsqcup@="1346     \def\bigsqcup{\DOTSB\bigsqcup@\slimits@}
\message{integrals,}
\def\LimitsOnInts{\global\let\ilimits@\displaylimits}
\def\NoLimitsOnInts{\global\let\ilimits@\nolimits}
\NoLimitsOnInts
\def\int{\DOTSI\intop\ilimits@}
\def\oint{\DOTSI\ointop\ilimits@}
\def\intic@{\mathchoice{\hskip.5em}{\hskip.4em}{\hskip.4em}{\hskip.4em}}
\def\negintic@{\mathchoice
 {\hskip-.5em}{\hskip-.4em}{\hskip-.4em}{\hskip-.4em}}
\def\intkern@{\mathchoice{\!\!\!}{\!\!}{\!\!}{\!\!}}
\def\intdots@{\mathchoice{\plaincdots@}
 {{\cdotp}\mkern1.5mu{\cdotp}\mkern1.5mu{\cdotp}}
 {{\cdotp}\mkern1mu{\cdotp}\mkern1mu{\cdotp}}
 {{\cdotp}\mkern1mu{\cdotp}\mkern1mu{\cdotp}}}
\newcount\intno@
\def\iint{\DOTSI\intno@\tw@\FN@\ints@}
\def\iiint{\DOTSI\intno@\thr@@\FN@\ints@}
\def\iiiint{\DOTSI\intno@4 \FN@\ints@}
\def\idotsint{\DOTSI\intno@\z@\FN@\ints@}
\def\ints@{\findlimits@\ints@@}
\newif\iflimtoken@
\newif\iflimits@
\def\findlimits@{\limtoken@true\ifx\next\limits\limits@true
 \else\ifx\next\nolimits\limits@false\else
 \limtoken@false\ifx\ilimits@\nolimits\limits@false\else
 \ifinner\limits@false\else\limits@true\fi\fi\fi\fi}
\def\multint@{\int\ifnum\intno@=\z@\intdots@                                
 \else\intkern@\fi                                                          
 \ifnum\intno@>\tw@\int\intkern@\fi                                         
 \ifnum\intno@>\thr@@\int\intkern@\fi                                       
 \int}                                                                      
\def\multintlimits@{\intop\ifnum\intno@=\z@\intdots@\else\intkern@\fi
 \ifnum\intno@>\tw@\intop\intkern@\fi
 \ifnum\intno@>\thr@@\intop\intkern@\fi\intop}
\def\ints@@{\iflimtoken@                                                    
 \def\ints@@@{\iflimits@\negintic@\mathop{\intic@\multintlimits@}\limits    
  \else\multint@\nolimits\fi                                                
  \eat@}                                                                    
 \else                                                                      
 \def\ints@@@{\iflimits@\negintic@
  \mathop{\intic@\multintlimits@}\limits\else
  \multint@\nolimits\fi}\fi\ints@@@}
\def\LimitsOnNames{\global\let\nlimits@\displaylimits}
\def\NoLimitsOnNames{\global\let\nlimits@\nolimits@}
\LimitsOnNames
\def\nolimits@{\relaxnext@
 \DN@{\ifx\next\limits\DN@\limits{\nolimits}\else
  \let\next@\nolimits\fi\next@}%
 \FN@\next@}
\message{operator names,}
\def\newmcodes@{\mathcode`\'"27\mathcode`\*"2A\mathcode`\."613A%
 \mathcode`\-"2D\mathcode`\/"2F\mathcode`\:"603A }
\def\operatorname#1{\mathop{\newmcodes@\kern\z@\fam\z@#1}\nolimits@}
\def\operatornamewithlimits#1{\mathop{\newmcodes@\kern\z@\fam\z@#1}\nlimits@}
\def\qopname@#1{\mathop{\fam\z@#1}\nolimits@}
\def\qopnamewl@#1{\mathop{\fam\z@#1}\nlimits@}
\def\arccos{\qopname@{arccos}}
\def\arcsin{\qopname@{arcsin}}
\def\arctan{\qopname@{arctan}}
\def\arg{\qopname@{arg}}
\def\cos{\qopname@{cos}}
\def\cosh{\qopname@{cosh}}
\def\cot{\qopname@{cot}}
\def\coth{\qopname@{coth}}
\def\csc{\qopname@{csc}}
\def\deg{\qopname@{deg}}
\def\det{\qopnamewl@{det}}
\def\dim{\qopname@{dim}}
\def\exp{\qopname@{exp}}
\def\gcd{\qopnamewl@{gcd}}
\def\hom{\qopname@{hom}}
\def\inf{\qopnamewl@{inf}}
\def\injlim{\qopnamewl@{inj\,lim}}
\def\ker{\qopname@{ker}}
\def\lg{\qopname@{lg}}
\def\lim{\qopnamewl@{lim}}
\def\liminf{\qopnamewl@{lim\,inf}}
\def\limsup{\qopnamewl@{lim\,sup}}
\def\ln{\qopname@{ln}}
\def\log{\qopname@{log}}
\def\max{\qopnamewl@{max}}
\def\min{\qopnamewl@{min}}
\def\Pr{\qopnamewl@{Pr}}
\def\projlim{\qopnamewl@{proj\,lim}}
\def\sec{\qopname@{sec}}
\def\sin{\qopname@{sin}}
\def\sinh{\qopname@{sinh}}
\def\sup{\qopnamewl@{sup}}
\def\tan{\qopname@{tan}}
\def\tanh{\qopname@{tanh}}
\def\varinjlim{\mathop{\vtop{\ialign{##\crcr
 \hfil\rm lim\hfil\crcr\noalign{\nointerlineskip}\rightarrowfill\crcr
 \noalign{\nointerlineskip\kern-\ex@}\crcr}}}}
\def\varprojlim{\mathop{\vtop{\ialign{##\crcr
 \hfil\rm lim\hfil\crcr\noalign{\nointerlineskip}\leftarrowfill\crcr
 \noalign{\nointerlineskip\kern-\ex@}\crcr}}}}
\def\varliminf{\mathop{\underline{\vrule height\z@ depth.2exwidth\z@
 \hbox{\rm lim}}}}

\newdimen\buffer@
\buffer@\fontdimen13 \tenex
\newdimen\buffer
\buffer\buffer@

\def\ResetBuffer{\fontdimen13 \tenex\buffer@\global\buffer\buffer@}
\def\shave#1{\mathop{\hbox{$\m@th\fontdimen13 \tenex\z@                     
 \displaystyle{#1}$}}\fontdimen13 \tenex\buffer}

\message{multilevel sub/superscripts,}
\Invalid@\\
\def\Let@{\relax\iffalse{\fi\let\\=\cr\iffalse}\fi}
\Invalid@\vspace
\def\vspace@{\def\vspace##1{\crcr\noalign{\vskip##1\relax}}}
\def\multilimits@{\bgroup\vspace@\Let@
 \baselineskip\fontdimen10 \scriptfont\tw@
 \advance\baselineskip\fontdimen12 \scriptfont\tw@
 \lineskip\thr@@\fontdimen8 \scriptfont\thr@@
 \lineskiplimit\lineskip
 \vbox\bgroup\ialign\bgroup\hfil$\m@th\scriptstyle{##}$\hfil\crcr}
\def\Sb{_\multilimits@}
\def\endSb{\crcr\egroup\egroup\egroup}
\def\Sp{^\multilimits@}

\def\spreadlines#1{\RIfMIfI@\onlydmatherr@\spreadlines\else
 \openup#1\relax\fi\else\onlydmatherr@\spreadlines\fi}
\def\Mathstrut@{\copy\Mathstrutbox@}
\newbox\Mathstrutbox@
\setbox\Mathstrutbox@\null
\setboxz@h{$\m@th($}
\ht\Mathstrutbox@\ht\z@
\dp\Mathstrutbox@\dp\z@
\message{matrices,}
\newdimen\spreadmlines@
\def\spreadmatrixlines#1{\RIfMIfI@
 \onlydmatherr@\spreadmatrixlines\else
 \spreadmlines@#1\relax\fi\else\onlydmatherr@\spreadmatrixlines\fi}
\def\matrix{\null\,\vcenter\bgroup\Let@\vspace@
 \normalbaselines\openup\spreadmlines@\ialign
 \bgroup\hfil$\m@th##$\hfil&&\quad\hfil$\m@th##$\hfil\crcr
 \Mathstrut@\crcr\noalign{\kern-\baselineskip}}
\def\endmatrix{\crcr\Mathstrut@\crcr\noalign{\kern-\baselineskip}\egroup
 \egroup\,}
\def\format{\crcr\egroup\iffalse{\fi\ifnum`}=0 \fi\format@}
\newtoks\hashtoks@
\hashtoks@{#}
\def\format@#1\\{\def\preamble@{#1}%
 \def\l{$\m@th\the\hashtoks@$\hfil}%
 \def\c{\hfil$\m@th\the\hashtoks@$\hfil}%
 \def\r{\hfil$\m@th\the\hashtoks@$}%
 \edef\preamble@@{\preamble@}\ifnum`{=0 \fi\iffalse}\fi
 \ialign\bgroup\span\preamble@@\crcr}
\def\smallmatrix{\null\,\vcenter\bgroup\vspace@\Let@
 \baselineskip9\ex@\lineskip\ex@
 \ialign\bgroup\hfil$\m@th\scriptstyle{##}$\hfil&&\thickspace\hfil
 $\m@th\scriptstyle{##}$\hfil\crcr}
\def\endsmallmatrix{\crcr\egroup\egroup\,}

\newmuskip\dotsspace@
\dotsspace@1.5mu
\def\strip@#1 {#1}
\def\spacehdots#1\for#2{\multispan{#2}\xleaders
 \hbox{$\m@th\mkern\strip@#1 \dotsspace@.\mkern\strip@#1 \dotsspace@$}\hfill}
\def\hdotsfor#1{\spacehdots\@ne\for{#1}}
\def\multispan@#1{\omit\mscount#1\unskip\loop\ifnum\mscount>\@ne\sp@n\repeat}
\def\spaceinnerhdots#1\for#2\after#3{\multispan@{\strip@#2 }#3\xleaders
 \hbox{$\m@th\mkern\strip@#1 \dotsspace@.\mkern\strip@#1 \dotsspace@$}\hfill}
\def\innerhdotsfor#1\after#2{\spaceinnerhdots\@ne\for#1\after{#2}}
\def\cases{\bgroup\spreadmlines@\jot\left\{\,\matrix\format\l&\quad\l\\}
\def\endcases{\endmatrix\right.\egroup}
\message{multiline displays,}
\newif\ifinany@
\newif\ifinalign@
\newif\ifingather@
\def\strut@{\copy\strutbox@}
\newbox\strutbox@
\setbox\strutbox@\hbox{\vrule height8\p@ depth3\p@ width\z@}
\def\topaligned{\null\,\vtop\aligned@}
\def\botaligned{\null\,\vbox\aligned@}
\def\aligned{\null\,\vcenter\aligned@}
\def\aligned@{\bgroup\vspace@\Let@
 \ifinany@\else\openup\jot\fi\ialign
 \bgroup\hfil\strut@$\m@th\displaystyle{##}$&
 $\m@th\displaystyle{{}##}$\hfil\crcr}
\def\endaligned{\crcr\egroup\egroup}

\def\alignedat#1{\null\,\vcenter\bgroup\doat@{#1}\vspace@\Let@
 \ifinany@\else\openup\jot\fi\ialign\bgroup\span\preamble@@\crcr}
\newcount\atcount@
\def\doat@#1{\toks@{\hfil\strut@$\m@th
 \displaystyle{\the\hashtoks@}$&$\m@th\displaystyle
 {{}\the\hashtoks@}$\hfil}
 \atcount@#1\relax\advance\atcount@\m@ne                                    
 \loop\ifnum\atcount@>\z@\toks@=\expandafter{\the\toks@&\hfil$\m@th
 \displaystyle{\the\hashtoks@}$&$\m@th
 \displaystyle{{}\the\hashtoks@}$\hfil}\advance
  \atcount@\m@ne\repeat                                                     
 \xdef\preamble@{\the\toks@}\xdef\preamble@@{\preamble@}}

\def\gathered{\null\,\vcenter\bgroup\vspace@\Let@
 \ifinany@\else\openup\jot\fi\ialign
 \bgroup\hfil\strut@$\m@th\displaystyle{##}$\hfil\crcr}
\def\endgathered{\crcr\egroup\egroup}
\newif\iftagsleft@
\def\TagsOnLeft{\global\tagsleft@true}
\def\TagsOnRight{\global\tagsleft@false}
\TagsOnLeft
\newif\ifmathtags@
\def\TagsAsMath{\global\mathtags@true}
\def\TagsAsText{\global\mathtags@false}
\TagsAsText
\def\tagform@#1{\hbox{\rm(\ignorespaces#1\unskip)}}
\def\thetag{\leavevmode\tagform@}
\def\tag#1$${\iftagsleft@\leqno\else\eqno\fi                                
 \maketag@#1\maketag@                                                       
 $$}                                                                        
\def\maketag@{\FN@\maketag@@}
\def\maketag@@{\ifx\next"\expandafter\maketag@@@\else\expandafter\maketag@@@@
 \fi}
\def\maketag@@@"#1"#2\maketag@{\hbox{\rm#1}}                                
\def\maketag@@@@#1\maketag@{\ifmathtags@\tagform@{$\m@th#1$}\else
 \tagform@{#1}\fi}
\interdisplaylinepenalty\@M
\def\allowdisplaybreaks{\RIfMIfI@
 \onlydmatherr@\allowdisplaybreaks\else
 \interdisplaylinepenalty\z@\fi\else\onlydmatherr@\allowdisplaybreaks\fi}
\Invalid@\allowdisplaybreak
\Invalid@\displaybreak
\Invalid@\intertext
\def\allowdisplaybreak@{\def\allowdisplaybreak{\crcr\noalign{\allowbreak}}}
\def\displaybreak@{\def\displaybreak{\crcr\noalign{\break}}}
\def\intertext@{\def\intertext##1{\crcr\noalign{%
 \penalty\postdisplaypenalty \vskip\belowdisplayskip
 \vbox{\normalbaselines\noindent##1}%
 \penalty\predisplaypenalty \vskip\abovedisplayskip}}}
\newskip\centering@
\centering@\z@ plus\@m\p@
\def\align{\relax\ifingather@\DN@{\csname align (in
  \string\gather)\endcsname}\else
 \ifmmode\ifinner\DN@{\onlydmatherr@\align}\else
  \let\next@\align@\fi
 \else\DN@{\onlydmatherr@\align}\fi\fi\next@}
\newhelp\andhelp@
{An extra & here is so disastrous that you should probably exit^^J
and fix things up.}
\newif\iftag@
\newcount\and@
\def\align@{\inalign@true\inany@true
 \vspace@\allowdisplaybreak@\displaybreak@\intertext@
 \def\tag{\global\tag@true\ifnum\and@=\z@\DN@{&&}\else
          \DN@{&}\fi\next@}%
 \iftagsleft@\DN@{\csname align \endcsname}\else
  \DN@{\csname align \space\endcsname}\fi\next@}
\def\Tag@{\iftag@\else\errhelp\andhelp@\err@{Extra & on this line}\fi}
\newdimen\lwidth@
\newdimen\rwidth@
\newdimen\maxlwidth@
\newdimen\maxrwidth@
\newdimen\totwidth@
\def\measure@#1\endalign{\lwidth@\z@\rwidth@\z@\maxlwidth@\z@\maxrwidth@\z@
 \global\and@\z@                                                            
 \setbox@ne\vbox                                                            
  {\everycr{\noalign{\global\tag@false\global\and@\z@}}\Let@                
  \halign{\setboxz@h{$\m@th\displaystyle{\@lign##}$}
   \global\lwidth@\wdz@                                                     
   \ifdim\lwidth@>\maxlwidth@\global\maxlwidth@\lwidth@\fi                  
   \global\advance\and@\@ne                                                 
   &\setboxz@h{$\m@th\displaystyle{{}\@lign##}$}\global\rwidth@\wdz@        
   \ifdim\rwidth@>\maxrwidth@\global\maxrwidth@\rwidth@\fi                  
   \global\advance\and@\@ne                                                
   &\Tag@
   \eat@{##}\crcr#1\crcr}}
 \totwidth@\maxlwidth@\advance\totwidth@\maxrwidth@}                       
\def\displ@y@{\global\dt@ptrue\openup\jot
 \everycr{\noalign{\global\tag@false\global\and@\z@\ifdt@p\global\dt@pfalse
 \vskip-\lineskiplimit\vskip\normallineskiplimit\else
 \penalty\interdisplaylinepenalty\fi}}}
\def\black@#1{\noalign{\ifdim#1>\displaywidth
 \dimen@\prevdepth\nointerlineskip                                          
 \vskip-\ht\strutbox@\vskip-\dp\strutbox@                                   
 \vbox{\noindent\hbox to#1{\strut@\hfill}}
 \prevdepth\dimen@                                                          
 \fi}}
\expandafter\def\csname align \space\endcsname#1\endalign
 {\measure@#1\endalign\global\and@\z@                                       
 \ifingather@\everycr{\noalign{\global\and@\z@}}\else\displ@y@\fi           
 \Let@\tabskip\centering@                                                   
 \halign to\displaywidth
  {\hfil\strut@\setboxz@h{$\m@th\displaystyle{\@lign##}$}
  \global\lwidth@\wdz@\boxz@\global\advance\and@\@ne                        
  \tabskip\z@skip                                                           
  &\setboxz@h{$\m@th\displaystyle{{}\@lign##}$}
  \global\rwidth@\wdz@\boxz@\hfill\global\advance\and@\@ne                  
  \tabskip\centering@                                                       
  &\setboxz@h{\@lign\strut@\maketag@##\maketag@}
  \dimen@\displaywidth\advance\dimen@-\totwidth@
  \divide\dimen@\tw@\advance\dimen@\maxrwidth@\advance\dimen@-\rwidth@     
  \ifdim\dimen@<\tw@\wdz@\llap{\vtop{\normalbaselines\null\boxz@}}
  \else\llap{\boxz@}\fi                                                    
  \tabskip\z@skip                                                          
  \crcr#1\crcr                                                             
  \black@\totwidth@}}                                                      
\newdimen\lineht@
\expandafter\def\csname align \endcsname#1\endalign{\measure@#1\endalign
 \global\and@\z@
 \ifdim\totwidth@>\displaywidth\let\displaywidth@\totwidth@\else
  \let\displaywidth@\displaywidth\fi                                        
 \ifingather@\everycr{\noalign{\global\and@\z@}}\else\displ@y@\fi
 \Let@\tabskip\centering@\halign to\displaywidth
  {\hfil\strut@\setboxz@h{$\m@th\displaystyle{\@lign##}$}%
  \global\lwidth@\wdz@\global\lineht@\ht\z@                                 
  \boxz@\global\advance\and@\@ne
  \tabskip\z@skip&\setboxz@h{$\m@th\displaystyle{{}\@lign##}$}%
  \global\rwidth@\wdz@\ifdim\ht\z@>\lineht@\global\lineht@\ht\z@\fi         
  \boxz@\hfil\global\advance\and@\@ne
  \tabskip\centering@&\kern-\displaywidth@                                  
  \setboxz@h{\@lign\strut@\maketag@##\maketag@}%
  \dimen@\displaywidth\advance\dimen@-\totwidth@
  \divide\dimen@\tw@\advance\dimen@\maxlwidth@\advance\dimen@-\lwidth@
  \ifdim\dimen@<\tw@\wdz@
   \rlap{\vbox{\normalbaselines\boxz@\vbox to\lineht@{}}}\else
   \rlap{\boxz@}\fi
  \tabskip\displaywidth@\crcr#1\crcr\black@\totwidth@}}
\expandafter\def\csname align (in \string\gather)\endcsname
  #1\endalign{\vcenter{\align@#1\endalign}}
\Invalid@\endalign
\newif\ifxat@
\def\alignat{\RIfMIfI@\DN@{\onlydmatherr@\alignat}\else
 \DN@{\csname alignat \endcsname}\fi\else
 \DN@{\onlydmatherr@\alignat}\fi\next@}
\newif\ifmeasuring@
\newbox\savealignat@
\expandafter\def\csname alignat \endcsname#1#2\endalignat                   
 {\inany@true\xat@false
 \def\tag{\global\tag@true\count@#1\relax\multiply\count@\tw@
  \xdef\tag@{}\loop\ifnum\count@>\and@\xdef\tag@{&\tag@}\advance\count@\m@ne
  \repeat\tag@}%
 \vspace@\allowdisplaybreak@\displaybreak@\intertext@
 \displ@y@\measuring@true                                                   
 \setbox\savealignat@\hbox{$\m@th\displaystyle\Let@
  \attag@{#1}
  \vbox{\halign{\span\preamble@@\crcr#2\crcr}}$}%
 \measuring@false                                                           
 \Let@\attag@{#1}
 \tabskip\centering@\halign to\displaywidth
  {\span\preamble@@\crcr#2\crcr                                             
  \black@{\wd\savealignat@}}}                                               
\Invalid@\endalignat
\def\xalignat{\RIfMIfI@
 \DN@{\onlydmatherr@\xalignat}\else
 \DN@{\csname xalignat \endcsname}\fi\else
 \DN@{\onlydmatherr@\xalignat}\fi\next@}
\expandafter\def\csname xalignat \endcsname#1#2\endxalignat
 {\inany@true\xat@true
 \def\tag{\global\tag@true\def\tag@{}\count@#1\relax\multiply\count@\tw@
  \loop\ifnum\count@>\and@\xdef\tag@{&\tag@}\advance\count@\m@ne\repeat\tag@}%
 \vspace@\allowdisplaybreak@\displaybreak@\intertext@
 \displ@y@\measuring@true\setbox\savealignat@\hbox{$\m@th\displaystyle\Let@
 \attag@{#1}\vbox{\halign{\span\preamble@@\crcr#2\crcr}}$}%
 \measuring@false\Let@
 \attag@{#1}\tabskip\centering@\halign to\displaywidth
 {\span\preamble@@\crcr#2\crcr\black@{\wd\savealignat@}}}
\def\attag@#1{\let\Maketag@\maketag@\let\TAG@\Tag@                          
 \let\Tag@=0\let\maketag@=0
 \ifmeasuring@\def\llap@##1{\setboxz@h{##1}\hbox to\tw@\wdz@{}}%
  \def\rlap@##1{\setboxz@h{##1}\hbox to\tw@\wdz@{}}\else
  \let\llap@\llap\let\rlap@\rlap\fi                                         
 \toks@{\hfil\strut@$\m@th\displaystyle{\@lign\the\hashtoks@}$\tabskip\z@skip
  \global\advance\and@\@ne&$\m@th\displaystyle{{}\@lign\the\hashtoks@}$\hfil
  \ifxat@\tabskip\centering@\fi\global\advance\and@\@ne}
 \iftagsleft@
  \toks@@{\tabskip\centering@&\Tag@\kern-\displaywidth
   \rlap@{\@lign\maketag@\the\hashtoks@\maketag@}%
   \global\advance\and@\@ne\tabskip\displaywidth}\else
  \toks@@{\tabskip\centering@&\Tag@\llap@{\@lign\maketag@
   \the\hashtoks@\maketag@}\global\advance\and@\@ne\tabskip\z@skip}\fi      
 \atcount@#1\relax\advance\atcount@\m@ne
 \loop\ifnum\atcount@>\z@
 \toks@=\expandafter{\the\toks@&\hfil$\m@th\displaystyle{\@lign
  \the\hashtoks@}$\global\advance\and@\@ne
  \tabskip\z@skip&$\m@th\displaystyle{{}\@lign\the\hashtoks@}$\hfil\ifxat@
  \tabskip\centering@\fi\global\advance\and@\@ne}\advance\atcount@\m@ne
 \repeat                                                                    
 \xdef\preamble@{\the\toks@\the\toks@@}
 \xdef\preamble@@{\preamble@}
 \let\maketag@\Maketag@\let\Tag@\TAG@}                                      
\Invalid@\endxalignat
\def\xxalignat{\RIfMIfI@
 \DN@{\onlydmatherr@\xxalignat}\else\DN@{\csname xxalignat
  \endcsname}\fi\else
 \DN@{\onlydmatherr@\xxalignat}\fi\next@}
\expandafter\def\csname xxalignat \endcsname#1#2\endxxalignat{\inany@true
 \vspace@\allowdisplaybreak@\displaybreak@\intertext@
 \displ@y\setbox\savealignat@\hbox{$\m@th\displaystyle\Let@
 \xxattag@{#1}\vbox{\halign{\span\preamble@@\crcr#2\crcr}}$}%
 \Let@\xxattag@{#1}\tabskip\z@skip\halign to\displaywidth
 {\span\preamble@@\crcr#2\crcr\black@{\wd\savealignat@}}}
\def\xxattag@#1{\toks@{\tabskip\z@skip\hfil\strut@
 $\m@th\displaystyle{\the\hashtoks@}$&%
 $\m@th\displaystyle{{}\the\hashtoks@}$\hfil\tabskip\centering@&}%
 \atcount@#1\relax\advance\atcount@\m@ne\loop\ifnum\atcount@>\z@
 \toks@=\expandafter{\the\toks@&\hfil$\m@th\displaystyle{\the\hashtoks@}$%
  \tabskip\z@skip&$\m@th\displaystyle{{}\the\hashtoks@}$\hfil
  \tabskip\centering@}\advance\atcount@\m@ne\repeat
 \xdef\preamble@{\the\toks@\tabskip\z@skip}\xdef\preamble@@{\preamble@}}
\Invalid@\endxxalignat
\newdimen\gwidth@
\newdimen\gmaxwidth@
\def\gmeasure@#1\endgather{\gwidth@\z@\gmaxwidth@\z@\setbox@ne\vbox{\Let@
 \halign{\setboxz@h{$\m@th\displaystyle{##}$}\global\gwidth@\wdz@
 \ifdim\gwidth@>\gmaxwidth@\global\gmaxwidth@\gwidth@\fi
 &\eat@{##}\crcr#1\crcr}}}
\def\gather{\RIfMIfI@\DN@{\onlydmatherr@\gather}\else
 \ingather@true\inany@true\def\tag{&}%
 \vspace@\allowdisplaybreak@\displaybreak@\intertext@
 \displ@y\Let@
 \iftagsleft@\DN@{\csname gather \endcsname}\else
  \DN@{\csname gather \space\endcsname}\fi\fi
 \else\DN@{\onlydmatherr@\gather}\fi\next@}
\expandafter\def\csname gather \space\endcsname#1\endgather
 {\gmeasure@#1\endgather\tabskip\centering@
 \halign to\displaywidth{\hfil\strut@\setboxz@h{$\m@th\displaystyle{##}$}%
 \global\gwidth@\wdz@\boxz@\hfil&
 \setboxz@h{\strut@{\maketag@##\maketag@}}%
 \dimen@\displaywidth\advance\dimen@-\gwidth@
 \ifdim\dimen@>\tw@\wdz@\llap{\boxz@}\else
 \llap{\vtop{\normalbaselines\null\boxz@}}\fi
 \tabskip\z@skip\crcr#1\crcr\black@\gmaxwidth@}}
\newdimen\glineht@
\expandafter\def\csname gather \endcsname#1\endgather{\gmeasure@#1\endgather
 \ifdim\gmaxwidth@>\displaywidth\let\gdisplaywidth@\gmaxwidth@\else
 \let\gdisplaywidth@\displaywidth\fi\tabskip\centering@\halign to\displaywidth
 {\hfil\strut@\setboxz@h{$\m@th\displaystyle{##}$}%
 \global\gwidth@\wdz@\global\glineht@\ht\z@\boxz@\hfil&\kern-\gdisplaywidth@
 \setboxz@h{\strut@{\maketag@##\maketag@}}%
 \dimen@\displaywidth\advance\dimen@-\gwidth@
 \ifdim\dimen@>\tw@\wdz@\rlap{\boxz@}\else
 \rlap{\vbox{\normalbaselines\boxz@\vbox to\glineht@{}}}\fi
 \tabskip\gdisplaywidth@\crcr#1\crcr\black@\gmaxwidth@}}
\newif\ifctagsplit@
\def\CenteredTagsOnSplits{\global\ctagsplit@true}
\def\TopOrBottomTagsOnSplits{\global\ctagsplit@false}
\TopOrBottomTagsOnSplits
\def\split{\relax\ifinany@\let\next@\insplit@\else
 \ifmmode\ifinner\def\next@{\onlydmatherr@\split}\else
 \let\next@\outsplit@\fi\else
 \def\next@{\onlydmatherr@\split}\fi\fi\next@}
\def\insplit@{\global\setbox\z@\vbox\bgroup\vspace@\Let@\ialign\bgroup
 \hfil\strut@$\m@th\displaystyle{##}$&$\m@th\displaystyle{{}##}$\hfill\crcr}
\def\endsplit{\crcr\egroup\egroup\iftagsleft@\expandafter\lendsplit@\else
 \expandafter\rendsplit@\fi}
\def\rendsplit@{\global\setbox9 \vbox
 {\unvcopy\z@\global\setbox8 \lastbox\unskip}
 \setbox@ne\hbox{\unhcopy8 \unskip\global\setbox\tw@\lastbox
 \unskip\global\setbox\thr@@\lastbox}
 \global\setbox7 \hbox{\unhbox\tw@\unskip}
 \ifinalign@\ifctagsplit@                                                   
  \gdef\split@{\hbox to\wd\thr@@{}&
   \vcenter{\vbox{\moveleft\wd\thr@@\boxz@}}}
 \else\gdef\split@{&\vbox{\moveleft\wd\thr@@\box9}\crcr
  \box\thr@@&\box7}\fi                                                      
 \else                                                                      
  \ifctagsplit@\gdef\split@{\vcenter{\boxz@}}\else
  \gdef\split@{\box9\crcr\hbox{\box\thr@@\box7}}\fi
 \fi
 \split@}                                                                   
\def\lendsplit@{\global\setbox9\vtop{\unvcopy\z@}
 \setbox@ne\vbox{\unvcopy\z@\global\setbox8\lastbox}
 \setbox@ne\hbox{\unhcopy8\unskip\setbox\tw@\lastbox
  \unskip\global\setbox\thr@@\lastbox}
 \ifinalign@\ifctagsplit@                                                   
  \gdef\split@{\hbox to\wd\thr@@{}&
  \vcenter{\vbox{\moveleft\wd\thr@@\box9}}}
  \else                                                                     
  \gdef\split@{\hbox to\wd\thr@@{}&\vbox{\moveleft\wd\thr@@\box9}}\fi
 \else
  \ifctagsplit@\gdef\split@{\vcenter{\box9}}\else
  \gdef\split@{\box9}\fi
 \fi\split@}
\def\outsplit@#1$${\align\insplit@#1\endalign$$}
\newdimen\multlinegap@
\multlinegap@1em
\newdimen\multlinetaggap@
\multlinetaggap@1em
\def\MultlineGap#1{\global\multlinegap@#1\relax}
\def\multlinegap#1{\RIfMIfI@\onlydmatherr@\multlinegap\else
 \multlinegap@#1\relax\fi\else\onlydmatherr@\multlinegap\fi}
\def\nomultlinegap{\multlinegap{\z@}}
\def\multline{\RIfMIfI@
 \DN@{\onlydmatherr@\multline}\else
 \DN@{\multline@}\fi\else
 \DN@{\onlydmatherr@\multline}\fi\next@}
\newif\iftagin@
\def\tagin@#1{\tagin@false\in@\tag{#1}\ifin@\tagin@true\fi}
\def\multline@#1$${\inany@true\vspace@\allowdisplaybreak@\displaybreak@
 \tagin@{#1}\iftagsleft@\DN@{\multline@l#1$$}\else
 \DN@{\multline@r#1$$}\fi\next@}
\newdimen\mwidth@
\def\rmmeasure@#1\endmultline{%
 \def\shoveleft##1{##1}\def\shoveright##1{##1}
 \setbox@ne\vbox{\Let@\halign{\setboxz@h
  {$\m@th\@lign\displaystyle{}##$}\global\mwidth@\wdz@
  \crcr#1\crcr}}}
\newdimen\mlineht@
\newif\ifzerocr@
\newif\ifonecr@
\def\lmmeasure@#1\endmultline{\global\zerocr@true\global\onecr@false
 \everycr{\noalign{\ifonecr@\global\onecr@false\fi
  \ifzerocr@\global\zerocr@false\global\onecr@true\fi}}
  \def\shoveleft##1{##1}\def\shoveright##1{##1}%
 \setbox@ne\vbox{\Let@\halign{\setboxz@h
  {$\m@th\@lign\displaystyle{}##$}\ifonecr@\global\mwidth@\wdz@
  \global\mlineht@\ht\z@\fi\crcr#1\crcr}}}
\newbox\mtagbox@
\newdimen\ltwidth@
\newdimen\rtwidth@
\def\multline@l#1$${\iftagin@\DN@{\lmultline@@#1$$}\else
 \DN@{\setbox\mtagbox@\null\ltwidth@\z@\rtwidth@\z@
  \lmultline@@@#1$$}\fi\next@}
\def\lmultline@@#1\endmultline\tag#2$${%
 \setbox\mtagbox@\hbox{\maketag@#2\maketag@}
 \lmmeasure@#1\endmultline\dimen@\mwidth@\advance\dimen@\wd\mtagbox@
 \advance\dimen@\multlinetaggap@                                            
 \ifdim\dimen@>\displaywidth\ltwidth@\z@\else\ltwidth@\wd\mtagbox@\fi       
 \lmultline@@@#1\endmultline$$}
\def\lmultline@@@{\displ@y
 \def\shoveright##1{##1\hfilneg\hskip\multlinegap@}%
 \def\shoveleft##1{\setboxz@h{$\m@th\displaystyle{}##1$}%
  \setbox@ne\hbox{$\m@th\displaystyle##1$}%
  \hfilneg
  \iftagin@
   \ifdim\ltwidth@>\z@\hskip\ltwidth@\hskip\multlinetaggap@\fi
  \else\hskip\multlinegap@\fi\hskip.5\wd@ne\hskip-.5\wdz@##1}
  \halign\bgroup\Let@\hbox to\displaywidth
   {\strut@$\m@th\displaystyle\hfil{}##\hfil$}\crcr
   \hfilneg                                                                 
   \iftagin@                                                                
    \ifdim\ltwidth@>\z@                                                     
     \box\mtagbox@\hskip\multlinetaggap@                                    
    \else
     \rlap{\vbox{\normalbaselines\hbox{\strut@\box\mtagbox@}%
     \vbox to\mlineht@{}}}\fi                                               
   \else\hskip\multlinegap@\fi}                                             
\def\multline@r#1$${\iftagin@\DN@{\rmultline@@#1$$}\else
 \DN@{\setbox\mtagbox@\null\ltwidth@\z@\rtwidth@\z@
  \rmultline@@@#1$$}\fi\next@}
\def\rmultline@@#1\endmultline\tag#2$${\ltwidth@\z@
 \setbox\mtagbox@\hbox{\maketag@#2\maketag@}%
 \rmmeasure@#1\endmultline\dimen@\mwidth@\advance\dimen@\wd\mtagbox@
 \advance\dimen@\multlinetaggap@
 \ifdim\dimen@>\displaywidth\rtwidth@\z@\else\rtwidth@\wd\mtagbox@\fi
 \rmultline@@@#1\endmultline$$}
\def\rmultline@@@{\displ@y
 \def\shoveright##1{##1\hfilneg\iftagin@\ifdim\rtwidth@>\z@
  \hskip\rtwidth@\hskip\multlinetaggap@\fi\else\hskip\multlinegap@\fi}%
 \def\shoveleft##1{\setboxz@h{$\m@th\displaystyle{}##1$}%
  \setbox@ne\hbox{$\m@th\displaystyle##1$}%
  \hfilneg\hskip\multlinegap@\hskip.5\wd@ne\hskip-.5\wdz@##1}%
 \halign\bgroup\Let@\hbox to\displaywidth
  {\strut@$\m@th\displaystyle\hfil{}##\hfil$}\crcr
 \hfilneg\hskip\multlinegap@}
\def\endmultline{\iftagsleft@\expandafter\lendmultline@\else
 \expandafter\rendmultline@\fi}
\def\lendmultline@{\hfilneg\hskip\multlinegap@\crcr\egroup}
\def\rendmultline@{\iftagin@                                                
 \ifdim\rtwidth@>\z@                                                        
  \hskip\multlinetaggap@\box\mtagbox@                                       
 \else\llap{\vtop{\normalbaselines\null\hbox{\strut@\box\mtagbox@}}}\fi     
 \else\hskip\multlinegap@\fi                                                
 \hfilneg\crcr\egroup}
\def\bmod{\mskip-\medmuskip\mkern5mu\mathbin{\fam\z@ mod}\penalty900
 \mkern5mu\mskip-\medmuskip}
\def\pmod#1{\allowbreak\ifinner\mkern8mu\else\mkern18mu\fi
 ({\fam\z@ mod}\,\,#1)}
\def\pod#1{\allowbreak\ifinner\mkern8mu\else\mkern18mu\fi(#1)}
\def\mod#1{\allowbreak\ifinner\mkern12mu\else\mkern18mu\fi{\fam\z@ mod}\,\,#1}
\message{continued fractions,}
\newcount\cfraccount@
\def\cfrac{\bgroup\bgroup\advance\cfraccount@\@ne\strut
 \iffalse{\fi\def\\{\over\displaystyle}\iffalse}\fi}
\def\lcfrac{\bgroup\bgroup\advance\cfraccount@\@ne\strut
 \iffalse{\fi\def\\{\hfill\over\displaystyle}\iffalse}\fi}
\def\rcfrac{\bgroup\bgroup\advance\cfraccount@\@ne\strut\hfill
 \iffalse{\fi\def\\{\over\displaystyle}\iffalse}\fi}
\def\gloop@#1\repeat{\gdef\body{#1}\iterate}
\def\endcfrac{\gloop@\ifnum\cfraccount@>\z@\global\advance\cfraccount@\m@ne
 \egroup\hskip-\nulldelimiterspace\egroup\repeat}
\message{compound symbols,}
\def\binrel@#1{\setboxz@h{\thinmuskip0mu
  \medmuskip\m@ne mu\thickmuskip\@ne mu$#1\m@th$}%
 \setbox@ne\hbox{\thinmuskip0mu\medmuskip\m@ne mu\thickmuskip
  \@ne mu${}#1{}\m@th$}%
 \setbox\tw@\hbox{\hskip\wd@ne\hskip-\wdz@}}
\def\overset#1\to#2{\binrel@{#2}\ifdim\wd\tw@<\z@
 \mathbin{\mathop{\kern\z@#2}\limits^{#1}}\else\ifdim\wd\tw@>\z@
 \mathrel{\mathop{\kern\z@#2}\limits^{#1}}\else
 {\mathop{\kern\z@#2}\limits^{#1}}{}\fi\fi}
\def\underset#1\to#2{\binrel@{#2}\ifdim\wd\tw@<\z@
 \mathbin{\mathop{\kern\z@#2}\limits_{#1}}\else\ifdim\wd\tw@>\z@
 \mathrel{\mathop{\kern\z@#2}\limits_{#1}}\else
 {\mathop{\kern\z@#2}\limits_{#1}}{}\fi\fi}
\def\oversetbrace#1\to#2{\overbrace{#2}^{#1}}
\def\undersetbrace#1\to#2{\underbrace{#2}_{#1}}
\def\sideset#1\and#2\to#3{%
 \setbox@ne\hbox{$\dsize{\vphantom{#3}}#1{#3}\m@th$}%
 \setbox\tw@\hbox{$\dsize{#3}#2\m@th$}%
 \hskip\wd@ne\hskip-\wd\tw@\mathop{\hskip\wd\tw@\hskip-\wd@ne
  {\vphantom{#3}}#1{#3}#2}}
\def\rightarrowfill@#1{\setboxz@h{$#1-\m@th$}\ht\z@\z@
  $#1\m@th\copy\z@\mkern-6mu\cleaders
  \hbox{$#1\mkern-2mu\box\z@\mkern-2mu$}\hfill
  \mkern-6mu\mathord\rightarrow$}
\def\leftarrowfill@#1{\setboxz@h{$#1-\m@th$}\ht\z@\z@
  $#1\m@th\mathord\leftarrow\mkern-6mu\cleaders
  \hbox{$#1\mkern-2mu\copy\z@\mkern-2mu$}\hfill
  \mkern-6mu\box\z@$}
\def\leftrightarrowfill@#1{\setboxz@h{$#1-\m@th$}\ht\z@\z@
  $#1\m@th\mathord\leftarrow\mkern-6mu\cleaders
  \hbox{$#1\mkern-2mu\box\z@\mkern-2mu$}\hfill
  \mkern-6mu\mathord\rightarrow$}
\def\overrightarrow{\mathpalette\overrightarrow@}
\def\overrightarrow@#1#2{\vbox{\ialign{##\crcr\rightarrowfill@#1\crcr
 \noalign{\kern-\ex@\nointerlineskip}$\m@th\hfil#1#2\hfil$\crcr}}}

\def\overleftarrow{\mathpalette\overleftarrow@}
\def\overleftarrow@#1#2{\vbox{\ialign{##\crcr\leftarrowfill@#1\crcr
 \noalign{\kern-\ex@\nointerlineskip}$\m@th\hfil#1#2\hfil$\crcr}}}
\def\overleftrightarrow{\mathpalette\overleftrightarrow@}
\def\overleftrightarrow@#1#2{\vbox{\ialign{##\crcr\leftrightarrowfill@#1\crcr
 \noalign{\kern-\ex@\nointerlineskip}$\m@th\hfil#1#2\hfil$\crcr}}}
\def\underrightarrow{\mathpalette\underrightarrow@}
\def\underrightarrow@#1#2{\vtop{\ialign{##\crcr$\m@th\hfil#1#2\hfil$\crcr
 \noalign{\nointerlineskip}\rightarrowfill@#1\crcr}}}

\def\underleftarrow{\mathpalette\underleftarrow@}
\def\underleftarrow@#1#2{\vtop{\ialign{##\crcr$\m@th\hfil#1#2\hfil$\crcr
 \noalign{\nointerlineskip}\leftarrowfill@#1\crcr}}}
\def\underleftrightarrow{\mathpalette\underleftrightarrow@}
\def\underleftrightarrow@#1#2{\vtop{\ialign{##\crcr$\m@th\hfil#1#2\hfil$\crcr
 \noalign{\nointerlineskip}\leftrightarrowfill@#1\crcr}}}
\message{various kinds of dots,}
\let\DOTSI\relax
\let\DOTSB\relax

\newif\ifmath@
{\uccode`7=`\\ \uccode`8=`m \uccode`9=`a \uccode`0=`t \uccode`!=`h
 \uppercase{\gdef\math@#1#2#3#4#5#6\math@{\global\math@false\ifx 7#1\ifx 8#2%
 \ifx 9#3\ifx 0#4\ifx !#5\xdef\meaning@{#6}\global\math@true\fi\fi\fi\fi\fi}}}
\newif\ifmathch@
{\uccode`7=`c \uccode`8=`h \uccode`9=`\"
 \uppercase{\gdef\mathch@#1#2#3#4#5#6\mathch@{\global\mathch@false
  \ifx 7#1\ifx 8#2\ifx 9#5\global\mathch@true\xdef\meaning@{9#6}\fi\fi\fi}}}
\newcount\classnum@
\def\getmathch@#1.#2\getmathch@{\classnum@#1 \divide\classnum@4096
 \ifcase\number\classnum@\or\or\gdef\thedots@{\dotsb@}\or
 \gdef\thedots@{\dotsb@}\fi}
\newif\ifmathbin@
{\uccode`4=`b \uccode`5=`i \uccode`6=`n
 \uppercase{\gdef\mathbin@#1#2#3{\relaxnext@
  \DNii@##1\mathbin@{\ifx\space@\next\global\mathbin@true\fi}%
 \global\mathbin@false\DN@##1\mathbin@{}%
 \ifx 4#1\ifx 5#2\ifx 6#3\DN@{\FN@\nextii@}\fi\fi\fi\next@}}}
\newif\ifmathrel@
{\uccode`4=`r \uccode`5=`e \uccode`6=`l
 \uppercase{\gdef\mathrel@#1#2#3{\relaxnext@
  \DNii@##1\mathrel@{\ifx\space@\next\global\mathrel@true\fi}%
 \global\mathrel@false\DN@##1\mathrel@{}%
 \ifx 4#1\ifx 5#2\ifx 6#3\DN@{\FN@\nextii@}\fi\fi\fi\next@}}}
\newif\ifmacro@
{\uccode`5=`m \uccode`6=`a \uccode`7=`c
 \uppercase{\gdef\macro@#1#2#3#4\macro@{\global\macro@false
  \ifx 5#1\ifx 6#2\ifx 7#3\global\macro@true
  \xdef\meaning@{\macro@@#4\macro@@}\fi\fi\fi}}}
\def\macro@@#1->#2\macro@@{#2}
\newif\ifDOTS@
\newcount\DOTSCASE@
{\uccode`6=`\\ \uccode`7=`D \uccode`8=`O \uccode`9=`T \uccode`0=`S
 \uppercase{\gdef\DOTS@#1#2#3#4#5{\global\DOTS@false\DN@##1\DOTS@{}%
  \ifx 6#1\ifx 7#2\ifx 8#3\ifx 9#4\ifx 0#5\let\next@\DOTS@@\fi\fi\fi\fi\fi
  \next@}}}
{\uccode`3=`B \uccode`4=`I \uccode`5=`X
 \uppercase{\gdef\DOTS@@#1{\relaxnext@
  \DNii@##1\DOTS@{\ifx\space@\next\global\DOTS@true\fi}%
  \DN@{\FN@\nextii@}%
  \ifx 3#1\global\DOTSCASE@\z@\else
  \ifx 4#1\global\DOTSCASE@\@ne\else
  \ifx 5#1\global\DOTSCASE@\tw@\else\DN@##1\DOTS@{}%
  \fi\fi\fi\next@}}}
\newif\ifnot@
{\uccode`5=`\\ \uccode`6=`n \uccode`7=`o \uccode`8=`t
 \uppercase{\gdef\not@#1#2#3#4{\relaxnext@
  \DNii@##1\not@{\ifx\space@\next\global\not@true\fi}%
 \global\not@false\DN@##1\not@{}%
 \ifx 5#1\ifx 6#2\ifx 7#3\ifx 8#4\DN@{\FN@\nextii@}\fi\fi\fi
 \fi\next@}}}
\newif\ifkeybin@
\def\keybin@{\keybin@true
 \ifx\next+\else\ifx\next=\else\ifx\next<\else\ifx\next>\else\ifx\next-\else
 \ifx\next*\else\ifx\next:\else\keybin@false\fi\fi\fi\fi\fi\fi\fi}
\def\dots{\RIfM@\expandafter\mdots@\else\expandafter\tdots@\fi}
\def\tdots@{\unskip\relaxnext@
 \DN@{$\m@th\mathinner{\ldotp\ldotp\ldotp}\,
   \ifx\next,\,$\else\ifx\next.\,$\else\ifx\next;\,$\else\ifx\next:\,$\else
   \ifx\next?\,$\else\ifx\next!\,$\else$ \fi\fi\fi\fi\fi\fi}%
 \ \FN@\next@}
\def\mdots@{\FN@\mdots@@}
\def\mdots@@{\gdef\thedots@{\dotso@}
 \ifx\next\boldkey\gdef\thedots@\boldkey{\boldkeydots@}\else                
 \ifx\next\boldsymbol\gdef\thedots@\boldsymbol{\boldsymboldots@}\else       
 \ifx,\next\gdef\thedots@{\dotsc}
 \else\ifx\not\next\gdef\thedots@{\dotsb@}
 \else\keybin@
 \ifkeybin@\gdef\thedots@{\dotsb@}
 \else\xdef\meaning@{\meaning\next..........}\xdef\meaning@@{\meaning@}
  \expandafter\math@\meaning@\math@
  \ifmath@
   \expandafter\mathch@\meaning@\mathch@
   \ifmathch@\expandafter\getmathch@\meaning@\getmathch@\fi                 
  \else\expandafter\macro@\meaning@@\macro@                                 
  \ifmacro@                                                                
   \expandafter\not@\meaning@\not@\ifnot@\gdef\thedots@{\dotsb@}
  \else\expandafter\DOTS@\meaning@\DOTS@
  \ifDOTS@
   \ifcase\number\DOTSCASE@\gdef\thedots@{\dotsb@}%
    \or\gdef\thedots@{\dotsi}\else\fi                                      
  \else\expandafter\math@\meaning@\math@                                   
  \ifmath@\expandafter\mathbin@\meaning@\mathbin@
  \ifmathbin@\gdef\thedots@{\dotsb@}
  \else\expandafter\mathrel@\meaning@\mathrel@
  \ifmathrel@\gdef\thedots@{\dotsb@}
  \fi\fi\fi\fi\fi\fi\fi\fi\fi\fi\fi\fi
 \thedots@}
\def\plainldots@{\mathinner{\ldotp\ldotp\ldotp}}
\def\plaincdots@{\mathinner{\cdotp\cdotp\cdotp}}
\def\dotsi{\!\plaincdots@}
\let\dotsb@\plaincdots@
\newif\ifextra@
\newif\ifrightdelim@
\def\rightdelim@{\global\rightdelim@true                                    
 \ifx\next)\else                                                            
 \ifx\next]\else
 \ifx\next\rbrack\else
 \ifx\next\}\else
 \ifx\next\rbrace\else
 \ifx\next\rangle\else
 \ifx\next\rceil\else
 \ifx\next\rfloor\else
 \ifx\next\rgroup\else
 \ifx\next\rmoustache\else
 \ifx\next\right\else
 \ifx\next\bigr\else
 \ifx\next\biggr\else
 \ifx\next\Bigr\else                                                        
 \ifx\next\Biggr\else\global\rightdelim@false
 \fi\fi\fi\fi\fi\fi\fi\fi\fi\fi\fi\fi\fi\fi\fi}
\def\extra@{%
 \global\extra@false\rightdelim@\ifrightdelim@\global\extra@true            
 \else\ifx\next$\global\extra@true                                          
 \else\xdef\meaning@{\meaning\next..........}
 \expandafter\macro@\meaning@\macro@\ifmacro@                               
 \expandafter\DOTS@\meaning@\DOTS@
 \ifDOTS@
 \ifnum\DOTSCASE@=\tw@\global\extra@true                                    
 \fi\fi\fi\fi\fi}
\newif\ifbold@
\def\dotso@{\relaxnext@
 \ifbold@
  \let\next\delayed@
  \DNii@{\extra@\plainldots@\ifextra@\,\fi}%
 \else
  \DNii@{\DN@{\extra@\plainldots@\ifextra@\,\fi}\FN@\next@}%
 \fi
 \nextii@}
\def\extrap@#1{%
 \ifx\next,\DN@{#1\,}\else
 \ifx\next;\DN@{#1\,}\else
 \ifx\next.\DN@{#1\,}\else\extra@
 \ifextra@\DN@{#1\,}\else
 \let\next@#1\fi\fi\fi\fi\next@}
\def\ldots{\DN@{\extrap@\plainldots@}%
 \FN@\next@}
\def\cdots{\DN@{\extrap@\plaincdots@}%
 \FN@\next@}

\def\dotsc{\relaxnext@
 \DN@{\ifx\next;\plainldots@\,\else
  \ifx\next.\plainldots@\,\else\extra@\plainldots@
  \ifextra@\,\fi\fi\fi}%
 \FN@\next@}
\def\cdot{\mathchar"2201 }

\message{special superscripts,}
\def\dddot#1{{\mathop{#1}\limits^{\vbox to-1.4\ex@{\kern-\tw@\ex@
 \hbox{\rm...}\vss}}}}
\def\ddddot#1{{\mathop{#1}\limits^{\vbox to-1.4\ex@{\kern-\tw@\ex@
 \hbox{\rm....}\vss}}}}
\def\sphat{^{\mathchoice{}{}%
 {\,\,\botsmash{\hbox{\lower4\ex@\hbox{$\m@th\widehat{\null}$}}}}%
 {\,\botsmash{\hbox{\lower3\ex@\hbox{$\m@th\hat{\null}$}}}}}}

\def\spacute{^{\!\botsmash{\hbox{\lower\@ne ex\hbox{\'{}}}}}}
\def\spgrave{^{\mathchoice{}{}{}{\!}%
 \botsmash{\hbox{\lower\@ne ex\hbox{\`{}}}}}}
\def\spdot{^{\hbox{\raise\ex@\hbox{\rm.}}}}
\def\spddot{^{\hbox{\raise\ex@\hbox{\rm..}}}}
\def\spdddot{^{\hbox{\raise\ex@\hbox{\rm...}}}}
\def\spddddot{^{\hbox{\raise\ex@\hbox{\rm....}}}}
\def\spbreve{^{\!\botsmash{\hbox{\lower4\ex@\hbox{\u{}}}}}}

\message{\string\text,}
\def\textonlyfont@#1#2{\def#1{\RIfM@
 \Err@{Use \string#1\space only in text}\else#2\fi}}
\textonlyfont@\rm\tenrm
\textonlyfont@\it\tenit
\textonlyfont@\sl\tensl
\textonlyfont@\bf\tenbf
\def\oldnos#1{\RIfM@{\mathcode`\,="013B \fam\@ne#1}\else
 \leavevmode\hbox{$\m@th\mathcode`\,="013B \fam\@ne#1$}\fi}
\def\text{\RIfM@\expandafter\text@\else\expandafter\text@@\fi}
\def\text@@#1{\leavevmode\hbox{#1}}
\def\mathhexbox@#1#2#3{\text{$\m@th\mathchar"#1#2#3$}}
\def\dag{{\mathhexbox@279}}
\def\ddag{{\mathhexbox@27A}}
\def\S{{\mathhexbox@278}}
\def\P{{\mathhexbox@27B}}
\newif\iffirstchoice@
\firstchoice@true
\def\text@#1{\mathchoice
 {\hbox{\everymath{\displaystyle}\def\textfonti{\the\textfont\@ne}%
  \def\textfontii{\the\textfont\tw@}\textdef@@ T#1}}
 {\hbox{\firstchoice@false
  \everymath{\textstyle}\def\textfonti{\the\textfont\@ne}%
  \def\textfontii{\the\textfont\tw@}\textdef@@ T#1}}
 {\hbox{\firstchoice@false
  \everymath{\scriptstyle}\def\textfonti{\the\scriptfont\@ne}%
  \def\textfontii{\the\scriptfont\tw@}\textdef@@ S\rm#1}}
 {\hbox{\firstchoice@false
  \everymath{\scriptscriptstyle}\def\textfonti
  {\the\scriptscriptfont\@ne}%
  \def\textfontii{\the\scriptscriptfont\tw@}\textdef@@ s\rm#1}}}
\def\textdef@@#1{\textdef@#1\rm\textdef@#1\bf\textdef@#1\sl\textdef@#1\it}
\def\rmfam{0}
\def\textdef@#1#2{%
 \DN@{\csname\expandafter\eat@\string#2fam\endcsname}%
 \if S#1\edef#2{\the\scriptfont\next@\relax}%
 \else\if s#1\edef#2{\the\scriptscriptfont\next@\relax}%
 \else\edef#2{\the\textfont\next@\relax}\fi\fi}
\scriptfont\itfam\tenit \scriptscriptfont\itfam\tenit
\scriptfont\slfam\tensl \scriptscriptfont\slfam\tensl
\newif\iftopfolded@
\newif\ifbotfolded@
\def\topfoldedtext{\topfolded@true\botfolded@false\foldedtext@}
\def\botfoldedtext{\botfolded@true\topfolded@false\foldedtext@}
\def\foldedtext{\topfolded@false\botfolded@false\foldedtext@}
\Invalid@\foldedwidth
\def\foldedtext@{\relaxnext@
 \DN@{\ifx\next\foldedwidth\let\next@\nextii@\else
  \DN@{\nextii@\foldedwidth{.3\hsize}}\fi\next@}%
 \DNii@\foldedwidth##1##2{\setbox\z@\vbox
  {\normalbaselines\hsize##1\relax
  \tolerance1600 \noindent\ignorespaces##2}\ifbotfolded@\boxz@\else
  \iftopfolded@\vtop{\unvbox\z@}\else\vcenter{\boxz@}\fi\fi}%
 \FN@\next@}
\message{math font commands,}
\def\bold{\RIfM@\expandafter\bold@\else
 \expandafter\nonmatherr@\expandafter\bold\fi}
\def\bold@#1{{\bold@@{#1}}}
\def\bold@@#1{\fam\bffam\relax#1}
\def\slanted{\RIfM@\expandafter\slanted@\else
 \expandafter\nonmatherr@\expandafter\slanted\fi}
\def\slanted@#1{{\slanted@@{#1}}}
\def\slanted@@#1{\fam\slfam\relax#1}
\def\roman{\RIfM@\expandafter\roman@\else
 \expandafter\nonmatherr@\expandafter\roman\fi}
\def\roman@#1{{\roman@@{#1}}}
\def\roman@@#1{\fam\rmfam\relax#1}
\def\italic{\RIfM@\expandafter\italic@\else
 \expandafter\nonmatherr@\expandafter\italic\fi}
\def\italic@#1{{\italic@@{#1}}}
\def\italic@@#1{\fam\itfam\relax#1}
\def\Cal{\RIfM@\expandafter\Cal@\else
 \expandafter\nonmatherr@\expandafter\Cal\fi}
\def\Cal@#1{{\Cal@@{#1}}}
\def\Cal@@#1{\noaccents@\fam\tw@#1}
\mathchardef\Gamma="0000
\mathchardef\Delta="0001
\mathchardef\Theta="0002
\mathchardef\Lambda="0003
\mathchardef\Xi="0004
\mathchardef\Pi="0005
\mathchardef\Sigma="0006
\mathchardef\Upsilon="0007
\mathchardef\Phi="0008
\mathchardef\Psi="0009
\mathchardef\Omega="000A
\mathchardef\varGamma="0100
\mathchardef\varDelta="0101
\mathchardef\varTheta="0102
\mathchardef\varLambda="0103
\mathchardef\varXi="0104
\mathchardef\varPi="0105
\mathchardef\varSigma="0106
\mathchardef\varUpsilon="0107
\mathchardef\varPhi="0108
\mathchardef\varPsi="0109
\mathchardef\varOmega="010A
\let\alloc@@\alloc@
\def\hexnumber@#1{\ifcase#1 0\or 1\or 2\or 3\or 4\or 5\or 6\or 7\or 8\or
 9\or A\or B\or C\or D\or E\or F\fi}
\def\loadmsam{%
 \font@\tenmsa=msam10
 \font@\sevenmsa=msam7
 \font@\fivemsa=msam5
 \alloc@@8\fam\chardef\sixt@@n\msafam
 \textfont\msafam=\tenmsa
 \scriptfont\msafam=\sevenmsa
 \scriptscriptfont\msafam=\fivemsa
 \edef\next{\hexnumber@\msafam}%
 \mathchardef\dabar@"0\next39
 \edef\dashrightarrow{\mathrel{\dabar@\dabar@\mathchar"0\next4B}}%
 \edef\dashleftarrow{\mathrel{\mathchar"0\next4C\dabar@\dabar@}}%
 \let\dasharrow\dashrightarrow
 \edef\ulcorner{\delimiter"4\next70\next70 }%
 \edef\urcorner{\delimiter"5\next71\next71 }%
 \edef\llcorner{\delimiter"4\next78\next78 }%
 \edef\lrcorner{\delimiter"5\next79\next79 }%
 \edef\yen{{\noexpand\mathhexbox@\next55}}%
 \edef\checkmark{{\noexpand\mathhexbox@\next58}}%
 \edef\circledR{{\noexpand\mathhexbox@\next72}}%
 \edef\maltese{{\noexpand\mathhexbox@\next7A}}%
 \global\let\loadmsam\empty}%
\def\loadmsbm{%
 \font@\tenmsb=msbm10 \font@\sevenmsb=msbm7 \font@\fivemsb=msbm5
 \alloc@@8\fam\chardef\sixt@@n\msbfam
 \textfont\msbfam=\tenmsb
 \scriptfont\msbfam=\sevenmsb \scriptscriptfont\msbfam=\fivemsb
 \global\let\loadmsbm\empty
 }
\def\widehat#1{\ifx\undefined\msbfam \DN@{362}%
  \else \setboxz@h{$\m@th#1$}%
    \edef\next@{\ifdim\wdz@>\tw@ em%
        \hexnumber@\msbfam 5B%
      \else 362\fi}\fi
  \mathaccent"0\next@{#1}}
\def\widetilde#1{\ifx\undefined\msbfam \DN@{365}%
  \else \setboxz@h{$\m@th#1$}%
    \edef\next@{\ifdim\wdz@>\tw@ em%
        \hexnumber@\msbfam 5D%
      \else 365\fi}\fi
  \mathaccent"0\next@{#1}}
\message{\string\newsymbol,}
\def\newsymbol#1#2#3#4#5{\define#1{}%
  \count@#2\relax \advance\count@\m@ne 
 \ifcase\count@
   \ifx\undefined\msafam\loadmsam\fi \let\next@\msafam
 \or \ifx\undefined\msbfam\loadmsbm\fi \let\next@\msbfam
 \else  \Err@{\Invalid@@\string\newsymbol}\let\next@\tw@\fi
 \mathchardef#1="#3\hexnumber@\next@#4#5\space}

\def\Bbb{\RIfM@\expandafter\Bbb@\else
 \expandafter\nonmatherr@\expandafter\Bbb\fi}
\def\Bbb@#1{{\Bbb@@{#1}}}
\def\Bbb@@#1{\noaccents@\fam\msbfam\relax#1}
\message{bold Greek and bold symbols,}
\def\loadbold{%
 \font@\tencmmib=cmmib10 \font@\sevencmmib=cmmib7 \font@\fivecmmib=cmmib5
 \skewchar\tencmmib'177 \skewchar\sevencmmib'177 \skewchar\fivecmmib'177
 \alloc@@8\fam\chardef\sixt@@n\cmmibfam
 \textfont\cmmibfam\tencmmib
 \scriptfont\cmmibfam\sevencmmib \scriptscriptfont\cmmibfam\fivecmmib
 \font@\tencmbsy=cmbsy10 \font@\sevencmbsy=cmbsy7 \font@\fivecmbsy=cmbsy5
 \skewchar\tencmbsy'60 \skewchar\sevencmbsy'60 \skewchar\fivecmbsy'60
 \alloc@@8\fam\chardef\sixt@@n\cmbsyfam
 \textfont\cmbsyfam\tencmbsy
 \scriptfont\cmbsyfam\sevencmbsy \scriptscriptfont\cmbsyfam\fivecmbsy
 \let\loadbold\empty
}
\def\boldnotloaded#1{\Err@{\ifcase#1\or First\else Second\fi
       bold symbol font not loaded}}
\def\mathchari@#1#2#3{\ifx\undefined\cmmibfam
    \boldnotloaded@\@ne
  \else\mathchar"#1\hexnumber@\cmmibfam#2#3\space \fi}
\def\mathcharii@#1#2#3{\ifx\undefined\cmbsyfam
    \boldnotloaded\tw@
  \else \mathchar"#1\hexnumber@\cmbsyfam#2#3\space\fi}
\edef\bffam@{\hexnumber@\bffam}
\def\boldkey#1{\ifcat\noexpand#1A%
  \ifx\undefined\cmmibfam \boldnotloaded\@ne
  \else {\fam\cmmibfam#1}\fi
 \else
 \ifx#1!\mathchar"5\bffam@21 \else
 \ifx#1(\mathchar"4\bffam@28 \else\ifx#1)\mathchar"5\bffam@29 \else
 \ifx#1+\mathchar"2\bffam@2B \else\ifx#1:\mathchar"3\bffam@3A \else
 \ifx#1;\mathchar"6\bffam@3B \else\ifx#1=\mathchar"3\bffam@3D \else
 \ifx#1?\mathchar"5\bffam@3F \else\ifx#1[\mathchar"4\bffam@5B \else
 \ifx#1]\mathchar"5\bffam@5D \else
 \ifx#1,\mathchari@63B \else
 \ifx#1-\mathcharii@200 \else
 \ifx#1.\mathchari@03A \else
 \ifx#1/\mathchari@03D \else
 \ifx#1<\mathchari@33C \else
 \ifx#1>\mathchari@33E \else
 \ifx#1*\mathcharii@203 \else
 \ifx#1|\mathcharii@06A \else
 \ifx#10\bold0\else\ifx#11\bold1\else\ifx#12\bold2\else\ifx#13\bold3\else
 \ifx#14\bold4\else\ifx#15\bold5\else\ifx#16\bold6\else\ifx#17\bold7\else
 \ifx#18\bold8\else\ifx#19\bold9\else
  \Err@{\string\boldkey\space can't be used with #1}%
 \fi\fi\fi\fi\fi\fi\fi\fi\fi\fi\fi\fi\fi\fi\fi
 \fi\fi\fi\fi\fi\fi\fi\fi\fi\fi\fi\fi\fi\fi}
\def\boldsymbol#1{%
 \DN@{\Err@{You can't use \string\boldsymbol\space with \string#1}#1}%
 \ifcat\noexpand#1A%
   \let\next@\relax
   \ifx\undefined\cmmibfam \boldnotloaded\@ne
   \else {\fam\cmmibfam#1}\fi
 \else
  \xdef\meaning@{\meaning#1.........}%
  \expandafter\math@\meaning@\math@
  \ifmath@
   \expandafter\mathch@\meaning@\mathch@
   \ifmathch@
    \expandafter\boldsymbol@@\meaning@\boldsymbol@@
   \fi
  \else
   \expandafter\macro@\meaning@\macro@
   \expandafter\delim@\meaning@\delim@
   \ifdelim@
    \expandafter\delim@@\meaning@\delim@@
   \else
    \boldsymbol@{#1}%
   \fi
  \fi
 \fi
 \next@}
\def\mathhexboxii@#1#2{\ifx\undefined\cmbsyfam
    \boldnotloaded\tw@
  \else \mathhexbox@{\hexnumber@\cmbsyfam}{#1}{#2}\fi}
\def\boldsymbol@#1{\let\next@\relax\let\next#1%
 \ifx\next\cdot\mathcharii@201 \else
 \ifx\next\prime{{\null\mathcharii@030 \null}}\else
 \ifx\next\lbrack\mathchar"4\bffam@5B \else
 \ifx\next\rbrack\mathchar"5\bffam@5D \else
 \ifx\next\{\mathcharii@466 \else
 \ifx\next\lbrace\mathcharii@466 \else
 \ifx\next\}\mathcharii@567 \else
 \ifx\next\rbrace\mathcharii@567 \else
 \ifx\next\surd{{\mathcharii@170}}\else
 \ifx\next\S{{\mathhexboxii@78}}\else
 \ifx\next\P{{\mathhexboxii@7B}}\else
 \ifx\next\dag{{\mathhexboxii@79}}\else
 \ifx\next\ddag{{\mathhexboxii@7A}}\else
 \DN@{\Err@{You can't use \string\boldsymbol\space with \string#1}#1}%
 \fi\fi\fi\fi\fi\fi\fi\fi\fi\fi\fi\fi\fi}
\def\boldsymbol@@#1.#2\boldsymbol@@{\classnum@#1 \count@@@\classnum@        
 \divide\classnum@4096 \count@\classnum@                                    
 \multiply\count@4096 \advance\count@@@-\count@ \count@@\count@@@           
 \divide\count@@@\@cclvi \count@\count@@                                    
 \multiply\count@@@\@cclvi \advance\count@@-\count@@@                       
 \divide\count@@@\@cclvi                                                    
 \multiply\classnum@4096 \advance\classnum@\count@@                         
 \ifnum\count@@@=\z@                                                        
  \count@"\bffam@ \multiply\count@\@cclvi
  \advance\classnum@\count@
  \DN@{\mathchar\number\classnum@}%
 \else
  \ifnum\count@@@=\@ne                                                      
   \ifx\undefined\cmmibfam \DN@{\boldnotloaded\@ne}%
   \else \count@\cmmibfam \multiply\count@\@cclvi
     \advance\classnum@\count@
     \DN@{\mathchar\number\classnum@}\fi
  \else
   \ifnum\count@@@=\tw@                                                    
     \ifx\undefined\cmbsyfam
       \DN@{\boldnotloaded\tw@}%
     \else
       \count@\cmbsyfam \multiply\count@\@cclvi
       \advance\classnum@\count@
       \DN@{\mathchar\number\classnum@}%
     \fi
  \fi
 \fi
\fi}
\newif\ifdelim@
\newcount\delimcount@
{\uccode`6=`\\ \uccode`7=`d \uccode`8=`e \uccode`9=`l
 \uppercase{\gdef\delim@#1#2#3#4#5\delim@
  {\delim@false\ifx 6#1\ifx 7#2\ifx 8#3\ifx 9#4\delim@true
   \xdef\meaning@{#5}\fi\fi\fi\fi}}}
\def\delim@@#1"#2#3#4#5#6\delim@@{\if#32%
\let\next@\relax
 \ifx\undefined\cmbsyfam \boldnotloaded\@ne
 \else \mathcharii@#2#4#5\space \fi\fi}
\def\vert{\delimiter"026A30C }
\def\Vert{\delimiter"026B30D }
\let\|\Vert
\def\backslash{\delimiter"026E30F }
\def\boldkeydots@#1{\bold@true\let\next=#1\let\delayed@=#1\mdots@@
 \boldkey#1\bold@false}  
\def\boldsymboldots@#1{\bold@true\let\next#1\let\delayed@#1\mdots@@
 \boldsymbol#1\bold@false}
\message{Euler fonts,}

\def\frak{\mathfont@\frak}

\def\loadmathfont#1{%
   \expandafter\font@\csname ten#1\endcsname=#110
   \expandafter\font@\csname seven#1\endcsname=#17
   \expandafter\font@\csname five#1\endcsname=#15
   \edef\next{\noexpand\alloc@@8\fam\chardef\sixt@@n
     \expandafter\noexpand\csname#1fam\endcsname}%
   \next
   \textfont\csname#1fam\endcsname \csname ten#1\endcsname
   \scriptfont\csname#1fam\endcsname \csname seven#1\endcsname
   \scriptscriptfont\csname#1fam\endcsname \csname five#1\endcsname
   \expandafter\def\csname #1\expandafter\endcsname\expandafter{%
      \expandafter\mathfont@\csname#1\endcsname}%
 \expandafter\gdef\csname load#1\endcsname{}%
}
\def\mathfont@#1{\RIfM@\expandafter\mathfont@@\expandafter#1\else
  \expandafter\nonmatherr@\expandafter#1\fi}
\def\mathfont@@#1#2{{\mathfont@@@#1{#2}}}
\def\mathfont@@@#1#2{\noaccents@
   \fam\csname\expandafter\eat@\string#1fam\endcsname
   \relax#2}
\message{math accents,}
\def\accentclass@{7}
\def\noaccents@{\def\accentclass@{0}}
\def\makeacc@#1#2{\def#1{\mathaccent"\accentclass@#2 }}
\makeacc@\hat{05E}
\makeacc@\check{014}
\makeacc@\tilde{07E}
\makeacc@\acute{013}
\makeacc@\grave{012}
\makeacc@\dot{05F}
\makeacc@\ddot{07F}
\makeacc@\breve{015}
\makeacc@\bar{016}

\newcount\skewcharcount@
\newcount\familycount@
\def\theskewchar@{\familycount@\@ne
 \global\skewcharcount@\the\skewchar\textfont\@ne                           
 \ifnum\fam>\m@ne\ifnum\fam<16
  \global\familycount@\the\fam\relax
  \global\skewcharcount@\the\skewchar\textfont\the\fam\relax\fi\fi          
 \ifnum\skewcharcount@>\m@ne
  \ifnum\skewcharcount@<128
  \multiply\familycount@256
  \global\advance\skewcharcount@\familycount@
  \global\advance\skewcharcount@28672
  \mathchar\skewcharcount@\else
  \global\skewcharcount@\m@ne\fi\else
 \global\skewcharcount@\m@ne\fi}                                            
\newcount\pointcount@
\def\getpoints@#1.#2\getpoints@{\pointcount@#1 }
\newdimen\accentdimen@
\newcount\accentmu@
\def\dimentomu@{\multiply\accentdimen@ 100
 \expandafter\getpoints@\the\accentdimen@\getpoints@
 \multiply\pointcount@18
 \divide\pointcount@\@m
 \global\accentmu@\pointcount@}
\def\Makeacc@#1#2{\def#1{\RIfM@\DN@{\mathaccent@
 {"\accentclass@#2 }}\else\DN@{\nonmatherr@{#1}}\fi\next@}}
\def\unbracefonts@{\let\Cal@\Cal@@\let\roman@\roman@@\let\bold@\bold@@
 \let\slanted@\slanted@@}
\def\mathaccent@#1#2{\ifnum\fam=\m@ne\xdef\thefam@{1}\else
 \xdef\thefam@{\the\fam}\fi                                                 
 \accentdimen@\z@                                                           
 \setboxz@h{\unbracefonts@$\m@th\fam\thefam@\relax#2$}
 \ifdim\accentdimen@=\z@\DN@{\mathaccent#1{#2}}
  \setbox@ne\hbox{\unbracefonts@$\m@th\fam\thefam@\relax#2\theskewchar@$}
  \setbox\tw@\hbox{$\m@th\ifnum\skewcharcount@=\m@ne\else
   \mathchar\skewcharcount@\fi$}
  \global\accentdimen@\wd@ne\global\advance\accentdimen@-\wdz@
  \global\advance\accentdimen@-\wd\tw@                                     
  \global\multiply\accentdimen@\tw@
  \dimentomu@\global\advance\accentmu@\@ne                                 
 \else\DN@{{\mathaccent#1{#2\mkern\accentmu@ mu}%
    \mkern-\accentmu@ mu}{}}\fi                                             
 \next@}\Makeacc@\Hat{05E}
\Makeacc@\Check{014}
\Makeacc@\Tilde{07E}
\Makeacc@\Acute{013}
\Makeacc@\Grave{012}
\Makeacc@\Dot{05F}
\Makeacc@\Ddot{07F}
\Makeacc@\Breve{015}
\Makeacc@\Bar{016}
\def\Vec{\RIfM@\DN@{\mathaccent@{"017E }}\else
 \DN@{\nonmatherr@\Vec}\fi\next@}
\def\accentedsymbol#1#2{\csname newbox\expandafter\endcsname
  \csname\expandafter\eat@\string#1@box\endcsname
 \expandafter\setbox\csname\expandafter\eat@
  \string#1@box\endcsname\hbox{$\m@th#2$}\define
  #1{\copy\csname\expandafter\eat@\string#1@box\endcsname{}}}
\message{roots,}
\def\sqrt#1{\radical"270370 {#1}}
\let\underline@\underline
\let\overline@\overline
\def\underline#1{\underline@{#1}}
\def\overline#1{\overline@{#1}}
\Invalid@\leftroot
\Invalid@\uproot
\newcount\uproot@
\newcount\leftroot@
\def\root{\relaxnext@
  \DN@{\ifx\next\uproot\let\next@\nextii@\else
   \ifx\next\leftroot\let\next@\nextiii@\else
   \let\next@\plainroot@\fi\fi\next@}%
  \DNii@\uproot##1{\uproot@##1\relax\FN@\nextiv@}%
  \def\nextiv@{\ifx\next\space@\DN@. {\FN@\nextv@}\else
   \DN@.{\FN@\nextv@}\fi\next@.}%
  \def\nextv@{\ifx\next\leftroot\let\next@\nextvi@\else
   \let\next@\plainroot@\fi\next@}%
  \def\nextvi@\leftroot##1{\leftroot@##1\relax\plainroot@}%
   \def\nextiii@\leftroot##1{\leftroot@##1\relax\FN@\nextvii@}%
  \def\nextvii@{\ifx\next\space@
   \DN@. {\FN@\nextviii@}\else
   \DN@.{\FN@\nextviii@}\fi\next@.}%
  \def\nextviii@{\ifx\next\uproot\let\next@\nextix@\else
   \let\next@\plainroot@\fi\next@}%
  \def\nextix@\uproot##1{\uproot@##1\relax\plainroot@}%
  \bgroup\uproot@\z@\leftroot@\z@\FN@\next@}
\def\plainroot@#1\of#2{\setbox\rootbox\hbox{$\m@th\scriptscriptstyle{#1}$}%
 \mathchoice{\r@@t\displaystyle{#2}}{\r@@t\textstyle{#2}}
 {\r@@t\scriptstyle{#2}}{\r@@t\scriptscriptstyle{#2}}\egroup}
\def\r@@t#1#2{\setboxz@h{$\m@th#1\sqrt{#2}$}%
 \dimen@\ht\z@\advance\dimen@-\dp\z@
 \setbox@ne\hbox{$\m@th#1\mskip\uproot@ mu$}\advance\dimen@ 1.667\wd@ne
 \mkern-\leftroot@ mu\mkern5mu\raise.6\dimen@\copy\rootbox
 \mkern-10mu\mkern\leftroot@ mu\boxz@}
\def\boxed#1{\setboxz@h{$\m@th\displaystyle{#1}$}\dimen@.4\ex@
 \advance\dimen@3\ex@\advance\dimen@\dp\z@
 \hbox{\lower\dimen@\hbox{%
 \vbox{\hrule height.4\ex@
 \hbox{\vrule width.4\ex@\hskip3\ex@\vbox{\vskip3\ex@\boxz@\vskip3\ex@}%
 \hskip3\ex@\vrule width.4\ex@}\hrule height.4\ex@}%
 }}}
\message{commutative diagrams,}
\let\ampersand@\relax
\newdimen\minaw@
\minaw@11.11128\ex@
\newdimen\minCDaw@
\minCDaw@2.5pc
\def\minCDarrowwidth#1{\RIfMIfI@\onlydmatherr@\minCDarrowwidth
 \else\minCDaw@#1\relax\fi\else\onlydmatherr@\minCDarrowwidth\fi}
\newif\ifCD@
\def\CD{\bgroup\vspace@\relax\let\ampersand@&\iffalse}\fi
 \CD@true\vcenter\bgroup\Let@\tabskip\z@skip\baselineskip20\ex@
 \lineskip3\ex@\lineskiplimit3\ex@\halign\bgroup
 &\hfill$\m@th##$\hfill\crcr}
\def\endCD{\crcr\egroup\egroup\egroup}
\newdimen\bigaw@
\atdef@>#1>#2>{\ampersand@                                                  
 \setboxz@h{$\m@th\ssize\;{#1}\;\;$}
 \setbox@ne\hbox{$\m@th\ssize\;{#2}\;\;$}
 \setbox\tw@\hbox{$\m@th#2$}
 \ifCD@\global\bigaw@\minCDaw@\else\global\bigaw@\minaw@\fi                 
 \ifdim\wdz@>\bigaw@\global\bigaw@\wdz@\fi
 \ifdim\wd@ne>\bigaw@\global\bigaw@\wd@ne\fi                                
 \ifCD@\enskip\fi                                                           
 \ifdim\wd\tw@>\z@
  \mathrel{\mathop{\hbox to\bigaw@{\rightarrowfill@\displaystyle}}%
    \limits^{#1}_{#2}}
 \else\mathrel{\mathop{\hbox to\bigaw@{\rightarrowfill@\displaystyle}}%
    \limits^{#1}}\fi                                                        
 \ifCD@\enskip\fi                                                          
 \ampersand@}                                                              
\atdef@<#1<#2<{\ampersand@\setboxz@h{$\m@th\ssize\;\;{#1}\;$}%
 \setbox@ne\hbox{$\m@th\ssize\;\;{#2}\;$}\setbox\tw@\hbox{$\m@th#2$}%
 \ifCD@\global\bigaw@\minCDaw@\else\global\bigaw@\minaw@\fi
 \ifdim\wdz@>\bigaw@\global\bigaw@\wdz@\fi
 \ifdim\wd@ne>\bigaw@\global\bigaw@\wd@ne\fi
 \ifCD@\enskip\fi
 \ifdim\wd\tw@>\z@
  \mathrel{\mathop{\hbox to\bigaw@{\leftarrowfill@\displaystyle}}%
       \limits^{#1}_{#2}}\else
  \mathrel{\mathop{\hbox to\bigaw@{\leftarrowfill@\displaystyle}}%
       \limits^{#1}}\fi
 \ifCD@\enskip\fi\ampersand@}
\begingroup
 \catcode`\~=\active \lccode`\~=`\@
 \lowercase{%
  \global\atdef@)#1)#2){~>#1>#2>}
  \global\atdef@(#1(#2({~<#1<#2<}}
\endgroup
\atdef@ A#1A#2A{\llap{$\m@th\vcenter{\hbox
 {$\ssize#1$}}$}\Big\uparrow\rlap{$\m@th\vcenter{\hbox{$\ssize#2$}}$}&&}
\atdef@ V#1V#2V{\llap{$\m@th\vcenter{\hbox
 {$\ssize#1$}}$}\Big\downarrow\rlap{$\m@th\vcenter{\hbox{$\ssize#2$}}$}&&}
\atdef@={&\enskip\mathrel
 {\vbox{\hrule width\minCDaw@\vskip3\ex@\hrule width
 \minCDaw@}}\enskip&}
\atdef@|{\Big\Vert&&}
\atdef@\vert{\Big\Vert&&}
\def\pretend#1\haswidth#2{\setboxz@h{$\m@th\scriptstyle{#2}$}\hbox
 to\wdz@{\hfill$\m@th\scriptstyle{#1}$\hfill}}
\message{poor man's bold,}
\def\pmb{\RIfM@\expandafter\mathpalette\expandafter\pmb@\else
 \expandafter\pmb@@\fi}
\def\pmb@@#1{\leavevmode\setboxz@h{#1}%
   \dimen@-\wdz@
   \kern-.5\ex@\copy\z@
   \kern\dimen@\kern.25\ex@\raise.4\ex@\copy\z@
   \kern\dimen@\kern.25\ex@\box\z@
}
\def\binrel@@#1{\ifdim\wd2<\z@\mathbin{#1}\else\ifdim\wd\tw@>\z@
 \mathrel{#1}\else{#1}\fi\fi}
\newdimen\pmbraise@
\def\pmb@#1#2{\setbox\thr@@\hbox{$\m@th#1{#2}$}%
 \setbox4\hbox{$\m@th#1\mkern.5mu$}\pmbraise@\wd4\relax
 \binrel@{#2}%
 \dimen@-\wd\thr@@
   \binrel@@{%
   \mkern-.8mu\copy\thr@@
   \kern\dimen@\mkern.4mu\raise\pmbraise@\copy\thr@@
   \kern\dimen@\mkern.4mu\box\thr@@
}}
\def\documentstyle#1{\W@{}\input #1.sty\relax}
\message{syntax check,}
\font\dummyft@=dummy
\fontdimen1 \dummyft@=\z@
\fontdimen2 \dummyft@=\z@
\fontdimen3 \dummyft@=\z@
\fontdimen4 \dummyft@=\z@
\fontdimen5 \dummyft@=\z@
\fontdimen6 \dummyft@=\z@
\fontdimen7 \dummyft@=\z@
\fontdimen8 \dummyft@=\z@
\fontdimen9 \dummyft@=\z@
\fontdimen10 \dummyft@=\z@
\fontdimen11 \dummyft@=\z@
\fontdimen12 \dummyft@=\z@
\fontdimen13 \dummyft@=\z@
\fontdimen14 \dummyft@=\z@
\fontdimen15 \dummyft@=\z@
\fontdimen16 \dummyft@=\z@
\fontdimen17 \dummyft@=\z@
\fontdimen18 \dummyft@=\z@
\fontdimen19 \dummyft@=\z@
\fontdimen20 \dummyft@=\z@
\fontdimen21 \dummyft@=\z@
\fontdimen22 \dummyft@=\z@
\def\fontlist@{\\{\tenrm}\\{\sevenrm}\\{\fiverm}\\{\teni}\\{\seveni}%
 \\{\fivei}\\{\tensy}\\{\sevensy}\\{\fivesy}\\{\tenex}\\{\tenbf}\\{\sevenbf}%
 \\{\fivebf}\\{\tensl}\\{\tenit}}
\def\font@#1=#2 {\rightappend@#1\to\fontlist@\font#1=#2 }
\def\dodummy@{{\def\\##1{\global\let##1\dummyft@}\fontlist@}}
\def\nopages@{\output{\setbox\z@\box\@cclv \deadcycles\z@}%
 \alloc@5\toks\toksdef\@cclvi\output}
\let\galleys\nopages@
\newif\ifsyntax@
\newcount\countxviii@
\def\syntax{\syntax@true\dodummy@\countxviii@\count18
 \loop\ifnum\countxviii@>\m@ne\textfont\countxviii@=\dummyft@
 \scriptfont\countxviii@=\dummyft@\scriptscriptfont\countxviii@=\dummyft@
 \advance\countxviii@\m@ne\repeat                                           
 \dummyft@\tracinglostchars\z@\nopages@\frenchspacing\hbadness\@M}
\def\first@#1#2\end{#1}
\def\printoptions{\W@{Do you want S(yntax check),
  G(alleys) or P(ages)?}%
 \message{Type S, G or P, followed by <return>: }%
 \begingroup 
 \endlinechar\m@ne 
 \read\m@ne to\ans@
 \edef\ans@{\uppercase{\def\noexpand\ans@{%
   \expandafter\first@\ans@ P\end}}}%
 \expandafter\endgroup\ans@
 \if\ans@ P
 \else \if\ans@ S\syntax
 \else \if\ans@ G\galleys
 \else\message{? Unknown option: \ans@; using the `pages' option.}%
 \fi\fi\fi}
\def\alloc@#1#2#3#4#5{\global\advance\count1#1by\@ne
 \ch@ck#1#4#2\allocationnumber=\count1#1
 \global#3#5=\allocationnumber
 \ifalloc@\wlog{\string#5=\string#2\the\allocationnumber}\fi}
\def\document{\def\alloclist@{}\def\fontlist@{}}
\let\enddocument\bye

\let\proclaim\undefined
\let\footnote\undefined
\let\=\undefined
\let\>\undefined

\catcode`\@=\active
\message{... finished}

\expandafter\ifx\csname mathdefs.tex\endcsname\relax
  \expandafter\gdef\csname mathdefs.tex\endcsname{}
\else \message{Hey!  Apparently you were trying to
  \string\input{mathdefs.tex} twice.   This does not make sense.} 
\errmessage{Please edit your file (probably \jobname.tex) and remove
any duplicate ``\string\input'' lines}\endinput\fi




\catcode`\X=12\catcode`\@=11

\def\n@wcount{\alloc@0\count\countdef\insc@unt}
\def\n@wwrite{\alloc@7\write\chardef\sixt@@n}
\def\n@wread{\alloc@6\read\chardef\sixt@@n}
\def\r@s@t{\relax}\def\v@idline{\par}\def\@mputate#1/{#1}
\def\l@c@l#1X{\firstpart.#1}\def\gl@b@l#1X{#1}\def\t@d@l#1X{{}}

\def\crossrefs#1{\ifx\all#1\let\tr@ce=\all\else\def\tr@ce{#1,}\fi
   \n@wwrite\cit@tionsout\openout\cit@tionsout=\jobname.cit 
   \write\cit@tionsout{\tr@ce}\expandafter\setfl@gs\tr@ce,}
\def\setfl@gs#1,{\def\@{#1}\ifx\@\empty\let\next=\relax
   \else\let\next=\setfl@gs\expandafter\xdef
   \csname#1tr@cetrue\endcsname{}\fi\next}
\def\m@ketag#1#2{\expandafter\n@wcount\csname#2tagno\endcsname
     \csname#2tagno\endcsname=0\let\tail=\all\xdef\all{\tail#2,}
   \ifx#1\l@c@l\let\tail=\r@s@t\xdef\r@s@t{\csname#2tagno\endcsname=0\tail}\fi
   \expandafter\gdef\csname#2cite\endcsname##1{\expandafter
     \ifx\csname#2tag##1\endcsname\relax?\else\csname#2tag##1\endcsname\fi
     \expandafter\ifx\csname#2tr@cetrue\endcsname\relax\else
     \write\cit@tionsout{#2tag ##1 cited on page \folio.}\fi}
   \expandafter\gdef\csname#2page\endcsname##1{\expandafter
     \ifx\csname#2page##1\endcsname\relax?\else\csname#2page##1\endcsname\fi
     \expandafter\ifx\csname#2tr@cetrue\endcsname\relax\else
     \write\cit@tionsout{#2tag ##1 cited on page \folio.}\fi}
   \expandafter\gdef\csname#2tag\endcsname##1{\expandafter
      \ifx\csname#2check##1\endcsname\relax
      \expandafter\xdef\csname#2check##1\endcsname{}%
      \else\immediate\write16{Warning: #2tag ##1 used more than once.}\fi
      \multit@g{#1}{#2}##1/X%
      \write\t@gsout{#2tag ##1 assigned number \csname#2tag##1\endcsname\space
      on page \number\count0.}%
   \csname#2tag##1\endcsname}}

\def\multit@g#1#2#3/#4X{\def\t@mp{#4}\ifx\t@mp\empty%
      \global\advance\csname#2tagno\endcsname by 1 
      \expandafter\xdef\csname#2tag#3\endcsname
      {#1\number\csname#2tagno\endcsnameX}%
   \else\expandafter\ifx\csname#2last#3\endcsname\relax
      \expandafter\n@wcount\csname#2last#3\endcsname
      \global\advance\csname#2tagno\endcsname by 1 
      \expandafter\xdef\csname#2tag#3\endcsname
      {#1\number\csname#2tagno\endcsnameX}
      \write\t@gsout{#2tag #3 assigned number \csname#2tag#3\endcsname\space
      on page \number\count0.}\fi
   \global\advance\csname#2last#3\endcsname by 1
   \def\t@mp{\expandafter\xdef\csname#2tag#3/}%
   \expandafter\t@mp\@mputate#4\endcsname
   {\csname#2tag#3\endcsname\lastpart{\csname#2last#3\endcsname}}\fi}
\def\t@gs#1{\def\all{}\m@ketag#1e\m@ketag#1s\m@ketag\t@d@l p
\let\realscite\scite
\let\realstag\stag
   \m@ketag\gl@b@l r \n@wread\t@gsin
   \openin\t@gsin=\jobname.tgs \re@der \closein\t@gsin
   \n@wwrite\t@gsout\openout\t@gsout=\jobname.tgs }
\outer\def\localtags{\t@gs\l@c@l}
\outer\def\globaltags{\t@gs\gl@b@l}
\outer\def\newlocaltag#1{\m@ketag\l@c@l{#1}}
\outer\def\newglobaltag#1{\m@ketag\gl@b@l{#1}}

\newif\ifpr@ 
\def\m@kecs #1tag #2 assigned number #3 on page #4.%
   {\expandafter\gdef\csname#1tag#2\endcsname{#3}
   \expandafter\gdef\csname#1page#2\endcsname{#4}
   \ifpr@\expandafter\xdef\csname#1check#2\endcsname{}\fi}
\def\re@der{\ifeof\t@gsin\let\next=\relax\else
   \read\t@gsin to\t@gline\ifx\t@gline\v@idline\else
   \expandafter\m@kecs \t@gline\fi\let \next=\re@der\fi\next}
\def\pretags#1{\pr@true\pret@gs#1,,}
\def\pret@gs#1,{\def\@{#1}\ifx\@\empty\let\n@xtfile=\relax
   \else\let\n@xtfile=\pret@gs \openin\t@gsin=#1.tgs \message{#1} \re@der 
   \closein\t@gsin\fi \n@xtfile}

\newcount\sectno\sectno=0\newcount\subsectno\subsectno=0
\newif\ifultr@local \def\ultralocal{\ultr@localtrue}
\def\firstpart{\number\sectno}
\def\lastpart#1{\ifcase#1 \or a\or b\or c\or d\or e\or f\or g\or h\or 
   i\or k\or l\or m\or n\or o\or p\or q\or r\or s\or t\or u\or v\or w\or 
   x\or y\or z \fi}

\def\resetall{\global\advance\sectno by 1\subsectno=0
   \gdef\firstpart{\number\sectno}\r@s@t}
\def\resetsub{\global\advance\subsectno by 1
   \gdef\firstpart{\number\sectno.\number\subsectno}\r@s@t}
\def\newsection#1\par{\resetall\vskip0pt plus.3\vsize\penalty-250
   \vskip0pt plus-.3\vsize\bigskip\bigskip
   \message{#1}\leftline{\bf#1}\nobreak\bigskip}
\def\subsection#1\par{\ifultr@local\resetsub\fi
   \vskip0pt plus.2\vsize\penalty-250\vskip0pt plus-.2\vsize
   \bigskip\smallskip\message{#1}\leftline{\bf#1}\nobreak\medskip}


\newdimen\marginshift

\newdimen\margindelta
\newdimen\marginmax
\newdimen\marginmin

\def\margininit{       
\marginmax=3 true cm                  
				      
\margindelta=0.1 true cm              
\marginmin=0.1true cm                 
\marginshift=\marginmin
}    

\def\t@gsjj#1,{\def\@{#1}\ifx\@\empty\let\next=\relax\else\let\next=\t@gsjj
   \def\@@{p}\ifx\@\@@\else
   \expandafter\gdef\csname#1cite\endcsname##1{\citejj{##1}}
   \expandafter\gdef\csname#1page\endcsname##1{?}
   \expandafter\gdef\csname#1tag\endcsname##1{\tagjj{##1}}\fi\fi\next}
\newif\ifshowstuffinmargin
\showstuffinmarginfalse
\def\jjtags{\ifx\shlhetal\relax 
  \else
\ifx\shlhetal\undefinedcontrolseq
\else
\showstuffinmargintrue
\ifx\all\relax\else\expandafter\t@gsjj\all,\fi\fi \fi
}

\def\tagjj#1{\realstag{#1}\mginpar{\zeigen{#1}}}
\def\citejj#1{\rechnen{#1}\mginpar{\zeigen{#1}}}     

\def\rechnen#1{\expandafter\ifx\csname stag#1\endcsname\relax ??\else
                           \csname stag#1\endcsname\fi}

\newdimen\theight

\def\marginfont{\sevenrm}

\def\trymarginbox#1{\setbox0=\hbox{\marginfont\hskip\marginshift #1}%
		\global\marginshift\wd0 
		\global\advance\marginshift\margindelta}

\def \mginpar#1{%
\ifvmode\setbox0\hbox to \hsize{\hfill\rlap{\marginfont\quad#1}}%
\ht0 0cm
\dp0 0cm
\box0\vskip-\baselineskip
\else 
             \vadjust{\trymarginbox{#1}%
		\ifdim\marginshift>\marginmax \global\marginshift\marginmin
			\trymarginbox{#1}%
                \fi
             \theight=\ht0
             \advance\theight by \dp0    \advance\theight by \lineskip
             \kern -\theight \vbox to \theight{\rightline{\rlap{\box0}}%
\vss}}\fi}


\def\t@gsoff#1,{\def\@{#1}\ifx\@\empty\let\next=\relax\else\let\next=\t@gsoff
   \def\@@{p}\ifx\@\@@\else
   \expandafter\gdef\csname#1cite\endcsname##1{\zeigen{##1}}
   \expandafter\gdef\csname#1page\endcsname##1{?}
   \expandafter\gdef\csname#1tag\endcsname##1{\zeigen{##1}}\fi\fi\next}
\def\verbatimtags{\showstuffinmarginfalse
\ifx\all\relax\else\expandafter\t@gsoff\all,\fi}
\def\zeigen#1{\hbox{$\langle$}#1\hbox{$\rangle$}}
\def\margincite#1{\ifshowstuffinmargin\mginpar{\rechnen{#1}}\fi}
\def\margintag#1{\ifshowstuffinmargin\mginpar{\zeigen{#1}}\fi}

\def\(#1){\edef\dot@g{\ifmmode\ifinner(\hbox{\noexpand\etag{#1}})
   \else\noexpand\eqno(\hbox{\noexpand\etag{#1}})\fi
   \else(\noexpand\ecite{#1})\fi}\dot@g}

\newif\ifbr@ck
\def\eat#1{}
\def\[#1]{\br@cktrue[\br@cket#1'X]}
\def\br@cket#1'#2X{\def\temp{#2}\ifx\temp\empty\let\next\eat
   \else\let\next\br@cket\fi
   \ifbr@ck\br@ckfalse\br@ck@t#1,X\else\br@cktrue#1\fi\next#2X}
\def\br@ck@t#1,#2X{\def\temp{#2}\ifx\temp\empty\let\neext\eat
   \else\let\neext\br@ck@t\def\temp{,}\fi
   \def\teemp{#1}\ifx\teemp\empty\else\rcite{#1}\fi\temp\neext#2X}
\def\resetbr@cket{\gdef\[##1]{[\rtag{##1}]}}
\def\references{\resetbr@cket\newsection References\par}

\newtoks\symb@ls\newtoks\s@mb@ls\newtoks\p@gelist\n@wcount\ftn@mber
    \ftn@mber=1\newif\ifftn@mbers\ftn@mbersfalse\newif\ifbyp@ge\byp@gefalse
\def\defm@rk{\ifftn@mbers\n@mberm@rk\else\symb@lm@rk\fi}
\def\n@mberm@rk{\xdef\m@rk{{\the\ftn@mber}}%
    \global\advance\ftn@mber by 1 }
\def\rot@te#1{\let\temp=#1\global#1=\expandafter\r@t@te\the\temp,X}
\def\r@t@te#1,#2X{{#2#1}\xdef\m@rk{{#1}}}
\def\b@@st#1{{$^{#1}$}}\def\str@p#1{#1}
\def\symb@lm@rk{\ifbyp@ge\rot@te\p@gelist\ifnum\expandafter\str@p\m@rk=1 
    \s@mb@ls=\symb@ls\fi\write\f@nsout{\number\count0}\fi \rot@te\s@mb@ls}
\def\byp@ge{\byp@getrue\n@wwrite\f@nsin\openin\f@nsin=\jobname.fns 
    \n@wcount\currentp@ge\currentp@ge=0\p@gelist={0}
    \re@dfns\closein\f@nsin\rot@te\p@gelist
    \n@wread\f@nsout\openout\f@nsout=\jobname.fns }
\def\m@kelist#1X#2{{#1,#2}}
\def\re@dfns{\ifeof\f@nsin\let\next=\relax\else\read\f@nsin to \f@nline
    \ifx\f@nline\v@idline\else\let\t@mplist=\p@gelist
    \ifnum\currentp@ge=\f@nline
    \global\p@gelist=\expandafter\m@kelist\the\t@mplistX0
    \else\currentp@ge=\f@nline
    \global\p@gelist=\expandafter\m@kelist\the\t@mplistX1\fi\fi
    \let\next=\re@dfns\fi\next}
\def\symbols#1{\symb@ls={#1}\s@mb@ls=\symb@ls} 
\def\bigsymbol{\textstyle}
\symbols{\bigsymbol\ast,\dagger,\ddagger,\sharp,\flat,\natural,\star}
\def\ftnumbers{\ftn@mberstrue} \def\ftsymbols{\ftn@mbersfalse}
\def\paginal{\byp@ge} \def\resetftnumbers{\ftn@mber=1}
\def\ftnote#1{\defm@rk\expandafter\expandafter\expandafter\footnote
    \expandafter\b@@st\m@rk{#1}}

\long\def\jump#1\endjump{}
\def\ssum{\mathop{\lower .1em\hbox{$\textstyle\Sigma$}}\nolimits}

\def\qed{\nobreak\kern 1em \vrule height .5em width .5em depth 0em}
\def\newneq{\hbox{\rlap{\hbox to 1\wd9{\hss$=$\hss}}\raise .1em 
   \hbox to 1\wd9{\hss$\scriptscriptstyle/$\hss}}}
\def\subsetne{\setbox9 = \hbox{$\subset$}\mathrel{\hbox{\rlap
   {\lower .4em \newneq}\raise .13em \hbox{$\subset$}}}}
\def\supsetne{\setbox9 = \hbox{$\subset$}\mathrel{\hbox{\rlap
   {\lower .4em \newneq}\raise .13em \hbox{$\supset$}}}}

\def\vbar{\mathchoice{\vrule height6.3ptdepth-.5ptwidth.8pt\kern-.8pt}
   {\vrule height6.3ptdepth-.5ptwidth.8pt\kern-.8pt}
   {\vrule height4.1ptdepth-.35ptwidth.6pt\kern-.6pt}
   {\vrule height3.1ptdepth-.25ptwidth.5pt\kern-.5pt}}
\def\f@dge{\mathchoice{}{}{\mkern.5mu}{\mkern.8mu}}
\def\b@c#1#2{{\rm \mkern#2mu\vbar\mkern-#2mu#1}}
\def\b@b#1{{\rm I\mkern-3.5mu #1}}
\def\b@a#1#2{{\rm #1\mkern-#2mu\f@dge #1}}
\def\bb#1{{\count4=`#1 \advance\count4by-64 \ifcase\count4\or\b@a A{11.5}\or
   \b@b B\or\b@c C{5}\or\b@b D\or\b@b E\or\b@b F \or\b@c G{5}\or\b@b H\or
   \b@b I\or\b@c J{3}\or\b@b K\or\b@b L \or\b@b M\or\b@b N\or\b@c O{5} \or
   \b@b P\or\b@c Q{5}\or\b@b R\or\b@a S{8}\or\b@a T{10.5}\or\b@c U{5}\or
   \b@a V{12}\or\b@a W{16.5}\or\b@a X{11}\or\b@a Y{11.7}\or\b@a Z{7.5}\fi}}

\catcode`\X=11 \catcode`\@=12


\expandafter\ifx\csname citeadd.tex\endcsname\relax
\expandafter\gdef\csname citeadd.tex\endcsname{}
\else \message{Hey!  Apparently you were trying to
\string\input{citeadd.tex} twice.   This does not make sense.} 
\errmessage{Please edit your file (probably \jobname.tex) and remove
any duplicate ``\string\input'' lines}\endinput\fi

\sectno=-1   
\localtags
\jjtags
\NoBlackBoxes
\define\mr{\medskip\roster}
\define\sn{\smallskip\noindent}
\define\mn{\medskip\noindent}
\define\bn{\bigskip\noindent}
\define\ub{\underbar}
\define\wilog{\text{without loss of generality}}
\define\ermn{\endroster\medskip\noindent}

\define \nl{\newline}
\magnification=\magstep 1
\documentstyle{amsppt}

{    
\catcode`@11

\ifx\alicetwothousandloaded@\relax
  \endinput\else\global\let\alicetwothousandloaded@\relax\fi

\gdef\subjclass{\let\savedef@\subjclass
 \def\subjclass##1\endsubjclass{\let\subjclass\savedef@
   \toks@{\def\usualspace{{\rm\enspace}}\eightpoint}%
   \toks@@{##1\unskip.}%
   \edef\thesubjclass@{\the\toks@
     \frills@{{\noexpand\rm2000 {\noexpand\it Mathematics Subject
       Classification}.\noexpand\enspace}}%
     \the\toks@@}}%
  \nofrillscheck\subjclass}
} 


\expandafter\ifx\csname alice2jlem.tex\endcsname\relax
  \expandafter\xdef\csname alice2jlem.tex\endcsname{\the\catcode`@}
\else \message{Hey!  Apparently you were trying to
\string\input{alice2jlem.tex}  twice.   This does not make sense.}
\errmessage{Please edit your file (probably \jobname.tex) and remove
any duplicate ``\string\input'' lines}\endinput\fi

\expandafter\ifx\csname bib4plain.tex\endcsname\relax
  \expandafter\gdef\csname bib4plain.tex\endcsname{}
\else \message{Hey!  Apparently you were trying to \string\input
  bib4plain.tex twice.   This does not make sense.}
\errmessage{Please edit your file (probably \jobname.tex) and remove
any duplicate ``\string\input'' lines}\endinput\fi

\def\renewcommand{\newcommand}	       
\edef\cite{\the\catcode`@}%
\catcode`@ = 11
\let\@oldatcatcode = \cite
\chardef\@letter = 11
\chardef\@other = 12
%
%
%
%
\def\@innerdef#1#2{\edef#1{\expandafter\noexpand\csname #2\endcsname}}%
%
%
\@innerdef\@innernewcount{newcount}%
\@innerdef\@innernewdimen{newdimen}%
\@innerdef\@innernewif{newif}%
\@innerdef\@innernewwrite{newwrite}%
%
%
%
\def\@gobble#1{}%
%
%
%
\ifx\inputlineno\@undefined
   \let\@linenumber = \empty 
\else
   \def\@linenumber{\the\inputlineno:\space}%
\fi
%
%
%
\def\@futurenonspacelet#1{\def\cs{#1}%
   \afterassignment\@stepone\let\@nexttoken=
}%
\begingroup 
\def\\{\global\let\@stoken= }%
\\ 
\endgroup
\def\@stepone{\expandafter\futurelet\cs\@steptwo}%
\def\@steptwo{\expandafter\ifx\cs\@stoken\let\@@next=\@stepthree
   \else\let\@@next=\@nexttoken\fi \@@next}%
\def\@stepthree{\afterassignment\@stepone\let\@@next= }%
%
%
%
\def\@getoptionalarg#1{%
   \let\@optionaltemp = #1%
   \let\@optionalnext = \relax
   \@futurenonspacelet\@optionalnext\@bracketcheck
}%
%
%
\def\@bracketcheck{%
   \ifx [\@optionalnext
      \expandafter\@@getoptionalarg
   \else
      \let\@optionalarg = \empty
      \expandafter\@optionaltemp
   \fi
}%
\def\@@getoptionalarg[#1]{%
   \def\@optionalarg{#1}%
   \@optionaltemp
}%
%
%
%
\def\@nnil{\@nil}%
\def\@fornoop#1\@@#2#3{}%
\def\@for#1:=#2\do#3{%
   \edef\@fortmp{#2}%
   \ifx\@fortmp\empty \else
      \expandafter\@forloop#2,\@nil,\@nil\@@#1{#3}%
   \fi
}%
\def\@forloop#1,#2,#3\@@#4#5{\def#4{#1}\ifx #4\@nnil \else
       #5\def#4{#2}\ifx #4\@nnil \else#5\@iforloop #3\@@#4{#5}\fi\fi
}%
\def\@iforloop#1,#2\@@#3#4{\def#3{#1}\ifx #3\@nnil
       \let\@nextwhile=\@fornoop \else
      #4\relax\let\@nextwhile=\@iforloop\fi\@nextwhile#2\@@#3{#4}%
}%
%
%
%
\@innernewif\if@fileexists
\def\@testfileexistence{\@getoptionalarg\@finishtestfileexistence}%
\def\@finishtestfileexistence#1{%
   \begingroup
      \def\extension{#1}%
      \immediate\openin0 =
         \ifx\@optionalarg\empty\jobname\else\@optionalarg\fi
         \ifx\extension\empty \else .#1\fi
         \space
      \ifeof 0
         \global\@fileexistsfalse
      \else
         \global\@fileexiststrue
      \fi
      \immediate\closein0
   \endgroup
}%
%
%
%
%
\def\bibliographystyle#1{%
   \@readauxfile
   \@writeaux{\string\bibstyle{#1}}%
}%
\let\bibstyle = \@gobble
%
%
\let\bblfilebasename = \jobname
\def\bibliography#1{%
   \@readauxfile
   \@writeaux{\string\bibdata{#1}}%
   \@testfileexistence[\bblfilebasename]{bbl}%
   \if@fileexists
      \nobreak
      \@readbblfile
   \fi
}%
\let\bibdata = \@gobble
%
%
\def\nocite#1{%
   \@readauxfile
   \@writeaux{\string\citation{#1}}%
}%
\@innernewif\if@notfirstcitation
%
%
\def\cite{\@getoptionalarg\@cite}%
%
%
\def\@cite#1{%
   \let\@citenotetext = \@optionalarg
   \printcitestart
   \nocite{#1}%
   \@notfirstcitationfalse
   \@for \@citation :=#1\do
   {%
      \expandafter\@onecitation\@citation\@@
   }%
   \ifx\empty\@citenotetext\else
      \printcitenote{\@citenotetext}%
   \fi
   \printcitefinish
}%
\newif\ifweareinprivate
\weareinprivatetrue
\ifx\shlhetal\undefinedcontrolseq\weareinprivatefalse\fi
\ifx\shlhetal\relax\weareinprivatefalse\fi
\def\@onecitation#1\@@{%
   \if@notfirstcitation
      \printbetweencitations
   \fi
   \expandafter \ifx \csname\@citelabel{#1}\endcsname \relax
      \if@citewarning
         \message{\@linenumber Undefined citation `#1'.}%
      \fi
     \ifweareinprivate
      \expandafter\gdef\csname\@citelabel{#1}\endcsname{%
\strut 
\vadjust{\vskip-\dp\strutbox
\vbox to 0pt{\vss\parindent0cm \leftskip=\hsize 
\advance\leftskip3mm
\advance\hsize 4cm\strut\openup-4pt 
\rightskip 0cm plus 1cm minus 0.5cm ?  #1 ?\strut}}
         {\tt
            \escapechar = -1
            \nobreak\hskip0pt\pfeilsw
            \expandafter\string\csname#1\endcsname
\pfeilso
            \nobreak\hskip0pt
         }%
      }%
     \else  
      \expandafter\gdef\csname\@citelabel{#1}\endcsname{%
            {\tt\expandafter\string\csname#1\endcsname}
      }%
     \fi  
   \fi
   \csname\@citelabel{#1}\endcsname
   \@notfirstcitationtrue
}%
%
%
\def\@citelabel#1{b@#1}%
%
%
\def\@citedef#1#2{\expandafter\gdef\csname\@citelabel{#1}\endcsname{#2}}%
%
%
%
\def\@readbblfile{%
   \ifx\@itemnum\@undefined
      \@innernewcount\@itemnum
   \fi
   \begingroup
      \def\begin##1##2{%
         \setbox0 = \hbox{\biblabelcontents{##2}}%
         \biblabelwidth = \wd0
      }%
      \def\end##1{}
      %
      %
      \@itemnum = 0
      \def\bibitem{\@getoptionalarg\@bibitem}%
      \def\@bibitem{%
         \ifx\@optionalarg\empty
            \expandafter\@numberedbibitem
         \else
            \expandafter\@alphabibitem
         \fi
      }%
      \def\@alphabibitem##1{%
         \expandafter \xdef\csname\@citelabel{##1}\endcsname {\@optionalarg}%
         \ifx\biblabelprecontents\@undefined
            \let\biblabelprecontents = \relax
         \fi
         \ifx\biblabelpostcontents\@undefined
            \let\biblabelpostcontents = \hss
         \fi
         \@finishbibitem{##1}%
      }%
      \def\@numberedbibitem##1{%
         \advance\@itemnum by 1
         \expandafter \xdef\csname\@citelabel{##1}\endcsname{\number\@itemnum}%
         \ifx\biblabelprecontents\@undefined
            \let\biblabelprecontents = \hss
         \fi
         \ifx\biblabelpostcontents\@undefined
            \let\biblabelpostcontents = \relax
         \fi
         \@finishbibitem{##1}%
      }%
      \def\@finishbibitem##1{%
         \biblabelprint{\csname\@citelabel{##1}\endcsname}%
         \@writeaux{\string\@citedef{##1}{\csname\@citelabel{##1}\endcsname}}%
         \ignorespaces
      }%
      %
      %
      \let\em = \bblem
      \let\newblock = \bblnewblock
      \let\sc = \bblsc
      \frenchspacing
      \clubpenalty = 4000 \widowpenalty = 4000
      \tolerance = 10000 \hfuzz = .5pt
      \everypar = {\hangindent = \biblabelwidth
                      \advance\hangindent by \biblabelextraspace}%
      \bblrm
      \parskip = 1.5ex plus .5ex minus .5ex
      \biblabelextraspace = .5em
      \bblhook
      \input \bblfilebasename.bbl
   \endgroup
}%
%
%
\@innernewdimen\biblabelwidth
\@innernewdimen\biblabelextraspace
%
%
%
\def\biblabelprint#1{%
   \noindent
   \hbox to \biblabelwidth{%
      \biblabelprecontents
      \biblabelcontents{#1}%
      \biblabelpostcontents
   }%
   \kern\biblabelextraspace
}%
%
%
%
\def\biblabelcontents#1{{\bblrm [#1]}}%
%
%
\def\bblrm{\rm}%
%
%
\def\bblem{\it}%
%
%
\def\bblsc{\ifx\@scfont\@undefined
              \font\@scfont = cmcsc10
           \fi
           \@scfont
}%
%
%
\def\bblnewblock{\hskip .11em plus .33em minus .07em }%
%
%
\let\bblhook = \empty
%
%
%
\def\printcitestart{[}
\def\printcitefinish{]}
\def\printbetweencitations{, }
\def\printcitenote#1{, #1}
%
%
%
\let\citation = \@gobble
%
%
%
\@innernewcount\@numparams
%
%
\def\newcommand#1{%
   \def\@commandname{#1}%
   \@getoptionalarg\@continuenewcommand
}%
%
%
\def\@continuenewcommand{%
   \@numparams = \ifx\@optionalarg\empty 0\else\@optionalarg \fi \relax
   \@newcommand
}%
%
%
\def\@newcommand#1{%
   \def\@startdef{\expandafter\edef\@commandname}%
   \ifnum\@numparams=0
      \let\@paramdef = \empty
   \else
      \ifnum\@numparams>9
         \errmessage{\the\@numparams\space is too many parameters}%
      \else
         \ifnum\@numparams<0
            \errmessage{\the\@numparams\space is too few parameters}%
         \else
            \edef\@paramdef{%
               \ifcase\@numparams
                  \empty  No arguments.
               \or ####1%
               \or ####1####2%
               \or ####1####2####3%
               \or ####1####2####3####4%
               \or ####1####2####3####4####5%
               \or ####1####2####3####4####5####6%
               \or ####1####2####3####4####5####6####7%
               \or ####1####2####3####4####5####6####7####8%
               \or ####1####2####3####4####5####6####7####8####9%
               \fi
            }%
         \fi
      \fi
   \fi
   \expandafter\@startdef\@paramdef{#1}%
}%
%
%
%
%
\def\@readauxfile{%
   \if@auxfiledone \else 
      \global\@auxfiledonetrue
      \@testfileexistence{aux}%
      \if@fileexists
         \begingroup
            \endlinechar = -1
            \catcode`@ = 11
            \input \jobname.aux
         \endgroup
      \else
         \message{\@undefinedmessage}%
         \global\@citewarningfalse
      \fi
      \immediate\openout\@auxfile = \jobname.aux
   \fi
}%
%
%
\newif\if@auxfiledone
\ifx\noauxfile\@undefined \else \@auxfiledonetrue\fi
%
%
%
%
\@innernewwrite\@auxfile
\def\@writeaux#1{\ifx\noauxfile\@undefined \write\@auxfile{#1}\fi}%
%
%
%
\ifx\@undefinedmessage\@undefined
   \def\@undefinedmessage{No .aux file; I won't give you warnings about
                          undefined citations.}%
\fi
%
%
\@innernewif\if@citewarning
\ifx\noauxfile\@undefined \@citewarningtrue\fi
%
%
%
\catcode`@ = \@oldatcatcode

\def\pfeilso{\leavevmode
            \vrule width 1pt height9pt depth 0pt\relax
           \vrule width 1pt height8.7pt depth 0pt\relax
           \vrule width 1pt height8.3pt depth 0pt\relax
           \vrule width 1pt height8.0pt depth 0pt\relax
           \vrule width 1pt height7.7pt depth 0pt\relax
            \vrule width 1pt height7.3pt depth 0pt\relax
            \vrule width 1pt height7.0pt depth 0pt\relax
            \vrule width 1pt height6.7pt depth 0pt\relax
            \vrule width 1pt height6.3pt depth 0pt\relax
            \vrule width 1pt height6.0pt depth 0pt\relax
            \vrule width 1pt height5.7pt depth 0pt\relax
            \vrule width 1pt height5.3pt depth 0pt\relax
            \vrule width 1pt height5.0pt depth 0pt\relax
            \vrule width 1pt height4.7pt depth 0pt\relax
            \vrule width 1pt height4.3pt depth 0pt\relax
            \vrule width 1pt height4.0pt depth 0pt\relax
            \vrule width 1pt height3.7pt depth 0pt\relax
            \vrule width 1pt height3.3pt depth 0pt\relax
            \vrule width 1pt height3.0pt depth 0pt\relax
            \vrule width 1pt height2.7pt depth 0pt\relax
            \vrule width 1pt height2.3pt depth 0pt\relax
            \vrule width 1pt height2.0pt depth 0pt\relax
            \vrule width 1pt height1.7pt depth 0pt\relax
            \vrule width 1pt height1.3pt depth 0pt\relax
            \vrule width 1pt height1.0pt depth 0pt\relax
            \vrule width 1pt height0.7pt depth 0pt\relax
            \vrule width 1pt height0.3pt depth 0pt\relax}

\def\pfeilsw{ \leavevmode 
            \vrule width 1pt height0.3pt depth 0pt\relax
            \vrule width 1pt height0.7pt depth 0pt\relax
            \vrule width 1pt height1.0pt depth 0pt\relax
            \vrule width 1pt height1.3pt depth 0pt\relax
            \vrule width 1pt height1.7pt depth 0pt\relax
            \vrule width 1pt height2.0pt depth 0pt\relax
            \vrule width 1pt height2.3pt depth 0pt\relax
            \vrule width 1pt height2.7pt depth 0pt\relax
            \vrule width 1pt height3.0pt depth 0pt\relax
            \vrule width 1pt height3.3pt depth 0pt\relax
            \vrule width 1pt height3.7pt depth 0pt\relax
            \vrule width 1pt height4.0pt depth 0pt\relax
            \vrule width 1pt height4.3pt depth 0pt\relax
            \vrule width 1pt height4.7pt depth 0pt\relax
            \vrule width 1pt height5.0pt depth 0pt\relax
            \vrule width 1pt height5.3pt depth 0pt\relax
            \vrule width 1pt height5.7pt depth 0pt\relax
            \vrule width 1pt height6.0pt depth 0pt\relax
            \vrule width 1pt height6.3pt depth 0pt\relax
            \vrule width 1pt height6.7pt depth 0pt\relax
            \vrule width 1pt height7.0pt depth 0pt\relax
            \vrule width 1pt height7.3pt depth 0pt\relax
            \vrule width 1pt height7.7pt depth 0pt\relax
            \vrule width 1pt height8.0pt depth 0pt\relax
            \vrule width 1pt height8.3pt depth 0pt\relax
            \vrule width 1pt height8.7pt depth 0pt\relax
            \vrule width 1pt height9pt depth 0pt\relax
      }


\def\widestnumber#1#2{}

\def\citewarning#1{\ifx\shlhetal\relax 
    \else
    \par{#1}\par
    \fi
}

\def\rm{\fam0 \tenrm}

\def\fakesubhead#1\endsubhead{\bigskip\noindent{\bf#1}\par}



%
%
%

%

\font\textrsfs=rsfs10
\font\scriptrsfs=rsfs7
\font\scriptscriptrsfs=rsfs5

\newfam\rsfsfam
\textfont\rsfsfam=\textrsfs
\scriptfont\rsfsfam=\scriptrsfs
\scriptscriptfont\rsfsfam=\scriptscriptrsfs

\edef\oldcatcodeofat{\the\catcode`\@}
\catcode`\@11

\def\Cal@@#1{\noaccents@ \fam \rsfsfam #1}

\catcode`\@\oldcatcodeofat


\expandafter\ifx \csname margininit\endcsname \relax\else\margininit\fi

\pageheight{8.5truein}
\topmatter
\title{What majority decisions are possible \\
 Sh816} \endtitle
\author {Saharon Shelah \thanks {\null\newline 
Partially supported by the United States-Israel Binational Science
Foundation \null\newline
I would like to thank 
Alice Leonhardt for the beautiful typing. \null\newline
 First Typed - 02/Jan/16 \null\newline 
 Latest Revision - 03/Mar/6} \endthanks} \endauthor 
\affil{Institute of Mathematics\\
 The Hebrew University\\
 Jerusalem, Israel
 \medskip
 Rutgers University\\
 Mathematics Department\\
 New Brunswick, NJ  USA} \endaffil
\endtopmatter
\document  
 
\newpage

\head {Annotated Content} \endhead  \resetall 
\bn
\S0 Introduction
\bn
\S1 Basic definitions and facts
\bn
\S2 When every majority choice is possible: a characterization
\bn
\S3 Balanced choice functions
\mr
\item "{${{}}$}" [We characterize what the majority choice can be on
pr-$c \ell({\Cal D})$ for ${\Cal D} \subseteq {\frak C}$ which is
balanced, i.e., does not fall under \S2.  We get the full answer.]
\endroster
\bn
\S4 Generalization
\mr
\item "{${{}}$}" [We consider allowing the original voters to abstain
(i.e., allow a draw) and more.]
\endroster
\newpage

\head {\S0 Introduction} \endhead  \resetall \sectno=0
\bigskip

Condorcet's ``paradox'' demonstrates that given three candidates A, B and C,
the majority rule may result in the 
society preferring A to B , B to C and C to A.
McGarvey \cite {McG53} proved a far-reaching extension of Condorcet's paradox:
For every asymmetric relation $R$ on a finite set $M$ of
candidates there is a strict-preferences
(linear orders, no ties) voter profile that has the relation $R$
as its strict simple majority relation. In other words,
for every assymetric relation (equivalently, a tournament)
$R$ on a set $M$ of $m$ elements
there are $n$ order relations on $M$, $R_1,R_2, \dots, R_n$
such that for every $a,b \in M$, $aRb$ if and only if
$$|\{ i: aR_ib\}|> n/2.$$ McGarvey's proof gave $n = m(m-1)$.
Stearns \cite {Ste59} found a construction with
$n = m$ and noticed that a simple 
counting argument implies that $n$ must be at least  $m/\log m$.
Erd\H{o}s and Moser \cite {ErMo64} were able to give a construction
with $n=O(m/\log m)$. Alon \cite {Alo02} showed that for some constant
$c_1>0$ we can find $R_1,\dots, R_n$
with $$|\{ i: aR_ib\}|> (1/2+ c_1/\sqrt {n})n,$$ and
that this is no longer the case if $c_1$ is
replaced with another constant $c_2>0$.
\sn
Gil Kalai asked to what extent the assertion of McGarvey's theorem
holds if we replace the set of order relations
by an arbitrary isomorphism class of choice functions on pairs of
elements (or equivalently tournaments).
Namely, under which conditions (A) below holds (i.e., question \scite{1.3}).
\nl
Instead of choice functions we can speak on tournaments. \nl

The main result is (follows by \scite{2.1})
\demo{\stag{0.1} Conclusion}  Let $X$ be a finite set and ${\frak D}$
be a non empty family of choice functions for $\binom X 2$ closed under
permutation of $X$.  Then the following conditions are equivalent:
\mr
\item "{$(A)$}"  for any choice function $c$ on $\binom X2$ we can
find a finite set $J$ and $c_j \in {\frak D}$ for $j \in J$ such that
for any $x \ne y \in X$:
$$
c\{x,y\} = y \text{ \ub{iff} } |J|/2 < \{j \in J:c_j\{x,y\} \setminus y\}
$$
(so equality never occurs)
\sn
\item "{$(B)$}"  for some $c \in {\frak D}$ and $x \in X$ we have
$|\{y:c\{x,y\} = y\}| \ne (|X|-1)/2$.
\ermn
Gil Kalai further asks \nl
\margintag{0.q}\ub{\stag{0.q} Question}:  1) In \scite{0.1} can we bound $|J|$? 
\nl
2) What is the result of demanding a ``nontrivial majority"? (say 51\%?)
\mn
Under \scite{0.1} it seems reasonable to ask to characterize the two place
relation $\bold R$ on the family of choice for $\binom X2$, see
\scite{1.4}.  We then give a complete solution: what is the closure of
a set of choice functions by majority; in fact, there are just two
possibilities (see \scite{b.5}) we then discuss a generalization (in
\S4). \nl
I thank Gil for the simulating discussion and writing the historical
background. 
\enddemo
\bn
\margintag{0.2}\ub{\stag{0.2} Notation}:  Let $n,m,k,\ell,i,j$ denote natural numbers. \nl
Let $r,s,t,a,b$ denote real numbers.  \nl
Let $x,y,z,u,v,w$ denote members of the finite set $X$. \nl
Let $\binom Xk$ be the family of subsets of $X$ with exactly $k$
members. \nl
Let $c,d$ denote choice functions on $\binom X2$. \nl
Let conv$(A)$ be the convex hull of $A$, here for $A \subseteq \Bbb R
\times \Bbb R$. \nl
Let Per$(X)$ be the set of permutations of $X$.
\newpage

\head {\S1 Basic definitions and facts} \endhead  \resetall \sectno=1
\bigskip

\demo{\stag{1.1} Hypothesis}  Assume
\mr
\item "{$(a)$}"  $X$ is a (fixed) finite set with $\bold n \ge 3$
members, i.e. $\bold n = |X|$
\sn
\item "{$(b)$}"  ${\frak C} = {\frak C}_X$ is the set of choice
functions on $\binom X2$
\sn
\item "{$(c)$}"  ${\frak D}$ is fixed as a nonempty subset of
${\frak C}$ which is symmetric where
\endroster
\enddemo
\bigskip

\definition{\stag{1.1A} Definition}  1) ${\Cal C} \subseteq {\frak
C}$ is symmetric if it is closed under permutations of $X$ 
(i.e. for every $\pi \in \text{ Per}(X)$ the permutation $\hat \pi$ maps
${\Cal C}$ onto itself where $\pi$ induces $\hat \pi$, 
a permutation of ${\frak C}$, that is we can write $c_1 = c^\pi_2$ or $c_1
= \hat \pi c_2$ where: $x_1 = \pi(x_2),y_1 = \pi(y_2)$ implies
$c_1\{x_1,y_1\} = y_1 \Leftrightarrow c_2\{x_2,y_2\} = y_2$).
\enddefinition
\bigskip

\definition{\stag{1.2} Definition}  For ${\Cal D} \subseteq {\frak
C}$ let maj-$cl({\Cal D})$ be the set of $d \in {\frak C}$ such that
for some real numbers $r_c \in [0,1]_{\Bbb R}$ for $c \in {\Cal D}$
satisfying $\underset{c \in {\Cal D}} {}\to \Sigma r_c = 1$ we have
\footnote{note that there is no reason to assume that ${\Cal D}_2 =
\text{ maj}-c \ell({\Cal D}_1)$ implies ${\Cal D}_2 = \text{ maj}-c
\ell({\Cal D}_2)$}

$$
d(\{x,y\}) = x \Leftrightarrow \frac 12 < \Sigma\{r_c:c\{x,y\} = x
\text{ and } c \in {\Cal D}\}
$$
(so the sum is never $\frac 12$).
\enddefinition
\bn
So Kalai's question was \nl
\margintag{1.3}\ub{\stag{1.3} Question}:  Do we have (for $|X|$ large enough),
maj-cl$({\Cal D}^*) = {\frak C}$?
\bigskip

\definition{\stag{1.4} Definition}  1) Let Dis $= \text{ Dis}(X) = \{\mu:\mu$ a
distribution on ${\frak C}_X\}$; of course, ``$\mu$ a distribution on
${\frak C}$" means $\mu$ is a function from ${\frak C}$ into $[0,1]_{\Bbb
R}$ such that $\Sigma\{\mu(c):c \in {\frak C}\} = 1$, and we let
$\mu({\Cal C}) = \Sigma\{\mu(c):c \in {\Cal C}\}$. \nl
2) For ${\Cal C} \subseteq {\frak C}$ and $\mu \in \text{Dis}({\frak C})$
let $\mu({\Cal C}) = \Sigma\{\mu(c):c \in {\Cal C}\}$ so $\mu(\{{\Cal C}\})
\ge 0,\mu({\frak C}) = 1$. \nl
3) For ${\Cal D} \subseteq {\frak C}$ let Dis$_{\Cal D} = \{\mu \in
\text{ Dis}:\mu({\Cal D}) =1\}$. \nl
4) Let pr$({\frak C}) = \{\bar t:\bar t = \langle t_{x,y}:x \ne y
\in X \rangle\}$ such that $t_{x,y} \in [0,1]_{\Bbb R},t_{y,x} =
1-t_{x,y}$, we may write $\bar t(x,y)$; pr stands for probability. \nl
5) For $T \subseteq \text{ pr}({\frak C})$ let pr-cl$(T)$ be the convex hull
of $T$. \nl
6) For $d \in {\frak C}$ let $\bar t[d] = \langle t_{x,y}[d]:x \ne y
\in X \rangle$ be defined by $[t_{x,y}[d] = 1 \Leftrightarrow d\{x,y\} = y
\Leftrightarrow t_{x,y}[d] \ne 0]$. \nl
7) Let pr-cl$({\Cal D})$ for ${\Cal D} \subseteq {\frak C}$ be
pr-cl$(\{\bar t[d]:d \in {\Cal D}\})$, prd$({\Cal D}) = \{\bar t[c]:c
\in {\Cal D}\}$. \nl
8) For ${\Cal C} \subseteq {\frak C}$ we let sym-$c \ell({\Cal C})$ be
the minimal ${\Cal D} \subseteq {\frak C}$ which is symmetric and
includes ${\Cal C}$.
\nl
9) For $T \subseteq \text{ pr}({\frak C})$ let maj$(T) = \{c \in {\frak
C}$: for some $\bar t \in T$ for any $x \ne y$ from $X$ we have
$c\{x,y\} = y \Leftrightarrow t_{x,y} > \frac 12 \Leftrightarrow
t_{x,y} \ge \{ \frac 12\}\}$. 
\enddefinition
\bigskip

\proclaim{\stag{1.8} Claim}  
1) For $d \in {\frak C}$ we have $\bar t[d] \in \,{\text{\rm pr\/}}({\frak
C})$. \nl
2) For ${\Cal D} \subseteq {\frak C}$ we have 
${\text{\rm Dis\/}}_{\Cal D} \subseteq { \text{\rm Dis\/}}$. \nl
3) ${\text{\rm prd\/}}({\frak C}) \subseteq { \text{\rm pr\/}}({\frak C})$ 
and if ${\Cal D} \subseteq {\frak C}$ then 
${\text{\rm prd\/}}({\Cal D}) \subseteq {\text{\rm pr-cl\/}}({\Cal D})
\subseteq { \text{\rm Dis\/}}$. \nl
4) If ${\Cal C} \subseteq {\frak C}$ then ${\Cal C} 
\subseteq { \text{\rm sym\/}}-c \ell({\Cal C}) \subseteq {\frak C}$. \nl
5) If $T \subseteq { \text{\rm pr\/}}({\frak C})$ \ub{then} 
${\text{\rm maj\/}}(T) \subseteq {\frak C}$. \nl
6) For ${\Cal D} \subseteq {\frak C}$ we have ${\text{\rm maj\/}}-c
\ell({\Cal D}) = {\text{\rm maj(pr-cl\/}}({\Cal D}))$.
\endproclaim
\bigskip

\demo{Proof}  Obvious.
\enddemo
\bigskip

\proclaim{\stag{1.9} Claim}  (G. Kalai) If ${\Cal C} \subseteq {\frak
C}$ and for every $c \in {\Cal C}$ and $x \in X$, the in-valency and
out-valency are equal, (i.e., ${\text{\rm val\/}}_c(x) = (|X|-1)/2$,
see below) \ub{then} every $d \in {\text{\rm maj-cl\/}}({\Cal C})$
satisfies:
\mr
\item "{$(*)$}"  if $\emptyset \ne Y \subsetneqq X$ then there are
edges from $Y$ to $X \backslash Y$ and from $X \backslash Y$ to $Y$, see
\scite{b.2}(1),(2). 
\endroster
\endproclaim
\bigskip

\definition{\stag{1.10} Definition}  1) For $d \in {\frak C}$ and $x
\in X$ let val$_d(x)$, the valency of $x$ for $d$, be $|\{y:y \in X,y
\ne x$ and $d\{x,y\} = y\}|$ so val$_d(x) \in \{0,\dotsc,|\bold n|-1\}$. \nl
We also call val$_d(x)$ the out-valency \footnote{natural under the
tournament interpretation} of $x$ in $d$ and $\bold n -
\text{ val}_d(x)-1$ the in-valency of $x$ in $d$. \nl
2) For $d \in {\frak C}$ let Val$(d) = \{\text{val}_d(x):x \in X\}$.
\nl
3) For $d \in {\frak C}$ and $\ell \in \{0,1\}$ let $V_\ell(d) =
\{(\text{val}_d(x)$,val$_d(y)):x \ne y \in X$ and $d\{x,y\} = y
\Leftrightarrow \ell = 1\}$. \nl
4) For $d \in {\frak C}$ and $\ell \in \{0,1\}$ let $V^*_\ell(d) =
\{\bar k-(\ell,1-\ell):\bar k \in V_\ell(d)\}$ and let $V^*(d) = V^*_0(d) \cup
V^*_1(d)$.
\nl
5) For $c \in {\frak C}$ let ${\text{\rm
dual\/}} (c) \in {\frak C}$ be ${\text{\rm dual\/}}(c)\{x,y\} 
\in \{x,y\} \backslash \{c\{x,y\}\}$; similarly $\bar t' = \text{ dual}
(\bar t)$ for $\bar t \in \text{ pr}({\frak C})$ means that $t'_{x,y}
= 1 - t_{x,y}$.
\enddefinition
\bigskip

\proclaim{\stag{1.11} Claim}  1) 
\mr
\item "{$(\alpha)$}"  $c_1 \in { \text {\rm sym-cl\/}}\{c_2\} 
\,{\text {\rm \ub{iff} dual\/}}(c_1) \in 
{ \text{\rm sym-cl\/}}\{{\text{\rm dual\/}}\{c_2\}\}$ 
\sn
\item "{$(\beta)$}"    $c_1 \in { \text{\rm maj-cl\/}}
({\text{\rm sym-cl\/}}\{c_2\}) 
\,{\text {\rm \ub{iff}\, dual\/}}(c_1) \in 
{\text{\rm maj-cl\/}}({\text{\rm sym-cl\/}}\{{\text{\rm dual\/}}(c_2)\})$.
\ermn
2) $(k_0,k_1) \in V_0(d) \Rightarrow (k_1,k_0) \in V_1(d)$. \nl
3) ``$c_2 \in { \text{\rm sym-cl\/}}\{c_2\}$" is an equivalence
relation on ${\frak C}$ and it implies $V_\ell(c_1) = V_\ell(c_2)$ for
$\ell=0,1$.  
\endproclaim
\bigskip

\demo{Proof}  Easy.
\enddemo
\newpage

\head {\S2 When every majority choice is possible: a characterization} \endhead  \resetall \sectno=2
\bigskip

The following is the main part of the solution
(probably $(c) \Leftrightarrow (g)$
should be stated as the main conclusion here).
\proclaim{\stag{2.1} Main Claim}  The following conditions on ${\frak D}
\subseteq {\frak C}$ which is symmetric and not empty, (i.e.,
${\frak D}$ is a set of choice
functions on $\binom X2$ closed under permutation on $X$, 
${\frak D} \ne \emptyset$) and for simplicity assuming for the time
being that ${\frak D} = { \text{\rm sym-cl\/}}(d^*)$ for 
any $d^* \in {\frak D}$, are equivalent, where $x,y$ vary
on distinct members of $X$:
\mr
\widestnumber\item{$(b)_{x,y}$}
\item "{$(a)$}"  ${\text{\rm maj-cl\/}}({\frak D}) = {\frak C}$
\sn
\item "{$(b)_{x,y}$}"  there is $\bar t \in { \text{\rm pr-cl\/}}({\frak D})
\subseteq {\text{\rm pr\/}}({\frak C})$ such that 
{\roster
\itemitem{ $(i)$ }  $t_{x,y} > \frac 12$
\sn
\itemitem{ $(ii)$ }  $\{x,y\} \ne \{u,v\} \in \binom Xr 
\Rightarrow t_{u,v} = \frac 12$
\endroster}
\item "{$(c)$}"  for any $c \in {\frak C}$ we can find a finite set
$J$ and sequence $\langle d_j:j \in J \rangle$ such that $d_j \in
{\frak D}$ and: if $u \ne v \in X$ then \nl
$c\{u,v\} = v \Leftrightarrow (\{j \in J:d_j\{u,v\} = j\}| > |J|/2$
\sn
\item "{$(d)$}"  $(\frac 12,\frac 12)$ belongs to 
${\text{\rm Pr\/}}_{> \frac 12}({\frak D})$, see definition below
\sn
\item "{$(e)$}"  $(\frac 12,\frac 12) \in {\text{\rm Pr\/}}_{\ne 1/2}({\frak
D})$
\sn
\item "{$(f)$}"  $(\frac{\bold n}2 - 1,\frac{\bold n}2 -1)$ can be
represented as $r^*_0 \times \bar s_0 + r^*_1 \times \bar s_1$ where 
{\roster
\itemitem{ $(*)(i)$ }  $r^*_0,r^*_1 
\in [0,1]_{\Bbb R} \backslash \{\frac 12\}$
\sn
\itemitem{ $(ii)$ }  $1 = r^*_0 + r^*_1$
\sn
\itemitem{ $(iii)$ }  for $\ell =0,1$ the pair 
$\bar s_\ell \text{ belongs to the convex hull of } V^*_\ell(d^*)
\text{ for some } d^* \in {\frak D}$, see Definition
\scite{1.10}(3), but recall that by a hypothesis of the claim, the
choice of $d^*$ is immaterial
\endroster}
\sn
\item "{$(g)$}"  for some $(d^* \in {\frak D}$ and) $x \in X$ we have
${\text{\rm val\/}}_{d^*}(x) \ne \frac{\bold n-1}2$. 
\endroster
\endproclaim
\bigskip

\demo{Proof}  $(b)_{x,y} \Leftrightarrow (b)_{x',y'}$:

(So $x,y,x',y' \in X$ and $x \ne y,x' \ne y'$).  
Trivial as ${\frak D}$ is closed under permutations of $X$.
\mn
\ub{$(b)_{x,y} \Rightarrow (a)$}:

Let $c \in {\frak C}^*$. \nl
Let $\{(u_i,v_i):i < i(*)\}$ list the pairs $(u,v)$ of distinct
members of $X$ such that $c\{u,v\} = v$ so $i(*) = \binom{|X|}{2}$.  For
each $i < i(*)$ as $(b)_{x,y} \Rightarrow (b)_{u_i,v_i}$ clearly
there is $\bar t^i \in \text{ pr-cl}({\frak D}^*)$
such that

$$
t^i_{u_i,v_i} > \frac 12 \text{ so } t^i_{v_i,u_i} = 1 - t_{u_i,v_i} <
\frac 12
$$

$$
\{u_i,v_i\} \ne \{u,v\} \in \binom X2 \Rightarrow t^i_{u,v} = \frac
12.
$$
\mn
Let $\bar t^* = \langle t^*_{u,v}:u \ne v \in X \rangle$ be defined
to

$$
t^*_{u,v} = \Sigma\{t^i_{u,v}:i < i(*)\}/i(*).
$$
\mn
As pr-cl$({\frak D})$ is convex and $i < i(*) \Rightarrow \bar t^i
\in \text{ pr-cl}({\frak D})$ clearly $\bar t \in \text{ pr-cl}
({\frak D})$.  Now for each $j < i(*),t^i_{u_j,v_j}$ is 
$\frac 12$ if $i \ne j$ and is $> \frac 12$ if $i=j$.   Hence
$t^*_{u_j,v_j}$ being the average of $\langle t^i_{u_j,v_j}:i < i(*)
\rangle$ is $> \frac 12$.  Hence $t^*_{v_j,u_j} = 1-t^*_{u_j,v_j} <
\frac 12$.  So $t^*_{u,v} > \frac 12$ \ub{iff} $t^*_{u,v} \ge \frac
12$ \ub{iff} $c\{u,v\} = v$ by the
choice of $\langle(u_i,v_i):i < i(*) \rangle$.  So $\bar t^*$ witness
$c \in$ maj-cl$({\frak D})$ as required in clause (a).
\mn
\ub{$(a) \Rightarrow (b)_{x,y}$}: 

By clause (a), for every $d \in {\frak C}$ there is $\langle r_c:c \in
{\frak D} \rangle$ as in Definition \scite{1.2}, so some
$\varepsilon_d > 0,d(\{x,y\}) = y \Rightarrow \frac 12 + \varepsilon_d
< \Sigma\{r_c:c \in {\frak D}$ and $c\{x,y\} = y\}$.  Hence
$\varepsilon = \text{ Min}\{\varepsilon_d:d \in {\frak D}\}$ 
is a real $>0$. \nl
Let $T = \{\bar t:\bar t \in \text{ pr-cl}({\frak D})$ and $t_{x,y}
\ge \frac 12 + \varepsilon$ and $u \ne v \in X \Rightarrow t_{u,v} \le
t_{x,y}\}$, so
\mr
\item "{$(*)_1$}"  $T \ne \emptyset$
\nl
[Why?  As ${\frak D}$ is symmetric clearly ${\frak D}_{x,y} = \{c \in {\frak
D}:c\{x,y\} = y\}$ is non empty, and $\bar t[d] \in T$ for every $d
\in {\frak D}_{x,y}$]
\sn
\item "{$(*)_2$}"    $T$ is convex and closed.
\ermn
[Why?  Trivial.] \nl
For $\bar t \in T$ define 
\mr
\item "{$\boxtimes$}"  err$(T) = \max\{|t_{u,v} - \frac 12|:u \ne
v \in X \text{ and } \{u,v\} \ne \{x,y\}\}$
\sn
\item "{$(*)_3$}"   if $\bar t \in T$, err$(\bar t) > 0$ \ub{then} we
can find $\bar t' \in T$ 
such that err$(\bar t') \le \frac 12 \text{ err}(t)$ and $t'_{x,y} \ge
(t_{x,y} + \frac 12 + \varepsilon)/2$.
\ermn
Why?  Choose $c \in {\frak C}$ such that $c\{x,y\} = y$ and

$$
u \ne v \in X \and \{u,v\} \ne \{x,y\} \and t_{u,v} > \frac 12
\Rightarrow c\{u,v\} = u
$$
\mn
(so if $t_{u,v} = t_{v,u} = \frac 12$ it does not matter).

So $c$ is ``a try to correct $\bar t$".

As we are assuming clause (a) and the choice of $\varepsilon_d$,
we can find $\bar r^* = \langle r^*_d:d
\in {\frak D} \rangle$ with $r^*_d \in [0,1]_{\Bbb R}$ and $1 =
\Sigma\{r^*_d:d \in {\frak D}\}$ such that  

$$
\frac 12 + \varepsilon_d < \Sigma\{r^*_d:d\{x,y\} = y\}
$$

$$
c\{u,v\} = v \Rightarrow \frac 12 < \Sigma\{r^*_d:d\{u,v\} = v\}
$$
hence

$$
c\{u,v\} = u \Rightarrow \frac 12 < \Sigma\{r^*_d:d\{u,v\} = u\}.
$$
\mn
By the choice of $\varepsilon,\frac 12 + \varepsilon < \Sigma\{r^*_d:d\{x,y\}
= y\}$. \nl
Let $\bar s = \langle s_{u,v}:u \ne v \in X \rangle$ be defined by
$s_{u,v} = \Sigma\{r^*_d:d\{u,v\} = v\}$, so
\mr
\item "{$\circledast(i)$}"   $\bar s \in \text{ pr-cl}({\frak D}^*)$
\sn
\item "{$(ii)$}"  $s_{x,y} > \frac 12 + \varepsilon$ (so $s_{y,x} <
\frac 12$)
\sn
\item "{$(iii)$}"  if $t_{u,v} > \frac 12$ and $u \ne v \in X,\{u,v\}
\ne \{x,y\}$ then $c\{u,v\} = u$ hence $s_{u,v} < \frac 12$
\sn
\item "{$(iv)$}"  if $t_{u,v} < \frac 12$ and $u \ne v \in X,\{u,v\}
\ne \{x,y\}$ then $c\{u,v\} = v$ hence $s_{u,v} > \frac 12$.
\ermn
Let $\bar t' = \frac 12 \bar t + \frac 12 \bar s$, i.e. $t'_{u,v} =
\frac 12(t_{u,v} + s_{u,v})$ so clearly
\mr
\item "{$\circledast(i)$}"  $\bar t' \in \text{ pr-cl}({\frak D})$
\sn
\item "{$(ii)$}"  $t'_{x,y} = (t_{x,y} + \frac 12 + \varepsilon)/2$ (hence 
$\bar t' \in T$)
\sn
\item "{$(iii)$}"  if $ u \ne v \in X,\{u,v\} \ne \{x,y\}$ then
$|t'_{u,v} - \frac 12| \le \frac 12 |t_{u,v} - \frac 12|$.
\ermn
So we are done proving $(*)_3$.

As $T$ is closed (and included in a $\{\bar t:\bar t = \langle
t_{u,v}:u \ne v \in X \rangle,0 \le t_{u,v} \le 1\})$, clearly there
is $t \in T$ such that $u \ne v \in X \and \{u,v\} \ne \{x,y\}
\Rightarrow t_{u,v} = \frac 12$ as required.
\mn
\ub{$(c) \Rightarrow (a)$}:

Let $d \in {\frak C}$ and let $\langle c_j:j \in J \rangle$ witness
clause (c). \nl
Let $r_c = |\{j \in J:c_j = c\}|/|J|$ now $\langle r_c:c \in {\frak
D} \rangle$ witness clause (a), i.e., witness that $d \in$ 
maj-cl$({\frak D}^*)$.
\mn
\ub{$(a) \Rightarrow (c)$}:

Let $d \in {\frak C}$ and let $\langle r_c:c \in {\frak D}
\rangle$ be as guaranteed for $d$ by clause (a).  Let $n(*) > 0$ be large
enough and for $c \in {\frak D}$ let $k_c \in \{0,\dotsc,n(*)-1\}$ be
such that $c \in {\frak D} \Rightarrow k_c \le n(*) \times r_c < k_c +
1$; note that $k_c$ exists as $r_c \in [0,1]_{\Bbb R}$.
As $\dsize \sum_c \frac{k_c}{n(*)} \le 1 \le \dsize \sum_c
\frac{k_c+1}{n(*)}$, we can choose $m_c \in \{k_c,k_c+1\}$ such that $r'_c =
\frac{m_c}{n(*)}$ satisfies $\Sigma\{r'_c:c \in {\frak D}^*\} = 1$.
Let $J = \{(c,m):c \in {\frak D}$ and $m \in \{1,\dotsc,m_c\}\}$ and
we let $\bold c_{(d,m)} = d$ for $(d,m) \in J$.  
Now the ``majority" of $\langle \bold c_t:t
\in J \rangle$, see Definition \scite{1.2}, choose $d^*$ so clause (c) holds.
\enddemo
\bn
Before we deal with clauses (d),(e),(f) and (g) of \scite{2.1}, we define
\definition{\stag{2.2} Definition}  1) For ${\Cal D} \subseteq 
{\frak C}$ and $A \subseteq [0,1]_{\Bbb R}$ let 
Pr$_A({\Cal D})$ be the set of pairs $(s_0,s_1)$ of real
numbers $\in [0,1]_{\Bbb R}$ such that 
for some $\bar t \in \text{ pr-cl}({\Cal D})$ and $x \ne y \in X$ and
$a \in A$ we have $\bar t = \bar t \langle x,y,a,s_0,s_1 \rangle$
where \nl
2) $\bar t = \bar t\langle x,y,a,s_0,s_1 \rangle$ where
$x \ne y \in X,a \in [0,1]_{\Bbb R}$ and $s_0,s_1 \in [0,1]_{\Bbb R}$ and
$\bar t = \langle t_{u,v}:u \ne y \in X \rangle \in \text{ pr}({\frak
C})$ is defined by 
\mr
\item "{$(\alpha)$}"  $t_{x,y} = a$
\sn
\item "{$(\beta)$}"  if $z \in X \backslash \{x,y\}$ \ub{then}
$t_{y,z} = s_1$ (hence $t_{z,y} = 1 - s_1$)
\sn
\item "{$(\gamma)$}"  if $z \in X \backslash \{x,y\}$ \ub{then}
$t_{x,z} = s_0$ (hence $t_{z,x} = 1-s_0$)
\sn
\item "{$(\delta)$}"  if $z_1 \ne z_2 \in X \backslash \{x,y\}$
\ub{then} $t_{z_1,z_2} = \frac 12$.
\ermn
3) In Pr$_A({\Cal D})$ we may replace $A$ by $1,0 \ne,\frac 12,> \frac
12,< \frac 12$ if $A$ is
$\{1\},\{0\},[0,1]_{\Bbb R} \backslash \{\frac 12\},(\frac 12,1]_{\Bbb
R},[0,\frac 12)_{\Bbb R}$ respectively. \nl
4) For $\ell \in \{0,1\}$ and 
${\Cal D} \subseteq {\frak C}^*$ let Prd$_\ell({\Cal D})$ be the set
of pairs $\bar s = (s_0,s_1)$ of real (actually rational) 
numbers $\in [0,1]_{\Bbb R}$
such that for some $c \in {\Cal D}$ and $x \ne y \in X$ we have $\bar
s = \bar s^{c,x,y} = (s^{c,x,y}_0,s^{c,x,y}_1)$ where
\mr
\widestnumber\item{$(iii)$}
\item "{$(i)$}"  $s^{c,x,y}_1 = 
|\{z:z \in X \backslash \{x,y\} \text{ and } c\{y,z\} =
z\}|/(\bold n - 2)$
\sn
\item "{$(ii)$}"  $s^{c,x,y}_0 = 
|\{z:z \in X \backslash \{x,y\} \text{ and } c\{x,z\} = z\}|/(\bold n
- 2)$
\sn
\item "{$(iii)$}"  $\ell = 1 
\Leftrightarrow \ell \ne 0 \Leftrightarrow c\{x,y\} = y$.
\ermn
5) Prd$({\Cal D})$ is Prd$_0({\Cal D}) \cup \text{ Prd}_1({\Cal D})$
and ${\Cal D}_{x,y} = \{c \in {\Cal D}:c\{x,y\} = y\}$.
\enddefinition
\bigskip

\proclaim{\stag{2.3} Claim}   Let ${\Cal D} \subseteq {\frak C}$
\mr
\item  ${\text{\rm Pr\/}}_{A_1}
({\Cal D}_1) \subseteq { \text{\rm Pr\/}}_{A_2}({\Cal D}_2)$
if $A_1 \subseteq A_2 \subseteq [0,1]_{\Bbb R}$ and ${\Cal D}_1
\subseteq {\Cal D}_2 \subseteq {\frak C}$
\sn
\item   ${\text{\rm Pr\/}}_A({\Cal D})$ is a convex subset of
$[0,1]_{\Bbb R} \times [0,1]_{\Bbb R}$
\sn
\item  ${\text{\rm Prd\/}}_\ell
({\Cal D})$ is finite and its convex hull is $\subseteq
{\text{\rm Pr\/}}_{\{\ell\}}({\Cal D})$, increasing with ${\Cal D}(\subseteq
{\frak C})$ for $\ell =0,1$
\sn
\item  For $x \ne y \in X$ and $c \in {\frak C}$ satisfying $c\{x,y\}
= y \Leftrightarrow \ell =1 \Leftrightarrow \ell \ne 0$ we have (see
Claim \scite{2.2}(4))
$$
s^{c,x,y}_0 = ({\text{\rm val\/}}_c(x) - \ell)/(\bold n - 2)
$$

$$
s^{c,x,y}_1 = ({\text{\rm val\/}}_c(y) - (1 - \ell))/(\bold n - 2).
$$
\sn
\item  ${\text{\rm Pr\/}}_{A_1 \cup A_2}
({\Cal D}) = { \text{\rm Pr\/}}_{A_1}({\Cal D}) \cup
{ \text{\rm Pr\/}}_{A_2}({\Cal D})$, in fact ${\text{\rm Pr\/}}_A({\Cal D}) =
\cup\{{\text{\rm Pr\/}}_{\{a\}}({\Cal D}):a \in A\}$.
\endroster
\endproclaim
\bigskip

\demo{Proof}  Immediate.
\enddemo
\bigskip

\proclaim{\stag{2.3a} Claim}  1) If $x \ne y \in X$ and $c \in {\frak
C}$ and $\ell \in \{0,1\}$ satisfies $(\ell = 1) \equiv (c(\{x,y\})
=y)$ \ub{then} $\bar t \langle x,y,\ell,s^{c,x,y}_0,s^{c,x,y}_1
\rangle = \frac{1}{|\Pi_{x,y}|} \Sigma\{\bar t[\hat \pi(c)]:\pi \in
\Pi_{x,y}\}$ where $\Pi_{x,y} =: \{\pi \in { \text{\rm Per\/}}(X):\pi(x) = x$
and $\pi(y) = y\}$ hence $|\Pi_{x,y}| = (\bold n -2)!$. \nl
2) If ${\Cal D} \subseteq {\frak C}$ is symmetric and $\bar t \in 
{\text{\rm pr-cl\/}}({\Cal D})$ and $\bar t^* = \Sigma\{\bar t^\pi:\pi
\in \Pi_{x,y}\}/|\Pi_{x,y}|$ where $\bar t^\pi = \langle t^\pi_{u,v}:u
\ne v \in X \rangle,t^\pi_{u,v} = t_{\pi^{-1}(u),\pi^{-1}(v)}$
\ub{then} $t^* \in { \text{\rm pr-cl\/}}({\Cal D})$ and $\bar t^* \in
{ \text{\rm Pr\/}}_{\{a\}}(s_0,s_1)$ where $a = t_{x,y},s_0 =
\Sigma\{t_{x,z}:z \in X \backslash \{x,y\}\}/(\bold n-2)!,s_1 =
\Sigma\{t_{y;z}:z \in X \backslash \{y,z\}/(\bold n-2)!$.
\endproclaim
\bigskip

\demo{Proof}  Easy.
\enddemo
\bigskip

\proclaim{\stag{2.4} Claim}  For any symmetric ${\Cal D} \subseteq
{\frak C}$ (i.e., closed under permutations of $X$): \nl
1) For $\ell \in \{0,1\}$, the set ${\text{\rm Pr\/}}_\ell
({\Cal D})$ is the convex
hull of ${\text{\rm Prd\/}}_\ell({\Cal D})$ in $\Bbb R \times \Bbb R$. \nl
2) ${\text{\rm Pr\/}}_{[0,1]}({\Cal D})$ is the convex hull of 
${\text{\rm Prd\/}}({\Cal D})$. \nl
3) Let $s^*_0,s^*_1 \in [0,1]_{\Bbb R}$.  \ub{Then} $\bar t^* = \bar
t\langle x,y,a,s^*_0,s^*_1 \rangle \in {\text{\rm pr-cl\/}}({\Cal
D})$ \ub{iff} we can find $\langle r_{\bar s,\ell}:\ell \in \{0,1\}$
and $\bar s \in {\text{\rm Prd\/}}_\ell({\Cal D})\rangle$ such that
$r_{\bar s,\ell} \in [0,1]_{\Bbb R}$ and $1 = \Sigma\{r_{\bar
s,\ell}:\ell \in \{0,1\},\bar s \in 
{\text{\rm Prd\/}}_\ell({\Cal D})\}$ and $\bar s = \Sigma\{r_{\bar
s,\ell} \times \bar s:\ell \in \{0,1\},\bar s \in { \text{\rm
Prd\/}}_\ell({\Cal D})\}$
and $(\frac 12,\frac 12) = \Sigma\{r_{\bar s,\ell} \times \bar
s:\ell \in \{0,1\},\bar s \in {\text{\rm Prd\/}}_\ell({\Cal D})\}$ and
$a = \Sigma\{r_{\bar s,1}:\bar s \in {\text{\rm Prd\/}}_1({\Cal D})\}$.
\endproclaim
\bigskip

\demo{Proof}  1) By \scite{2.3}(3) we have one inclusion.

For the other direction assume $(s^*_0,s^*_1) \in \text{
Pr}_\ell({\Cal D})$ and we should prove that it belongs to the convex
hull of Prd$_\ell({\Cal D})$.
Fix $x \ne y \in X$ and let $\bar t^* = \bar t \langle
x,y,\ell,s^*_0,s^*_1 \rangle$ so
$\bar t^* = \langle t^*_{u,v}:u \ne v \in X \rangle$ is defined as
follows $t^*_{x,y} = \ell,t^*_{u,v} = \frac 12$ if 
$u \ne v \in X \backslash \{x,y\},t^*_{y,z} =
s^*_1$ if $z \in X \backslash \{x,y\},t_{x,z} = s^*_0$ if $z \in X
\backslash \{x,y\}$.
As $(s^*_0,s^*_1) \in \text{ Pr}_\ell({\Cal D})$ by 
Definition\scite{2.2} we know that $\bar t^* \in \text{ pr-cl}({\Cal D})$ and
let $\bar r = \langle r_c:c \in {\Cal D} \rangle$ be such that 
\mr
\item "{$\circledast_0$}"  $x,y,\bar r$ witness that 
$\bar t^* \in \text{ pr-cl}({\Cal D})$, so
$r_c \ge 0$ and $1 = \Sigma\{r_c:c \in {\Cal D}\}$ and

$$
\bar t^* = \Sigma\{r_c \times \bar t[c]:c \in {\Cal D}\}.
$$
\ermn
As $t^*_{x,y} = \ell$, necessarily
\mr
\item "{$\boxtimes_1$}"   $r_c \ne 0 \Rightarrow c \in {\Cal
D}_{x,y} =: \{c \in {\Cal D}:(c\{x,y\} = y) \equiv (\ell=1)\}$.
\ermn
To make the rest of the proof also a proof of part (3) let $a = \ell$
(as the real number $a$ maybe $\ne 0,1$ we need $m \in \{0,1\}$ below).
\nl
Let $\Pi_{x,y} = \{\pi \in \text{ Per}(X):\pi(x) = x,\pi(y) = y\}$ and
recall that for $\pi \in \text{ Per}(X),\hat \pi$ is the 
permutation of ${\frak C}$ which $\pi$ induces so $\hat
\pi$ maps ${\Cal D}_{x,y}$ onto ${\Cal D}_{x,y}$.  Clearly
$|\Pi_{x,y}| = (\bold n - 2)!$
\nl
For $(s_0,s_1) \in \text{ Prd}_\ell({\Cal D}) \subseteq 
\{(\frac{m_1}{(\bold n-2)!},\frac{m_2}{(\bold n-2)!}):m_1,m_2 \in 
\{0,1\dotsc,(\bold n-2)!\}\}$ let $r^*_{(s_0,s_1)} =
\Sigma\{r_c:\bar s^{c,x,y} = \bar s\}$ and for $m \in \{0,1\}$ we 
let $\bar r^*_{(s_0,s_1),m}
= \Sigma\{r_c:\bar s^{c,x,y} = \bar s$ and $c\{x,y\} = y \equiv (m=1)\}$.

Clearly $\pi \in \Pi_{x,y} \Rightarrow \langle t^*_{\pi(u),\pi(v)}:u \ne
v \in X \rangle = \bar t^*$, just check the definition, hence
recalling \scite{2.3}(4) we have
\mr
\item "{$\boxtimes_2$}"  $\bar t^* = \dsize \frac{1}{| \Pi_{x,y}|} \dsize
\sum_{\pi \in \Pi_{x,y}} \langle t^*_{\pi(u),\pi(v)}:u \ne v \in X
\rangle = \frac{1}{|\Pi_{x,y}|} \dsize \sum_{\pi \in \Pi_{x,y}} \,
\dsize \sum_{c \in {\Cal D}} r_c \times \bar t^*[\hat \pi(c)] =$ \nl 
$\dsize \sum_{c \in {\Cal D}} \, r_c(\frac{1}{|\Pi_{x,y}|} \, \dsize
\sum_{\pi \in \Pi_{x,y}} \bar t^*[\hat \pi(c)]) =
\dsize \sum_{c \in {\Cal D}} r_c \times \bar t \langle
x,y,a,s^{c,x,y}_0,s^{c,x,y}_1 \rangle =$ \nl
$\dsize \sum_{m \in \{0,1\}} \,\, 
\dsize \sum_{(s_0,s_1) \in \text{ Pr}_\ell({\Cal D})} 
\, (\Sigma\{r_c:c \in {\Cal D},c\{x,y\} = y \equiv (k=1)$ \nl

\hskip110pt and $\bar s^{c,x,y} = (s_0,s_1)\}) \times
\bar t \langle x,y,m,\bar s_0,s_1) =$ \nl
$\dsize \sum_{m \in \{0,1\}} \, 
\dsize \sum_{(s_0,s_1) \in \text{ Pr}_\ell({\Cal D})}
r^*_{(s_0,s_1),m} \bar t\langle x,y,m,s_0,s_1 \rangle$
(we have used \scite{2.3a}).  
\ermn
Now concentrate again on the case $a = \ell \in \{0,1\}$, so
$r^*_{\bar s,1-\ell} = 0$ by $\boxtimes_1$ and $r^*_{\bar s,\ell} = r^*_s$.
So clearly
\mr
\item "{$\circledast_1$}"  $r^*_{\bar s} \ge 0$ \nl
[Why?  As the sum of non negative integers]
\sn
\item "{$\circledast_2$}"   $1 =
\Sigma\{r^*_{\bar s}:\bar s \in \text{ Prd}_\ell({\Cal D})\}$ \nl
[Why?  As by Definition \scite{2.2}(2), $c \in {\Cal D} \and r_c > 0
\Rightarrow c \in {\Cal D}_{x,y} \Rightarrow \bar s^{c,x,y} \in
\text{prd}_\ell({\Cal D})$ and the definition of $r^*_{\bar s}$]
\sn
\item "{$\circledast_3$}"  we have
{\roster
\itemitem{ $(\alpha)$ }  $z \in X \backslash \{x,y\} \Rightarrow s^*_1
= t^*_{y,z} = \Sigma\{r^*_{(s_0,s_1)} \times s_1:(s_0,s_1) \in \text{
Prd}_\ell({\Cal D})\}$,
\sn
\itemitem{ $(\beta)$ }  $z \in X \backslash \{y,z\} \Rightarrow s^*_0
= t^*_{x,z} = \Sigma\{r^*_{(s_0,s_1)} \times s_0:(s_0,s_1) \in \text{ Prd}_\ell
({\Cal D})\}$ \nl
[Why?  By $\boxtimes_2$.]
\endroster}
\ermn
So $\langle r^*_{\bar s}:\bar s \in \text{ Prd}_\ell({\Cal D}) \rangle$
witness that $(s^*_0,s^*_1) \in$ convex hull of Prd$_\ell({\Cal
D})$. \nl
2) Similar proof (and not used). \nl
3) One direction is as in \scite{2.3}(3).  For the other, by the
hypothesis, $\circledast_0$ in the proof of part (1) with $\ell$
replaced by $a$ holds.  So by the part of the proof of part (1) after
$\boxtimes_1,r^*_{\bar s,m}$ are defined and $\boxtimes_2$ holds.  So
\mr
\item "{$\boxdot_1$}"  $r^*_{\bar s,m} \ge 0$
\sn
\item "{$\boxdot_2$}"  $1= \Sigma\{r^*_{\bar s,m}:m \in \{0,1\}$ and
$\bar s \in$ Prd$_m({\Cal D})\}$
\sn
\item "{$\boxdot_3$}"  $(s^*_0,s^*_1) = \Sigma\{r^*_{\bar s,m} \times
\bar s:m \in \{0,1\},\bar s \in$ Prd$_m({\Cal D})\}$ \nl
[Why?  By $\boxtimes_2$.]
\sn
\item "{$\boxdot_4$}"  $a = \Sigma\{r^*_{\bar s,1}:\bar s \in$
Prd$_1({\Cal D})\}$.
\ermn
So we are done.   \hfill$\square_{\scite{2.4}}$\margincite{2.4}
\enddemo
\bn
\ub{Continuation of the proof of \scite{2.1}}:
\sn
\ub{$(d) \Leftrightarrow (b)_{x,y}$}:

Read the definition \scite{2.2}(1) 
(and the symmetry).  \hfill$\square_{\scite{2.1}}$\margincite{2.1}
\bn
\ub{$(d) \Rightarrow (e)$}:

By \scite{2.3}(1).
\bn
\ub{$(e) \Rightarrow (d)$}:

Why?  If clause (e) holds, for some $a \in [0,1]_{\Bbb R} \backslash
\{\frac 12\}$ and $x \ne y \in X$ we have $\bar t^* =: \bar t \langle
x,y,a,\frac 12,\frac 12 \rangle \in$  pr-cl$({\frak D})$.  
If $a > \frac 12$ this witness $(\frac 12,\frac 12) \in \text{ Pr}_{>
1/2}({\Cal D})$, so assume $a < \frac 12$.  But trivially $\bar t
\langle y,x,1-a,\frac 12,\frac 12 \rangle$ is equal to $t^*$ hence (as
in \scite{1.11}) is in pr-cl$({\frak D})$ and we are done.
\bn
\ub{$(e) \Leftrightarrow (f)$}:

Clearly (e) means that
\mr
\item "{$(*)_0$}"  there are $r_c \in [0,1]_{\Bbb R}$ for 
$c \in {\frak D}$ such that 
$1 = \Sigma\{r_c:c \in {\frak D}\}$ and
$a \in [0,1]_{\Bbb R} \backslash \{\frac
12\}$ such that $\bar t \langle x,y,a, \frac 12,\frac 12 \rangle = 
\Sigma\{r_c \times \bar t[c]:c \in {\frak D}\}$. 
\ermn
By \scite{2.4}(3) we know that $(*)_0$ is equivalent to
\mr
\item "{$(*)_1$}"  there are $r_{\bar s,\ell} \in [0,1]_{\Bbb R}$ for $\bar
s \in \text{ Prd}_\ell({\frak D}),\ell \in \{0,1\}$ such that 
{\roster
\itemitem{ $(i)$ }  $1 = \Sigma\{r_{\bar s,\ell}:\bar s \in \text{
Prd}_\ell({\frak D}),\ell \in \{0,1\}\}$
\sn
\itemitem{ $(ii)$ }  $(\frac 12,\frac 12) = 
\Sigma\{r_{\bar s,\ell} \times \bar s:\bar s \in
\text{ Prd}_\ell({\frak D}),\ell \in \{0,1\}\}$
\sn
\itemitem{ $(iii)$ }  $\frac 12 \ne a = \Sigma\{r_{\bar s,1}:
\bar s \in \text{ Prd}_1({\frak D})\}$.
\endroster}
\ermn
But by \scite{2.3}(4) and Definition \scite{2.2}(4) for $\ell \in \{0,1\}$:

$$
\text{Prd}_\ell({\frak D}) = \bigl\{ 
\bigl( \frac{\text{val}_c(x) - \ell}{\bold n -2}, \,
\frac{\text{val}_c(y)-(1-\ell)}{\bold n-2} \bigr): c \in {\frak D} \text{ and }
x \ne y \text{ and } (c\{x,y\} =y) \equiv (\ell =1)\bigr\}.
$$
\mn
Let $d^* \in {\frak D}$ recall that ${\frak D} = 
\text{ sym-cl}(\{d^*\})$ by a hypothesis of \scite{2.1} and 
$V_\ell(d^*) = \{(k_1,k_2):\text{ for
some } x_1 \ne x_2 \in X,k_1 = \text{ val}_d(x_1),k_2 = \text{
val}_d(x_2)$ and $d\{x_1,x_2\} = x_{\ell +1}\}$.  So $(*)_1$ means
(recalling the definition of Prd$_\ell({\frak D}))$
\mr
\item "{$(*)_2$}"  there is a sequence $\langle r_{\bar k,\ell}:\bar k
\in V_\ell(d^*)$ and $\ell \in \{0,1\} \rangle$ such that
{\roster
\itemitem{ $(i)$ }  $r_{\bar k,\ell} \in [0,1]_{\Bbb R}$ and
\sn
\itemitem{ $(ii)$ }  $1 = \Sigma\{r_{\bar k,\ell}:\bar k \in
V_\ell(d^*)\}$ and
\sn
\itemitem{ $(iii)$ }  $(\frac 12,\frac 12) = \Sigma \bigl\{
r_{(k_1,k_2),\ell} \times 
\bigl( \frac{k_1 - \ell}{\bold n-2}, \, 
\frac{k_2 -(1-\ell)}{\bold n-2} \bigr):\ell \in
\{0,1\}$ and \nl

\hskip65pt $(k_1,k_2) \in V_\ell(d^*)\}$
\sn
\itemitem{ $(iv)$ }  $\frac 12 \ne \Sigma\{r_{\bar k,1}:\bar k \in
V_1(d)\}$.
\endroster}
\ermn
Let us analyze $(*)_2$.  
Let $r^*_\ell = \Sigma\{r_{\bar k,\ell}:\bar k \in V_\ell(d^*)\}$ for
$\ell \in \{0,1\}$.  So
\mr
\item "{$\circledast_1$}"  $r^*_\ell \in 
[0,1]_{\Bbb R}$ and $1 = r^*_0 + r^*_1$.
\ermn
If $r^*_\ell = 0$ for some $\ell \in \{0,1\}$ 
we can finish easily.   So we assume
\mr
\item "{$\circledast_2$}"  $r^*_0,r^*_1 \ne 0$.
\ermn
Now clause $(iii)$ of $(*)_2$  means $(iii)_1 + (iii)_2$ where
\mr
\item "{$(iii)_1$}"  $\frac 12 = \frac{1}{\bold n-2} \bigl(
\Sigma\{r_{\bar k,\ell} \times k_1:\bar k \in V_\ell(d^*)$ and $\ell \in
\{0,1\}\} \bigr) -$ \nl

\hskip30pt $\frac{1}{\bold n-2} \, \Sigma\{r_{\bar k,\ell} \times \ell:
\bar k \in V_\ell(d^*),\ell \in \{0,1\}\} = $ \nl

\hskip30pt $\frac{1}{\bold n -2} \, \Sigma\{r_{\bar k,\ell} \times k_1:
\bar k \in V_\ell(d^*)$ and $\ell \in \{0,1\}\} - \frac{r^*_1}{\bold n-2}$, \nl
i.e., 
\sn
\item "{$(iii)'_1$}"  $\frac{\bold n}{2} - 
(1 - r^*_1) = \frac{\bold n -2}{2} + r^*_1 = 
\Sigma\{r_{\bar k,\ell} \times k_1:\bar k \in V_\ell(d^*)$ and $\ell \in
\{0,1\}\}$
\sn
\item "{$(iii)_2$}"  $\frac 12 = \frac{1}{\bold n -2} \,
\Sigma\{r_{\bar k,\ell} \times k_2:\bar k \in V_\ell(d^*)$ and $\ell \in
\{0,1\}\} - \frac{1}{\bold n -2} \, \Sigma\{r_{\bar k,\ell}
\times(1-\ell):$ \nl

\hskip30pt $\bar k \in V_\ell(d^*)$ and $\ell \in \{0,1\}\} =$ \nl

\hskip30pt $\frac{1}{\bold n -2} \, 
\Sigma\{r_{\bar k,\ell} \times k_2:\bar k \in
V_\ell(d^*),\ell \in \{0,1\}\} - \frac{r^*_0}{\bold n -2}$, \nl
i.e., 
\sn
\item "{$(iii)'_2$}"   $\frac{\bold n}{2} -(1 - r^*_0) = \frac{\bold
n}{2} + r^*_1 = 
\Sigma\{r_{\bar k,\ell} \times k_2:\bar k \in V_\ell(d^*)$ and $\ell \in
\{0,1\}\}$.
\ermn
Together $(iii)$ of $(*)_2$ is equivalent to 
\mr
\item "{$(iii)^+$}"  $\bigl(\frac{\bold n}{2} - (1-r^*_1),\frac{\bold
n}{2} - (1-r^*_0) \bigr) =
\Sigma\{r_{\bar k,\ell} \times \bar k:\bar k \in V_\ell(d^*)$ and $\ell \in
\{0,1\}\}$.
\ermn
Let $\bar s_\ell = \Sigma\{r_{\bar k,\ell} \times \bar k:
\bar k \in V_\ell(d^*)\} / r^*_\ell$, so $(*)_2$ is equivalent to
($V_\ell(d^*)$ is from Definition \scite{1.10})
\mr
\item "{$(*)_3$}"  there are $\bar s_0,\bar s_1,r^*_0,r^*_1$ such that
{\roster
\itemitem{ $(i)$ }  $\bar s_\ell \in \text{ convex}(V_\ell(d^*))$ for
$\ell = 0,1$
\sn
\itemitem{ $(ii)$ }  $r^*_0,r^*_1 \in [0,1]_{\Bbb R}$ and $r = r^*_0 +
r^*_1$
\sn
\itemitem{ $(iii)$ }  $(\frac{\bold n}{2} - (1-r^*_1),\frac{\bold n}{2} -
(1-r^*_0))$ is $r^*_0 \times \bar s_0 + r^*_1 \times \bar s_1$
\sn
\itemitem{ $(iv)$ }  $r^*_\ell \ne \frac 12$ (by clause $(iv)$ in $(*)_2$
above).
\endroster}
\ermn
Clearly $(*)_3(iii)$ is equivalent to
\mr
\item "{$(*)_4(iii)'$}"  $(\frac{\bold n}{2} - 1,\frac{\bold n}{2} -
1)$ is $\bigl(r^*_0 \times (\bar s_0 - (0,1)),r^*_1 \times 
(\bar s_1 - (1,0)\bigr)$.
\ermn
Recalling the definition of $V^*_\ell(d^*)$, see Definition
\scite{1.10}(4), clearly $(*)_4(iii)'$ means
\mr
\item "{$(*)_4(iii)''$}"  $(\frac{\bold n}{2} -1,\frac{\bold n}{2} -1)
\in V^*_\ell(d^*)$.
\ermn
But this is clause (f).
\sn
\ub{$(g) \Rightarrow (f)$}:  By \scite{2.5}, \scite{2.6}, \scite{2.7}
below (i.e., they show $(g) + \neg(f)$ lead to contradiction).
\bn
\ub{$(f) \Rightarrow (g)$}:  As $\neg(g) \Rightarrow \neg(f)$
trivially.

We have proved $(b)_{x,y} \Leftrightarrow (b)_{x',y'},(b)_{x,y} \Rightarrow
(a) \Rightarrow (b)_{x,y},(c) \Rightarrow (a) \Rightarrow (c),(d)
\Leftrightarrow (b)_{x,y} (d) \Rightarrow (e) \Rightarrow (d),(e)
\Leftrightarrow (f),(g) \Rightarrow (f) \Rightarrow (g)$, 
so we are done proving \scite{2.1}.  \hfill$\square_{\scite{2.1}}$\margincite{2.1} 
\bigskip

\proclaim{\stag{2.5} Claim}  Assume that clause (f) of \scite{2.1}
fails, $d = d^* \in {\frak D}$ but clause (g) of \scite{2.1} holds
(equivalently $\langle {\text{\rm val\/}}_d(x):x \in X \rangle$ is not
constant).  \ub{Then} the following hold
\mr
\item "{$\boxdot_1$}"  there are no $\bar s^0 \in {\text{\rm
conv\/}}(V^*_0(d^*)),\bar s^1 \in {\text{\rm conv\/}}(V^*_1(d^*))$
such that $(\frac{\bold n}{2} -1,\frac{\bold n}{2} -1)$ lie on
${\text{\rm conv\/}}\{\bar s^0,\bar s^1\}$ and for some $\ell \in
\{0,1\}$ this set contains an interior point of $V^*_\ell(d^*)$
\sn
\item "{$\boxdot_2$}"   the lines $L^*_0 = 
\{(\frac{\bold n}{2} -1,y):y \in \Bbb R\},L^*_1 = 
\{(x,\frac{\bold n}{2} -1):x \in \Bbb R\}$ divides the plane; and
${\text{\rm conv\/}}(V^*(d^*))$ is
{\roster
\itemitem{ $(i)$ }  included in one of the four closed half planes \ub{or}
\sn
\itemitem{ $(ii)$ }  is disjoint to at least one of the closed
quarters minus $\{(\frac{\bold n}{2} -1,\frac{\bold n}{2} -1)\}$. 
\endroster}
\ermn
\endproclaim
\bigskip

\remark{\stag{2.5A} Remark}  1) Recall $V^*(d^*) = 
V^*_0(d^*) \cup V^*_1(d^*)$ and $V^*_\ell(d^*) = 
\{\bar s - (\ell,1-\ell):\bar s \in V_1(d^*)\}$. \nl
2) So
\mr
\widestnumber\item{$(iii)$}
\item "{$(i)$}"  $(k_1,k_2) \in V^*_0(d^*) \Leftrightarrow (k_1,k_2) +
(0,1) \in V_0(d^*) \Leftrightarrow (k_1,k_2+1) \in
V_0(d^*)$
\nl
[Why?  By Definition \scite{1.10}(4).]
\sn
\item "{$(ii)$}"  $(k_2,k_1) \in V^*_1(d^*) 
\Leftrightarrow (k_2,k_1) + (1,0) \in V_1(d^*) 
\Leftrightarrow (k_2+1,k_1) \in V_1(d^*)$ \nl
(see \scite{1.11}(2)) hence
\sn
\item "{$(iii)$}"   $(k_1,k_2) \in V^*_0(d^*) \Leftrightarrow
(k_1,k_2-1) \in V_0(d^*) \Leftrightarrow (k_2-1,k_1) \in V_1(d^*)
\Leftrightarrow (k_2,k_1) \in V^*_1(d^*)$. \nl
[Why? By the above $(i) + (ii)$ and \scite{1.11}(2).]
\endroster
\endremark
\bigskip

\demo{Proof}  Toward contradiction assume that $\boxdot_2$ or
$\boxdot_1$ in the claim fails.  So necessarily
\mr
\item "{$(*)_0$}"  $(\frac{\bold n}{2}-1,\frac{\bold n}{2} -1) \notin
V^*(d^*)$ \nl
[Why?  If it belongs to  $V^*_\ell(d^*)$ let $r^*_\ell =1,r^*_{1 - \ell}
= 0$ and we get clause (f) of \scite{2.1} which we are assuming fails]
\sn
\item "{$(*)'_0$}"  $(\frac{\bold n}{2}-1,\frac{\bold n}{2} -1) \notin
\text{ conv}(V^*_\ell(d^*))$ \nl
[Why?  As in the proof of $(*)_0$.]
\sn
\item "{$(*)_1$}"  $(\frac{\bold n}{2}-1,\frac{\bold n}{2}-1)$ belongs
to the convex hull of $V^*_0(d^*) \cup V^*_1(d^*)$ hence of
conv$(V^*_0(d^*)) \cup \text{ conv}(V^*_1(d^*))$ \nl
[Why?  Otherwise $\boxdot_1$ trivially holds; also there is a line $L$ through
$(\frac{\bold n}{2}-1,\frac{\bold n}{2}-1)$ such that $V^*(d^*)
\backslash L$ lie in one half plane of $L$, so easily clause (ii) of
$\boxdot_2$ holds so $\boxdot_2$ holds.]
\ermn
Let $E = \{(\bar s_0,\bar s_1):\bar s_\ell \in \text{
conv}(V^*_\ell(d^*))$ and $(\frac{\bold n}{2}-1,\frac{\bold n}{2}-1)$
belongs to the convex hull of $\{\bar s_0,\bar s_1\}\}$
\mr
\item "{$(*)_2$}"  $E \ne \emptyset$ \nl
[Why?  By $(*)_1$]
\sn
\item "{$(*)_3$}"  if $r_0,r_1 \in [0,1]_{\Bbb R},1 = r_0 + r_1,
(\frac{\bold n}{2}-1,\frac{\bold n}{2} -1) = r_0 \times \bar s_0 + r_1
\times \bar s_1$ and $\bar s_\ell \in \text{ conv}(V^*_\ell(d^*))$ for
$\ell =0,1$  \ub{then} $r_0 = r_1
= \frac 12$ \nl
[Why?  Otherwise clause (f) holds.]
\sn
\item "{$(*)_4$}"  if $(\bar s_0,\bar s_1) \in E$ 
then $(\frac{\bold n}{2}-1,\frac{\bold n}{2}-1) = \frac 12(\bar s_0 +
\bar s_1)$ \nl
[Why?  By $(*)_3$ and the definition of $E$)
\sn
\item "{$(*)_5$}"  if $(\bar s_0,\bar s_1) \in E,i \in \{0,1\}$,
\ub{then} $\bar s_\ell$ is the unique member of conv$(V^*_\ell(d^*))$
which lies on the line through $\{\bar s_0,\bar s_1\}$ \nl
[Why?  Otherwise let $\bar s'_\ell$ be a counterexample.  If
$(\frac{\bold n}{2}-1,\frac{\bold n}{2}-1) \in \text{ conv}\{\bar
s_\ell,\bar s'_\ell\}$ then it belongs to conv$(V^*_\ell(d^*))$ 
hence clause (f) holds, contradiction.  So letting
$\bar s'_{1-\ell} = \bar s_{1-\ell}$ we know that $(\bar s'_0,\bar
s'_1) \in E$ by $(*)_4$ hence $\frac 12(\bar s'_0,\bar s'_1) =
(\frac{\bold n}{2}-1,\frac{\bold n}{2}-1) = \frac 12(\bar s_0 + \bar
s_1)$ by $(*)_4$ hence subtracting the two equations, $\bar s_{1 - \ell}$ is
cancelled and we get $\bar s'_\ell = \bar s_\ell$, contradiction]
\sn
\item "{$(*)_6$}"  $\boxdot_1$ holds (so by the assumption towards
contradiction $\boxdot_2$ fails). \nl
[Why?  Assume $\bar s_0,\bar s_1$ are as there hence (by the
definition of $E$), $(\bar s_0,s_1) \in E$, now by $(*)_5$ the set (the
line through $\bar s_0,\bar s_1$) $\cap \text{ conv}(V^*_\ell(d))$ is
equal to $\{\bar s_\ell\}$.  So the conclusion of $\boxdot_1$ for
$\bar s_0,\bar s_1$ holds as
$V^*_\ell(d^*)$ is convex for $\ell \in \{0,1\}$.]
\ermn
Easily
\mr
\item "{$(*)_7$}"   $E_\ell = 
\{\bar s_\ell:(\bar s_0,\bar s_1) \in E\}$ is a convex subset
of conv$(V^*_1(d^*)) \subseteq \Bbb R^2$.
\ermn
Also
\mr
\item "{$(*)_8$}"  $(\frac{\bold n}{2} -1,\frac{\bold n}{2}) \notin
E_\ell$ \nl
[Why?  By $(*)_5 + (*)'_0$.]
\endroster
\enddemo
\bn
Now we split the proof to three cases which trivially exhausts all the
possibilities. \nl
\ub{Case 1}:  $E$ is not a singleton. \nl
This implies by $(*)_4$  that $V^*_\ell(d^*)$ is not a singleton for
$\ell=0,1$.  As $|E| \ge 2$ by $(*)_4$
clearly $|E_1| \ge 2$.  Also by $(*)_5$ if 
$\bar s_1 \in E_1$ so $\bar s_1 \ne (\frac{\bold n}{2} -1,\frac{\bold
n}{2} -1)$ by $(*)_8$, then $\bar s_1$ is
the unique member of conv$(V^*_1(d_1)) \cap$ (the line through $\bar
s_1, (\frac{\bold n}{2}-1,\frac{\bold n}{2}-1))$.  Also $E_1$ is
convex (by $(*)_7$) so necessarily $E_1$ lies on a
line $L_1$ to which by $(*)_5$, the point
$(\frac{\bold n}{2}-1,\frac{\bold n}{2}-1)$ does not
belong. \nl
As $E_1 \subseteq V^*_1(d^*) \cap L_1$ is a convex set with $\ge 2$
members and $(*)_5$ it follows that conv$(V^*_1(d^*))$ is included 
in this line $L_1$ and as $V^*_0(d^*) = \{(k_2,k_1):(k_1,k_2) \in
V_1(d^*)\}$ (by clause $(iii)$ of \scite{2.5A}(2)) it follows that
conv$(V^*_0(d^*))$ is included in the line $L_0 = \{(a_0,a_1):(a_1,a_0) \in
L_1\}$ to which $(\frac{\bold n}{2} -1,\frac{\bold n}{2} -1)$ does not belong.

But $E_0 = \{\bar s_0:(\bar s_0,\bar s_1) \in E\}$ is necessarily an
interval of $L_0$ and by $(*)_4$ we have
\mr
\item "{$(*)_5$}"  $L_0 = \{(a_0,a_1):
2(\frac{\bold n}{2}-1,\frac{\bold n}{2} -1) - (a_0,a_1) \in L_1\}$.
\ermn
As $L_1$ is a line, for some reals $r_0,r_1,r_2$ we have

$$
L_1 = \{(a_0,a_1) \in \Bbb R^2:r_0 a_0 + r_1 a_1+r_2 = 0\}
$$
\mn
and

$$ 
(r_0,r_1) \ne (0,0).
$$
\mn
Hence by the definition of $L_0$ we have

$$
L_0 = \{(a_0,a_1) \in \Bbb R^2:r_1a_0+r_0 a_1+r_2 = 0\}
$$
\mn
and by $(*)_4$ the line $L_0$ includes the interval $\{2(\frac{\bold
n}{2} -1,\frac{\bold n}{2} -1) - \bar s_1:\bar s_1 \in E_1\}$ so

$$
L_0 = \{(a_0,a_1):(-r_0)a_0+(-r_1)a_1+r'_2 = 0\}
$$
\mn
where $r'_2 = 
2r_0(\frac{\bold n}{2}-1) + 2r_1(\frac{\bold n}{2} -1) + r_2$.

So for some $s \in \Bbb R$ we have $r_0 = -sr_1,r_1 = -sr_0,r_2 =
sr'_2$ but $(r_0,r_1) \ne (0,0)$
hence $s \in \{1,-1\}$ hence $r_0 \in \{r_1,-r_1\}$, so \wilog
\, $r_0=1,r_1 \in \{1,-1\}$.
\bn
\ub{Subcase 1A}:  $r_1 = -1$.

So $d^*\{x,y\} = y \Rightarrow (\text{val}_{d^*}(x),\text{val}_{d^*}(y)) \in
V_1(d^*) \Rightarrow (\text{val}_{d^*}(x)-1,
\text{ val}_{d^*}(x)) \in L_1 \Rightarrow \text{ val}_{d^*}(x) 
- \text{ val}_{d^*}(y) = -r_2 +1$, i.e., 
is constant, is the same for any such pair $(x,y)$.  
But the directed graph $(X,\{(u,v):c\{u,v\} =
v\})$ necessarily contains a cycle, so necessarily $r_2+1=0$, 
so the val$_{d^*}(x)$ is the same for all $x \in X$ hence is
necessarily $(\frac{\bold n}{2}-1)$, which is an ``allowable" case, in
particular, contradict clause (g) of \scite{2.1} which we are assuming.
\bn
\ub{Subcase 1B}:  $r_1 = 1$.

Clearly $x \ne y \in X \and d^*\{x,y\} = y \Rightarrow
(\text{val}_{d^*}(x),\text{val}_{d^*}(y)) \in V_1(d^*) \Rightarrow
(\text{val}_{d^*}(x)-1,\text{val}_{d^*}(y)) \in V_1(d^*) \Rightarrow 
(\text{val}_{d^*}(x)-1,\text{val}_{d^*}(y)) \in V^*_1(d^*) \Rightarrow
(\text{val}_{d^*}(x)-1,\text{val}_{d^*}(y)) \in L_1 \Rightarrow 
\text{ val}_{d^*}(x) + \text{ val}_{d^*}(x) = -r_2-1$.  As $\bold n \ge 3$
there are distinct $x_0,x_1,x_2 \in X$ so val$_{d^*}
(x_{\ell_1}) + \text{ val}_{d^*}(x_{\ell_2}) = -r_2-1$ for 
$\{\ell_1,\ell_2\} \in \{\{0,1\},\{0,2\},\{1,2\}\}$, the order is not
important as $r_1=r_2$ hence are equal and $-r_2-1$ is twice their
value.  So for $y \in X \backslash \{x_1\}$, we have val$_{d^*}(x_1) +
\text{ val}_{d^*}(y) = -r_2-1$ so val$_{d^*}(y) = 
\text{ val}_{d^*}(x_2)$, so we are done as in case 1A.
\bn
\ub{Case 2}:  $E = \{(\bar s^*_0,\bar s^*_1)\}$ and $\bar s^*_0 \ne
\bar s^*_1$.

Let $L$ be the line through $\{\bar s^*_0,\bar s^*_1\}$ so let the
real $r_0,r_1,r_2$ be such that $L= 
\{(a_0,a_1):r_0a_0+r_1a_1+r_2=0\}$ and $(r_0,r_1) \ne (0,0)$.

So by $(*)_5$ the set 
conv$(V^*_\ell(d^*))$ intersect $L$ in the singleton $\{\bar
s^*_\ell\}$.

If one (closed) half plane for the line $L$ contains $V^*_0(d^*) \cup
V^*_1(d^*)$ then $\boxdot_2$ of \scite{2.5} holds (if $L$ is not
parallel to the $x$-axis and the $y$-axis (i.e., $r_0,r_1 \ne 0$) then
$\boxdot_2(ii)$ holds, otherwise $\boxdot_2(i)$ holds.  So by $(*)_6$ 
we are done.  

If there are $\ell \in \{0,1\}$ and $\bar k',\bar k'' \in
V^*_\ell(d^*) \backslash L$ on different sides of $L$, also
$V^*_{1-\ell}(d^*)$ has a member outside $L$ (by clause $(iii)$ of
\scite{2.5A}(2) and $(*)_5$).
Hence (whether there is such $\ell$ or not, recalling the previous paragraph)
there are $\bar k_0 \in V^*_0(d^*)
\backslash \{\bar s^*_0\}$ and $\bar k_1 \in V^*_1(d^*) \backslash
\{\bar s^*_1\}$ and they are outside $L$ in different sides.  As
$(\frac{\bold n}{2}-1,\frac{\bold n}{2}-1)$ lie in the open interval
spanned by $\bar s^*_0$ and $\bar s^*_1$, necessarily
$(\frac{\bold n}{2}-1,\frac{\bold n}{2}-1)$ is an interior point of
$\{\bar s^*_0,\bar k_0,\bar s^*_1,\bar k_1\}$, contradiction to the
case assumption.
\bn
\ub{Case 3}:  $E = \{(\bar s^*_0,\bar s^*_1)\}$ and $\bar s^*_0 = \bar s^*_1$.

So $\bar s^*_\ell = (\frac{\bold n}{2} -1,\frac{\bold n}{2}-1)$, but
this contradicts $(*)'_0$ above. \hfill$\square_{\scite{2.5}}$\margincite{2.5}
\bigskip

\proclaim{\stag{2.6} Claim}  In \scite{2.5}, clause (i) of $\boxdot_2$ 
is impossible.
\endproclaim
\bigskip

\demo{Proof}  We know that $\langle \text{val}_{d^*}(x):
x \in X \rangle$ is not constant (as we assume clause (g) of
\scite{2.1} holds).  As the average of val$_{d^*}(x),x \in X$ is
$\frac{\bold n-1}{2} = \frac{\bold n}{2}-1$ clearly
\mr
\item "{$(*)_1$}"  for some points $x \in X$ we have val$_{d^*}(x)$ is 
$< \frac{\bold n-1}{2} = \frac{\bold n}{2} - \frac 12$ and for some it is 
$> \frac{\bold n}{2} - \frac 12$.
\ermn
The assumption leaves us four possibilities: which half plan and what
side, so we have 4 cases.  \hfill$\square_{\scite{2.6}}$\margincite{2.6}
\enddemo
\bn
\ub{Case 1}:  For no $(k_0,k_1) \in V^*(d^*)$ do we have $k_0 >
\frac{\bold n}{2}-1$.

It follows that

$$
(k_0,k_1) \in V_0(d^*) \Rightarrow (k_0,k_1-1) \in V^*_0(d^*)
\subseteq V^*(d^*) \Rightarrow k_0 \le \frac{\bold n}{2}-1.
$$
\mn
So if $x \in X$ and for some $y \in X \backslash \{x\}$ we have
$d^*\{x,y\} = x$ (this means just that, val$_{d^*}(x) < \bold n-1)$ then
val$_{d^*}(x) \le \frac{\bold n}{2}-1$, so
\mr
\item "{$(*)_2$}"  if $x \in X$ and val$_{d^*}(x) < \bold n-1$ \ub{then}
val$_{d^*}(x) \le \frac{\bold n}{2}-1$.
\ermn
There can be at most one $x \in X$ with val$_{d^*}(x) = \bold n -1$; if
there is none then we have:

if $x \in X$ then val$_{d^*}(x) \le \frac{\bold n}{2} - 1$.

But the average valency is $\frac{\bold n -1}{2}$ which is $>
\frac{\bold n}{2} -1$, contradiction.  So 
there is $x^* \in X$ such that val$_{d^*}(x) = \bold n-1$,
of course, it is unique.  Now $(X \backslash \{x^*\},\{(y,z):y \ne z
\in X \backslash \{x^*\},d^*\{y,z\} = z\})$ is a tournament with
$\bold n-1$ points each has out-valency $\le \frac{\bold n}{2}-1 = 
\frac{(\bold n-1)-1}{2}$, so necessarily $\bold n$ is even and every
node in the tournament has out-valency exactly $\frac{\bold n}{2} - 1 = 
\frac{(\bold n -1)-1}{2}$.
So $\bold n$ is even and as $\bold n \ge 3$ we can choose $y \ne z \in
X \backslash \{x^*\}$ and \wilog \, $d^*(y,z) = z$.  Now
$(\frac{\bold n}{2}-1,\frac{\bold
n}{2}-1),(\bold n-1,\frac{\bold n}{2}-1) \in V_1(d^*)$ as witnessed by
the pairs $(y,z),(x^*,y)$ respectively, hence
$(\frac{\bold n}{2}-2,\frac{\bold n}{2}-1),(\bold n-2,\frac{\bold
n}{2}-1) \in V^*_1(d)$, again by Definition \scite{1.11}, 
hence (as $\bold n \ge 3$ so $\bold n -2 \ge \frac{\bold n}{2} -1$) \,
$(\frac{\bold n}{2} -1,\frac{\bold n}{2}- 1) \in 
\text{ conv}(V^*_1(d^*))$, contradiction to ``clause (f) of
\scite{2.1} fails" assumed in \scite{2.5}.
\bn
\ub{Case 2}:  For no $(k_0,k_1) \in V^*(d^*)$ do we have $k_0 <
\frac{\bold n}{2}-1$.  \nl
So

$$
(k_0,k_1) \in V_0(d^*) \Rightarrow (k_0,k_1-1) \in V^*_0(d^*)
\subseteq V^*(d^*) \Rightarrow k_0 \ge \frac{\bold n}{2}-1.
$$
\mn
So if $x \in X$ and for some $y \in X \backslash \{x\}$ we have
$d^*\{x,y\} = x$ (equivalently val$_{d^*}(x) < \bold n -1)$
then $(\text{val}_{d^*}(x),\text{val}_{d^*}(y)) \in V_0(d^*)
\Rightarrow (\text{val}_{d^*}(x),\text{val}_{d^*}(y)-1) \in V^*_1(d^*)$ hence
val$_{d^*}(x) \ge \frac{\bold n}{2}-1$, but $\bold n -1 \ge
\frac{\bold n}{2} -1$ so in any case
\mr
\item "{$(*)_2$}"  if $x \in X$ then val$_{d^*}(x) \ge \frac{\bold
n}{2}-1$.
\ermn
So $x \in X \Rightarrow \text{ val}_{d^*}(x) \ge \frac{\bold n}{2}-1$.

If $\bold n$ is odd we have $x \in X \Rightarrow \text{ val}_{d^*}(x)
\ge \frac{\bold n}{2} - \frac 12 = \frac{\bold n -1}{2}$, impossible
by $(*)_1$ so $\bold n$ is even.  Let $k = \frac{\bold n}{2}-1$.  The
average val$_{d^*}(x)$ is necessarily $k + \frac 12$ hence $Y =: \{x \in
X:\text{val}_{d^*}(x) \le k$ (equivalently $=k)\}$ has at least 
$k+1 = \frac{\bold n}{2}$ members.  If $x \in
X$,val$_{d^*}(x) = k+1$ then $x \notin Y,|\{y:d\{x,y\} = y\}| = k+1 =
\frac{\bold n}{2} > |X \backslash Y \backslash \{x\}|$ 
so there is $y \in Y$ such that $d\{x,y\} = y$, hence $(k,k) =
(\text{val}_{d^*}(x)-1$, val$_{d^*}(y)) \in V_1(d^*)$ so clause (f) of
\scite{2.1} holds, contradiction. So
\mr
\item "{$(*)_3$}"  $x \in X \Rightarrow \text{ val}_{d^*}(x) \ne k+1$.
\ermn
Now $|Y| = \bold n$ is impossible by $(*)_1$.  Also if $|Y| = \bold
n-1$ let $z$ be the unique element of $X$ outside $Y$ so in the
tournament $(Y,\{(x,y):x \ne y$ are from $Y$ and $d^*\{x,y\} = y\})$
each $x$ has out-valency $\le$ val$_{d^*}(x) = k = \frac{(|Y|-1)}{2}$,
but this is the average so equality holds and we get contradiction as 
in Subcase 1b.  Hence
\mr
\item "{$(*)_4$}"  $|Y| \le \bold n -2$.
\ermn
Clearly we can find $x_1 \in Y$ such that $|\{y \in Y:y \ne x_1,d\{x,y_1\} =
y_1\}| \le \frac{|Y|-1}{2}$ (considering the average valency inside
$Y$) but $\frac{|Y|-1}{2} \le \frac{\bold n}{2} - \frac 32 = k - \frac
12$ hence there is $x_2 \in X \backslash Y$ such that
$d^*\{x_1,x_2\} = x_2$.  Now
let $m = \text{ val}_{d^*}(x_2)$ so $m > k+1$ and $(k-1,m) \in
V^*_1(d^*)$.  As $\bold n \ge 3,|Y| \ge \frac{\bold n}{2}$
obviously $|Y| \ge 2$ hence also $(k-1,k) \in V^*_1(d^*)$ and as easily
there is $x_3 \in Y$ such that $d\{x_2,x_3\} = x_3$ also $(m-1,k) \in
V^*_1(d^*)$ so $(\frac{\bold n}{2} -1,\frac{\bold n}{2} -1) = (k,k)
\in \text{ conv}\{(k-1,k),(m-1,k)\})$ recalling $m > k+1$, 
contradiction to ``not clause (f) of \scite{2.1}" with $r^*_0=0$.  
\bn
\ub{Case 3}: For no $(k_0,k_1) \in V^*(d^*)$ do we have $k_1 >
\frac{\bold n}{2} -1$.

So

$$
(k_0,k_1) \in V_1(d^*) \Rightarrow (k_0 -1,k_1) \in V^*_1(d^*)
\subseteq V^*(d^*) \Rightarrow k_1 \le \frac{\bold n}{2} -1.
$$
\mn
So if $y \in X$ and for some $x \in X \backslash \{y\}$ we have
$d\{x,y\} = y$ then $(\text{val}_{d^*}(x),\text{val}_{d^*}(y)) \in
V_1(d^*)$ hence val$_{d^*}(y) \le \frac{\bold n}{2} -1$.  But there is
such $x$ iff val$_{d^*}(y) \ne \bold n -1$, that is
\mr
\item "{$(*)_5$}"  if $y \in X$ and val$_{d^*}(y) \ne \bold n -1$ then
val$_{d^*}(y) \le \frac{\bold n}{2} -1$.
\ermn
We continue as in Case 1, (or dualize see \scite{1.11}(1)).
\bn
\ub{Case 4}:  For no $(k_0,k_1) \in V^*(d^*)$ do we have $k_1 <
\frac{\bold n}{2} -1$.  

So $(k_0,k_1) \in V_1(d^*) \Rightarrow (k_0-1,k_1) \in V^*_1(d^*)
\subseteq V^*(d^*) \Rightarrow k_1 \ge \frac{\bold n}{2} -1$.  So if
$y \in X$ and for some $x \in X \backslash \{y\}$ we have $d\{x,y\} =
y$ then $(\text{val}_{d^*}(x),\text{val}_{d^*}(y)) \in V_1(d^*)
\Rightarrow (\text{val}_{d^*}(x)-1,\text{val}_{d^*}(y)) \in V^*_1(d^*)
\subseteq V^*(d^*) \Rightarrow \text{ val}_{d^*}(y) \ge \frac{\bold
n}{2} -1$.  So if $y \in X$ and for some $x \in X \backslash \{y\}$
we have $d^*\{x,y\} = y$ then val$_{d^*}(y) \ge \frac{\bold n}{2}-1$
so if val$_{d^*}(y) < \bold n-1$ then val$_{d^*}(y) \ge \frac{\bold
n}{2} -1$, but if val$_{d^*}(y) \ge \bold n -1$ we get the same
conclusion, so
\mr
\item "{$(*)_6$}"  val$_{d^*}(y) \ge \frac{\bold n}{2} -1$
\ermn
and we can continue as in case 2.  \hfill$\square_{\scite{2.6}}$\margincite{2.6}
\bigskip

\proclaim{\stag{2.7} Claim}  In \scite{2.5}, clause (ii) of $\boxdot_2$
is impossible.
\endproclaim
\bigskip

\demo{Proof}  Note that as we are assuming clause (f) of \scite{2.1}
\mr
\item "{$(*)_0$}"  $(\frac{\bold n}{2} -1,\frac{\bold n}{2} -1) \notin
V^*(d^*)$.
\ermn
Again we have four cases.
\sn
\ub{Case 1}:  If $a_0 \ge \frac{\bold n}{2} -1,a_1 \ge \frac{\bold
n}{2} -1$ but $(a_0,a_1) \ne (\frac{\bold n}{2} -1,\frac{\bold n}{2} -1)$
then $(a_0,a_1) \notin \text{ conv}(V^*(d^*))$. \nl
So
\mr
\item "{$(*)_1$}"  for at most one $x \in X$ we have val$_{d^*}(x) \ge
\frac{\bold n}{2}$.\nl
[Why?  If $x \ne y \in X$ and val$_{d^*}(x) \ge \frac{\bold n}{2}$,
val$_{d^*}(x) \ge \frac{\bold n}{2}$ \ub{then}
$(\text{val}_{d^*}(x)-1,\text{val}_{d^*}(y)) \in V^*_0(d^*) \subseteq
V^*(d^*)$ or $(\text{val}_{d^*}(x),\text{val}_{d^*}(y)-1) \in
V^*_0(d^*) \subseteq V^*(d^*)$, a contradiction to the case assumption
in both cases.]
\ermn
If there is no $x \in X$ with val$_{d^*}(x) \ge \frac{\bold n}{2}$ 
then $x \in X \Rightarrow \text{ val}_{d^*}(x) <
\frac{\bold n}{2}$ and so $x \in X \Rightarrow \text{ val}_{d^*}(x)
\le \frac{\bold n-1}{2}$ but this is the average valency, so always equality
holds, contradicting an assumption of \scite{2.5}.
\nl
So assume
\mr
\item "{$(*)_2$}"  $x_0 \in X$, val$_{d^*}(x_0) \ge \frac{\bold
n}{2}$.
\ermn
Now
\mr
\item "{$(*)_3$}"  $y \in X \backslash \{x_0\} \Rightarrow \text{
val}_{d^*}(y) < \frac{\bold n}{2} -1$ (hence $\le \frac{\bold n}{2} -
\frac 32$) \nl
[Why?  If $d^*\{x_0,y\} = x_0$ then $(\text{val}_{d^*}(x_0)-1,
\text{val}_{d^*} (y)) \in V^*_1(d^*) \subseteq V^*(d^*)$, now
$\text{val}_{d^*}(x_0) -1 \ge \frac{\bold n}{2} -1$ hence by the case
assumption $+ (*)_0$ we have $\text{val}_{d^*}(y) < \frac{\bold n}{2}
-1$.  If $d^*(x_0,y) = y$ then
$(\text{val}_{d^*}(x_0)-1,\text{val}_{d^*}(y)) \in V^*_0(d^*)
\subseteq V^*(d^*)$ and the proof if similar.]
\ermn
So letting $Y = X \backslash \{x_0\}$ we have
$(Y,\{\{y,z\},y \ne z \in Y,d^*(y,z) =
z\})$ is a tournament each node has out-valency $\le \frac{\bold
n-3}{2} < \frac{(\bold n-1)-1}{2} = \frac{|Y|-1}{2}$, contradiction.
\enddemo
\bn
\ub{Case 2}:  If $a_1 \le \frac{\bold n}{2}-1$ and $a_2 \le
\frac{\bold n}{2}-1$ and $(a_1,a_2) \ne (\frac{\bold n}{2}
-1,\frac{\bold n}{2} -1)$ \ub{then} $(a_1,a_2) \notin \text{ conv}(V^*(d^*))$.
\nl
Clearly, as above in the proof of $(*)_3$
\mr
\item "{$(*)'_1$}"  there is at most one $x \in X$ with val$_{d^*}(x)
\le \frac{\bold n}{2}-1$.
\ermn
If there is none then $x \in X \Rightarrow \text{ val}_{d^*}(x) \ge
\frac{\bold n}{2}-1 + \frac 12 = \frac{\bold n-1}{2}$, so considering the
average of val$_{d^*}(y)$ equality always holds so clause (g) of
\scite{2.1} holds contradicting an assumption of \scite{2.5}.  So
assume
\mr
\item "{$(*)'_2$}"  $x_0 \in X$, val$_{d^*}(x_0) \le \frac{\bold n}{2}-1$
\ermn
and we can show, as above that
\mr
\item "{$(*)'_3$}"  if $y \in X \backslash \{x_0\}$
then val$_{d^*}(y) > (\frac{\bold n}{2}-1) +1 = \frac{\bold n}{2}$ so
val$_{d^*}(y) \ge \frac{\bold n-1}{2}$.
\ermn
The directed graph $(X \backslash \{x_0\},\{(y,z\}:d\{y,z\} = z\})$
has $\bold n-1$ nodes, and the out-valency of very node is 
$> \frac{\bold n}{2} -1$ hence is $\ge \frac{\bold n}{2} -1 + \frac 12
= \frac{\bold n}{2} - \frac 12 = \frac{\bold n-1}{2}$, 
(check by cases) easy contradiction as in case (1).
\bn
\ub{Case 3}:  If $a_1 \ge \frac{\bold n}{2}-1$ and $a_2 \le
\frac{\bold n}{2}-1$ and $(a_1,a_2) \ne (\frac{\bold n}{2}
-1,\frac{\bold n}{2} -1)$, \ub{then} $(a_1,a_2) \notin \text{
conv}(V^*(d^*))$.
\nl
So (as in the proof of $(*)_3$)
\mr
\item "{$\bigodot_1$}"  there cannot be $x_0,x_1 \in X$ such that
val$_{d^*}(x_0) \ge \frac{\bold n}{2}$ and val$_{d^*}(x_1) \le
\frac{\bold n}{2}-1$ (the $x_0 \ne x_1$ follows)
\ermn
so one of the following subcases hold.
\bn
\ub{Subcase 3A}:  $x \in X \Rightarrow \text{ val}_{d^*}(x) <
\frac{\bold n}{2}$.

So $x \in X \Rightarrow \text{ val}_{d^*}(x) \le \frac{\bold n-1}{2}$
and (looking at average valency) equality holds, contradicting clause
(g) of \scite{2.1} which we are assuming.
\bn
\ub{Subcase 3B}:  $x \in X \Rightarrow \text{ val}_{d^*}(x) >
\frac{\bold n}{2}-1$.

So $x \in X \Rightarrow \text{ val}_{d^*}(x) \ge \frac{\bold n-1}{2}$,
and we finish as above.
\bn
\ub{Case 4}: If $a_1 \le \frac{\bold n}{2} -1$ and $a_2 
\ge \frac{\bold n}{2}-1$ then $(a_1,a_2) \notin \text{ conv}(V(d^*))$.
\nl
Similar (or dualize the situation by \scite{1.11}(1)).
\hfill$\square_{\scite{2.7}}$\margincite{2.7}
\newpage

\head {\S3 Balanced choice functions} \endhead  \resetall \sectno=3
\bigskip

\definition{\stag{b.1} Definition}  1) $c \in {\frak C}$ is called
\ub{balanced} if $x \in X \Rightarrow \text{ val}_c(x) = (\bold n-1)/2$,
let ${\frak C}^{\text{bl}} = \{c \in {\frak C}:c$ is balanced$\}$. \nl
2) $\bar t \in \text{ pr}({\frak C})$ is called balanced \ub{if} $x
\ne y \in X \Rightarrow \Sigma\{t_{x,y}:y \in X \backslash \{x\}\} =
(\bold n-1)/2$.  Let $\text{pr}^{\text{bl}}({\frak C})$ be the set of balanced
$\bar t \in \text{ pr}({\frak C})$. \nl
3) $c \in {\frak C}$ is called pseudo-balanced \ub{if} for some
balanced $\bar t \in \text{ pr}({\frak C})$ we have

$$
c\{x,y\} = y \Rightarrow t_{x,y} > \frac 12.
$$
\mn
4) We call ${\Cal D} \subseteq {\frak C}$ balanced \ub{if} every $c
\in {\Cal D}$ is balanced, similarly $T \subseteq \text{ pr}({\frak
C})$ is called balanced if every $\bar t \in T$ is. \nl
5) If $x,y,z \in X$ are distinct, let $\bar t = \bar t^{<x,y,z>}$
be defined by: 

$t_{u,v}$ \nl

is $1$ if $(u,v) \in \{(x,y),(y,z),(z,x)\}$. \nl

is $0$ if $(u,v) \in \{(y,x),(z,y),(x,z)\}$ \nl

is $\frac 12$ if otherwise. \nl
6) For a sequence 
$\bar x = (x_0,\dotsc,x_{k-1})$ with $x_\ell \in X$ with no repetitions
and $a \in [0,1]_{\Bbb R}$ let $\bar t = 
\bar t_{\bar x,a} \in \text{ pr}({\frak C})$ be 
defined by $t_{x_i,x_j} = a,t_{x_j,x_i} = 1-a$ if
$j = i+1$ mod $k$ and $t_{x,y} = \frac 12$ otherwise.  If $a=1$ we may
omit it.
\enddefinition
\bigskip

\proclaim{\stag{3.1a} Claim}  If ${\frak D} \subseteq {\frak C}$ is non
empty, symmetric and not balanced, \ub{then} ${\text{\rm
maj-cl\/}}({\frak D}) = {\frak C}$.
\endproclaim
\bigskip

\demo{Proof}  Choose $d^* \in {\frak D}$ which is not balanced, and
let ${\frak D}' =: \text{ sym-cl}(\{d^*\})$, so ${\frak D}'$ is as in
\scite{2.1} and it satisfies clause (g) there hence it satisfies clause
(a) there.  This means that maj-cl$({\frak D}') = {\frak C}$ but
${\frak D}' \subseteq {\frak D} \subseteq {\frak C}$ hence ${\frak C}
= \text{ maj-cl}({\frak D}') \subseteq \text{ maj-cl}({\frak D})
\subseteq {\frak C}$ so we are done.  \hfill$\square_{\scite{3.1a}}$\margincite{3.1a}
\enddemo
\bn
\margintag{b.2}\ub{\stag{b.2} Fact}  1) If $c \in {\frak C}^{\text{bl}}$ then $\bar t[c]$
belongs to $\text{pr}^{\text{bl}}({\frak C})$; if $\bar t \in \text{
pr}({\frak C})$ is balanced and $c \in \text{ maj}\{\bar t\}$,
\ub{then} $c$ is pseudo-balanced; (also if $c$ is pseudo balanced then
it satisfies the condition from \scite{1.9}; and it follows from (4a)). \nl
2)  If ${\Cal D}$ is balanced, \ub{then} pr-cl$({\Cal D})$ is
balanced hence every member of maj-cl$({\Cal D})$ is
pseudo-balanced.
\bn
\margintag{b.2a}\ub{\stag{b.2a} Fact}:  1) If $c \in {\frak C}$ 
is pseudo-balanced, \ub{then} every edge
belongs to some directed cycle. \nl
2) Assume that $\bar t \in \text{ pr}({\frak C})$ is balanced, then
\mr
\item "{$(a)$}" if $t_{u,v} > \frac 12$ then we can find $k \ge 3$ and
$\langle x_0,\dotsc,x_{k-1} \rangle \in X$ with no repetitions such
that $(x_0,x_1) = (u,v)$ and $j=i+1$ mod $k \Rightarrow t_{x_i,x_j} >
\frac 12$
\sn
\item "{$(b)$}"  we can find $m(*)$ and $\bar x_m,a_m$ as in \scite{b.1}(6) 
and $r_m \ge 0$ for $m < m(*)$ satisfying $\Sigma\{r_m:m < m(*)\} =1$ such
that $\bar t^* = \Sigma\{r_m \bar t_{{\bar x}_m,a_m}:m < m(*)\}$ is
quite similar to $\bar t$: for some real $s \ge 1$ we have $u \ne v
\in X \Rightarrow t_{u,v} - \frac 12 = s(t^*_{u,v} - \frac 12)$
\sn
\item "{$(c)$}"  in (b); if all $t_{u,v}$ are rational then 
we can add $r_m = \frac{1}{m(*)},s = m(*)$.
\ermn
3) In fact (a), (b), (c) are equivalent.
\bigskip

\demo{Proof of \scite{b.2}}  1), 2) Check. 
\enddemo
\bigskip

\demo{Proof of \scite{b.2a}}  1) By (a) of part (2). \nl
2) \ub{Clause (a)}: by clause (b).
\mn
\ub{Clause (b)}:  We define a directed graph $G = G(\bar t)$ as
follows: the set of nodes is $X$, the set of edge $E = E^G = 
\{(u,v):u \ne v \in X,t_{u,v} > \frac
12\}$ and we define the function 
$\bold w = w^{\bold t}:E^{G(\bar t)} \rightarrow \Bbb R^{>
0}$ by $\bold w(u,v) = t_{u,v} - \frac 12$.  Now $(G,\bold w)$ or
$(E^G,\bold w)$ is called balanced in the sense that
\mr
\item "{$(*)_{E,\bold w)}$}"  for every node the in-valency is equal
to the out-valency, i.e., 
$\Sigma\{w(u):(u,x) \in E^G\} = \Sigma\{w(x,v):(x,v)
\in E^G\}$. 
\ermn
For such pairs $(E,\bold w)$ we shall prove that for some $(\bar
x_m,a_m)$ as in \scite{b.1}(6) for $m=0,\dotsc,m(*)-1$ we have $\bold
w^{\bold t} = \Sigma\{r_m \bold w^{{\bar t}_{{\bar x}_m,a_m}}:m <
m(*)\}$ for some $r_m \ge 0$ satisfying $1 = \Sigma\{r_m:m <m(*)\}$.
We prove this by induction on $\{(u,v) \in E:\bold w(u,v) > 0\}$.  For
any such pair $(E,\bold w)$ if $E \ne \emptyset$, by the equality 
$(*)_{(E,\bold w)}$ for any node, the in-valency is 
$> 0$ iff the out-valency is $>0$ hence
there is a directed cycle $\bar x = \langle x_0,\dotsc,x_{k-1}
\rangle$; so $a = \text{ Min}\{\bold w(x_0,x_1),\bold
w(x_1,x_2),\dotsc,\bold w(x_{k-2},x_{k-1}),\bold w(x_{k-1},x_0)\}$ is
$>0$.  Define $E' = E \backslash \{(x_i,x_j):j = i+1$ mod $k$ and
$\bold w(x_i,x_j) =a\},w'(u,v)$ is $w(u,v)-a$ if $(u,v) =
(x_i,x_j),j=i+1$ mod $k$ and $\bold w'(u,v) = \bold w(u,v)$ otherwise.
Applying the induction hypothesis to $(E',\bold w')$ and adding $\bar
t_{\bar x,a}$ we are done.  \nl
3) Should be clear.   \hfill$\square_{\scite{b.2}}$\margincite{b.2}
\enddemo
\bigskip

\proclaim{\stag{b.3} Claim}  Assume $|X| \ge 3,{\frak D} \subseteq
{\frak C}$ is symmetric, non empty and balanced.  \ub{Then}, for any
distinct $x,y,z \in Z$ we have $\bar t_{<x,y,z>} \in 
{\text{\rm pr-cl\/}}({\frak D})$. 
\endproclaim
\bigskip

\demo{Proof}  Let $d \in {\frak D}$, now as $(X,\{(u,v):d\{u,v\} =
v\})$ is a directed graph even a tournament with 
equal out-valance and in-valance
for every node, it has directed cycle and hence it has 
a triangle, i.e., $x,y,z \in X$ distinct such that
\mr
\item "{$(*)_1$}"  $d\{x,y\} = y,d\{y,z\} = z,d\{z,x\} = x$.
\ermn
Let $\Pi_{x,y,z} = \{\pi \in \text{ Per}(X):\pi \restriction
\{x,y,z\}$ is the identity$\}$.  Let $\bar t = \Sigma\{\bar
t[c^\pi]:\pi \in \Pi_{x,y,z}\}/|\Pi_{x,y,z}|$. \nl
Clearly $c^\pi \in {\frak D}$ for $\pi \in \Pi_{x,y,z}$ hence $\bar t
\in \text{ pr-cl}({\frak D})$.  Also by $(*)_1$ and the definition
of $\Pi_{x,y,z}$
\mr
\item "{$(*)_2$}"  $t_{x,y} = t_{y,z} = t_{z,y} =1$.
\ermn
Also

$$
\align
|\{w:&w \in X \backslash \{x,y,z\} \text{ and } d\{x,w\} = w\}| = \\
  &|\{w:w \in X \backslash \{x\} \text{ and } d\{x,w\} = w\}| - \{w:w
\in \{y,z\} \text{ and } d\{x,w\} = w\} = \\
  &(|X|-1)/2 -1 = (|X| -3)/2 = \\
  &|\{w:w \in X \backslash \{x,y,z\} \text{ and } d\{x,w\} = x\}|
\endalign
$$
\mn
hence
\mr
\item "{$(*)_3$}"  $t_{x,w} = 1/2$ for $w \in X \backslash \{x,y,z\}$.
\ermn
Similarly
\mr
\item "{$(*)_4$}"  $t_{y,w} = 1/2$ for $w \in X \backslash \{x,y,z\}$
\sn
\item "{$(*)_5$}"  $t_{z,w} = 1/2$ for $w \in X \backslash \{x,y,z\}$
\ermn
and even easier (and as in \S2)
\mr
\item "{$(*)_6$}"  $t_{u,v} = 1/2$ if $u \ne v \in X 
\backslash \{x,y,z\}$.
\ermn
So we are done.  \hfill$\square_{\scite{b.3}}$\margincite{b.3}
\enddemo
\bigskip

\proclaim{\stag{b.4} Claim}  Assume ${\frak D} \subseteq {\frak C}$ is
symmetric non empty and $c^* \in {\frak C}$ is pseudo balanced
\ub{then} $c^* \in {\text{\rm maj-cl\/}}({\frak D})$.
\endproclaim
\bigskip

\demo{Proof}  Without loss of generality ${\frak D}$ is balanced
(otherwise use \scite{3.1a}).  So by \scite{b.3}
\mr
\item "{$\circledast$}"  if $x,y,z \in X$ are distinct then 
$\bar t_{<x,y,z>} \in \text{ pr-cl}({\frak D})$.
\ermn
Let $\bar t^* = \bar t[c^*]$ and let $\langle \bar x^i:i < i(*) \rangle$
list the set cyc$(d)$ of tuples $\bar x = \langle x_\ell:\ell \le k
\rangle$ such that:

$$
k \ge 2,x_\ell \in X,\ell_1 < \ell_2 \le k \Rightarrow x_{\ell_1} \ne
x_{\ell_2}
$$

$$
c\{x_\ell,x_{\ell +1}\} = x_{\ell +1} \text{ for } \ell < k
$$

$$
c\{x_k,x_0\} = x_0
$$
\mn
Note
\mr
\item "{$\otimes$}"  for every $\bar x \in \text{ Cyc}(c^*)$ for some
$\bar t = \bar t^{\bar x} \in \text{ pr-cl}({\frak D})$ we have
{\roster
\itemitem{ $(a)$ }  $t_{u,v} = \frac 12 + \frac{1}{\ell g(\bar x)-2} 
\,\,\text{ \ub{if}}$ \nl

\hskip40pt $(u,v) \in \{(x_{\ell_1},x_{\ell_2}):\ell_1 < \ell g(\bar x)-1 \and
\ell_2 = \ell_1 +1 \text{ or}$ \nl

\hskip40pt $\ell_1 = \ell g(\bar x)-1 \and \ell_2=0\}$
\sn
\itemitem{ $(b)$ }  $t_{u,v} = \frac 12 
- \frac{1}{\ell g(\bar x)-2} \,\, \text{ \ub{if} }
(v,u) \text{ is as above}$
\sn
\itemitem{ $(c)$ }  $t_{u,v} = \frac 12 \text{ \ub{if} otherwise}$. 
\endroster}
\ermn
[Why?  If $\bar x = \langle x_\ell:\ell \le k \rangle$, let $\bar t$
be the arithmetic average of \nl
$\langle t^{<x_0,x_1,x_2>},
\bar t^{<x_0,x_2,x_3>},\dotsc,\bar t^{<x_0,x_{k-1},x_k>}\rangle$.] 
\sn
Now let

$$
\bar t = \Sigma\{ \frac{1}{i(*)} \bar t^{\bar x^i}:i < i(*)\}.
$$
\mn
(In fact we just need that $c\{y_0,y_1\} = y_1 \Rightarrow (y_0,y_1)$
appears in at least one cycle $\bar x^i,i < i(*)$).
\hfill$\square_{\scite{b.4}}$\margincite{b.4} 
\enddemo
\bn
So now we can give a complete answer.
\demo{\stag{b.5} Conclusion}  Assume
\mr
\item "{$(a)$}"  ${\frak D} \subseteq {\frak C}$ is symmetric non
empty
\sn
\item "{$(b)$}"  $c^* \in {\frak D}$.
\ermn
\ub{Then} $c^* \in \text{ maj-cl}({\frak D})$ iff ${\frak D}$ has
a non balanced member or $c^*$ is pseudo balanced.
\enddemo
\bigskip

\demo{Proof}  If ${\frak D}$ has no nonbalanced member and $c^*$ is
not pseudo-balanced, by \scite{b.2}(2) we know $c^* \notin$
maj-cl$({\frak D})$.  For the other direction, if 
$d^* \in {\frak D}$ is not balanced use \scite{b.2a} that is
$(a) \Leftrightarrow (g)$ 
of claim \scite{2.1} for sym-cl$\{d^*\}$.
Otherwise ${\frak D}$ is balanced, $c^*$ is pseudo balanced and 
we use \scite{b.4}.  \hfill$\square_{\scite{b.5}}$\margincite{b.5} 
\enddemo
\newpage

\head {\S4 More flexible voters} \endhead  \resetall \sectno=4
\bigskip

We may now relook and note that maj$(\bar t)$ for $\bar t \in \text{
pr}({\frak C})$ naturally allow a draw, in the case $t_{u,v} = \frac
12$.  We may allow the ``voters" to abstain, i.e., replace ${\frak C}$ by
${\frak C}_p$, the set of partial choice functions (so $\bar t = \bar
t[c]$ will be defined such that $t_{u,v} = \frac 12 \Leftrightarrow
(t\{u,v\}$ undefined)).  We may further ask what can be
$\{\text{maj}(c,\bar t):\bar t \in \text{ conv}(T)\}$ for $T$ a finite
subset of pr$({\frak C})$.  We shall continue this elsewhere. 
\newpage

     \shlhetal 

\newpage
    
REFERENCES.  
\bibliographystyle{lit-plain}
\bibliography{lista,listb,listx,listf,liste}

\enddocument